\documentclass[reqno,12pt,letterpaper]{amsart}

\usepackage{amsmath,amssymb,amsthm,amsxtra}
\usepackage{esint}
\usepackage{mathrsfs}
\usepackage{mathtools}
\usepackage{longtable}

\usepackage{graphicx}
\usepackage{float}
\usepackage[font=small,labelfont=bf]{caption}
\usepackage{subcaption}

\usepackage[dvipsnames]{xcolor}
\usepackage{tikz}
\usepackage{url}
\usepackage[colorlinks=true,linkcolor=Red,citecolor=Green,hypertexnames=false]{hyperref}
\usepackage{cleveref}
\crefname{figure}{Figure}{Figures}
\Crefname{figure}{Figure}{Figures}

\newcommand{\globalregimefigure}{%
\begin{figure}[htbp]
  \centering
  \includegraphics[width=0.9\linewidth]{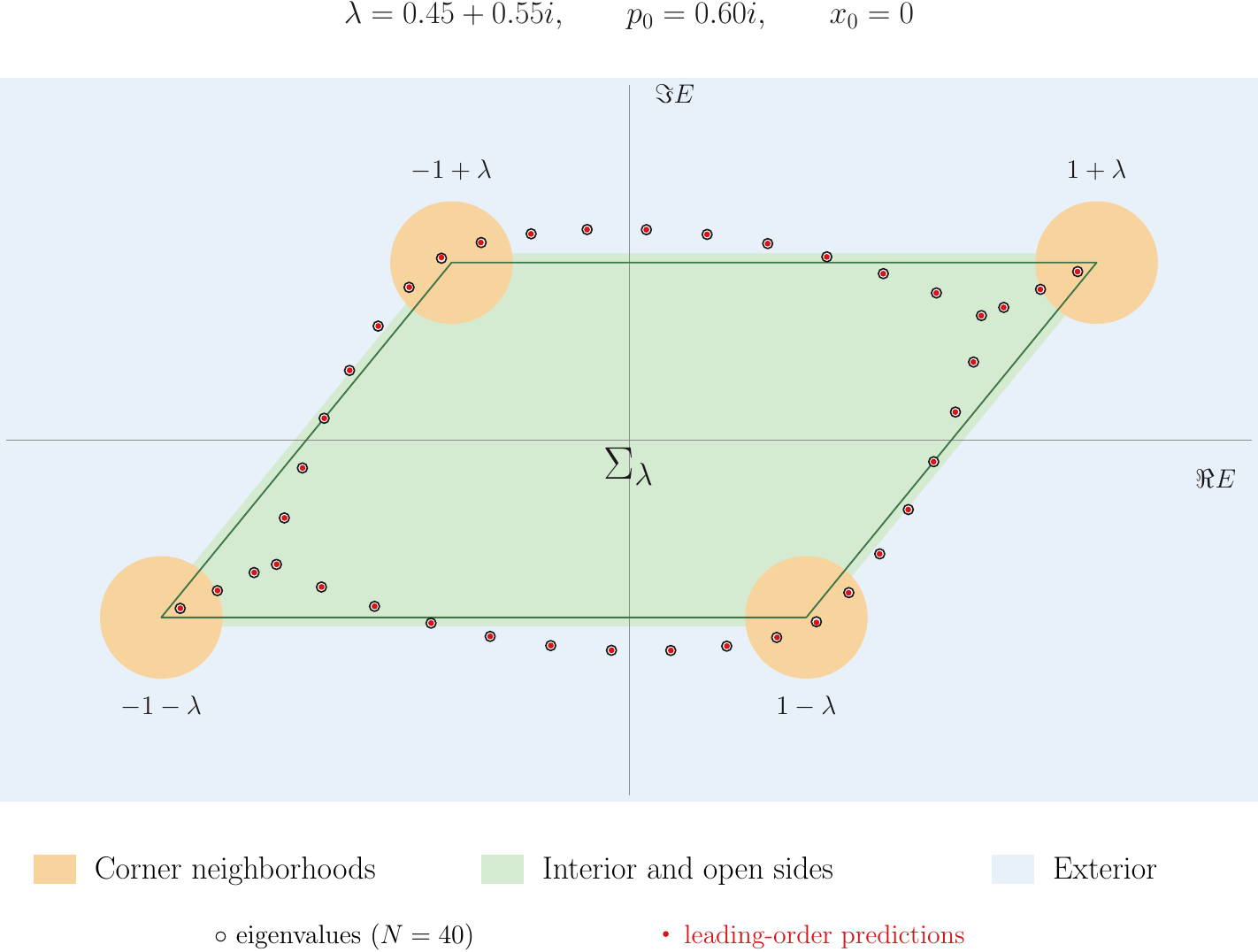}
  \caption{Computed eigenvalues (open circles) and predictions from
  \Cref{thm:main} and \Cref{thm:main-corner} with the correction terms omitted (red dots); the maximum
  absolute difference is approximately $2.34\times10^{-3}$.
  Shading indicates interior and open sides (green), exterior (blue),
  and corner neighborhoods (amber).}
  \label{fig:global-energy-regimes}
\end{figure}
}

\newcommand{\energydiagramfigure}{%
\begin{figure}[!htbp]
  \centering
  \includegraphics[width=\linewidth]{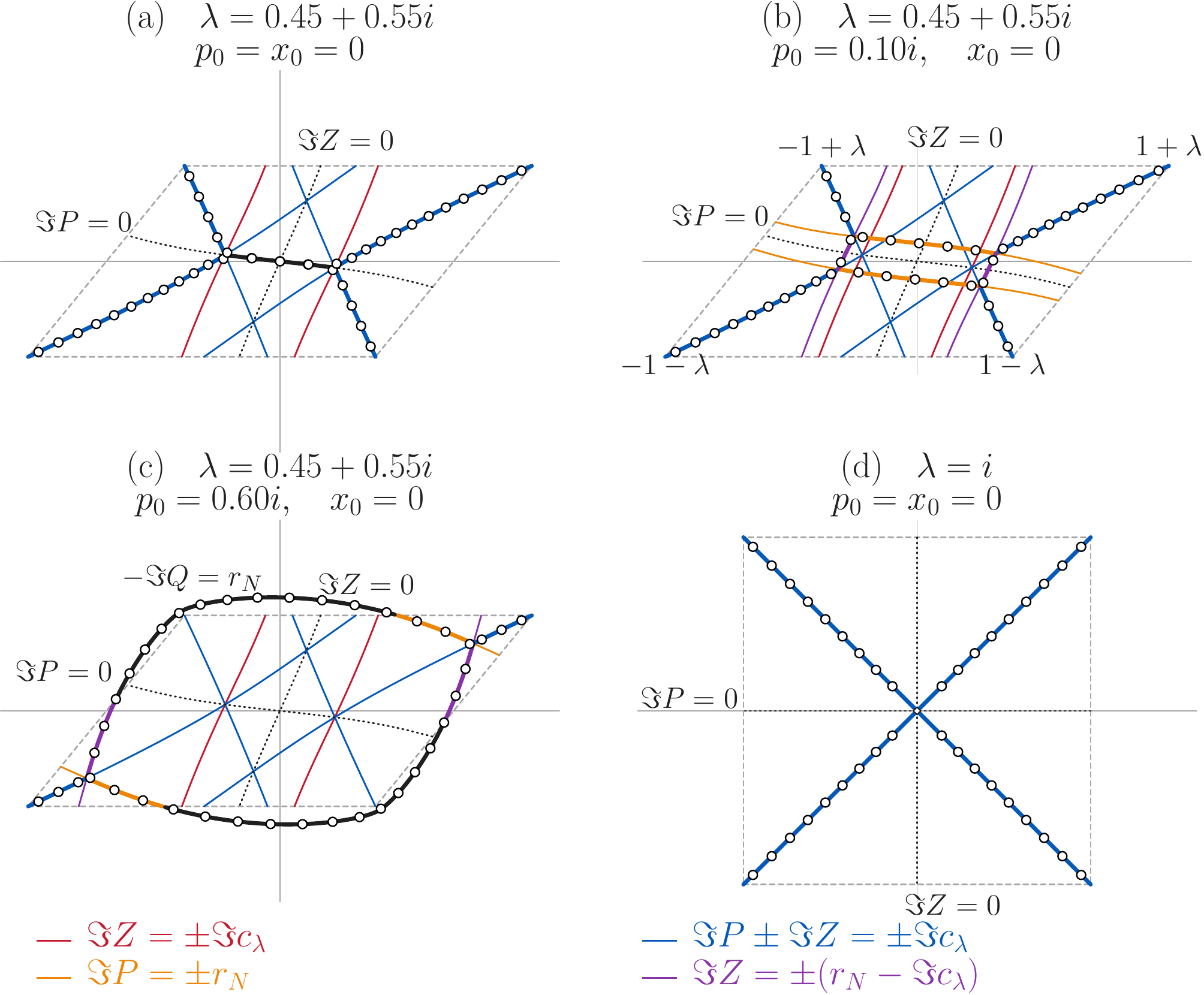}
  \caption{Curves of balancing classical actions and computed eigenvalues (circles),
  $N=40$, with $\Im c_\lambda=\log\lvert\lambda\rvert$ and
  ${r_N=N^{-1}\log\lvert 2R_N\rvert}$. Thick portions indicate dominant balances; dashed boundaries
  mark $\Sigma_\lambda$.}
  \label{fig:energy-diagram}
\end{figure}
}

\newcommand{\realmomentumstripfigure}{%
\begin{figure}[htbp]
  \centering
  \includegraphics[width=0.98\linewidth]{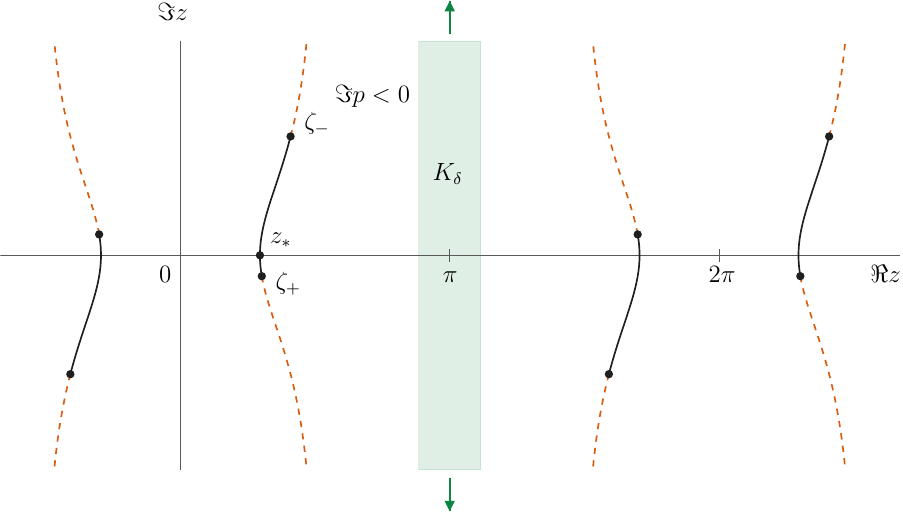}
  \caption{The initial WKB strip $K_\delta$ (green) and real-momentum
  arcs (black) in the $z$-plane. Arc endpoints are turning points;
  $z_*$ is the regular real-axis crossing. Dashed orange extensions
  have $\cos p\in\mathbb R\setminus[-1,1]$ and serve as cuts for the
  initial branch; the continuation domains may use different cuts.}
  \label{fig:real-momentum-strip}
\end{figure}
}

\newcommand{\interiorconstructionfigure}{%
\begin{figure}[htbp]
  \centering
  \includegraphics[width=\linewidth]{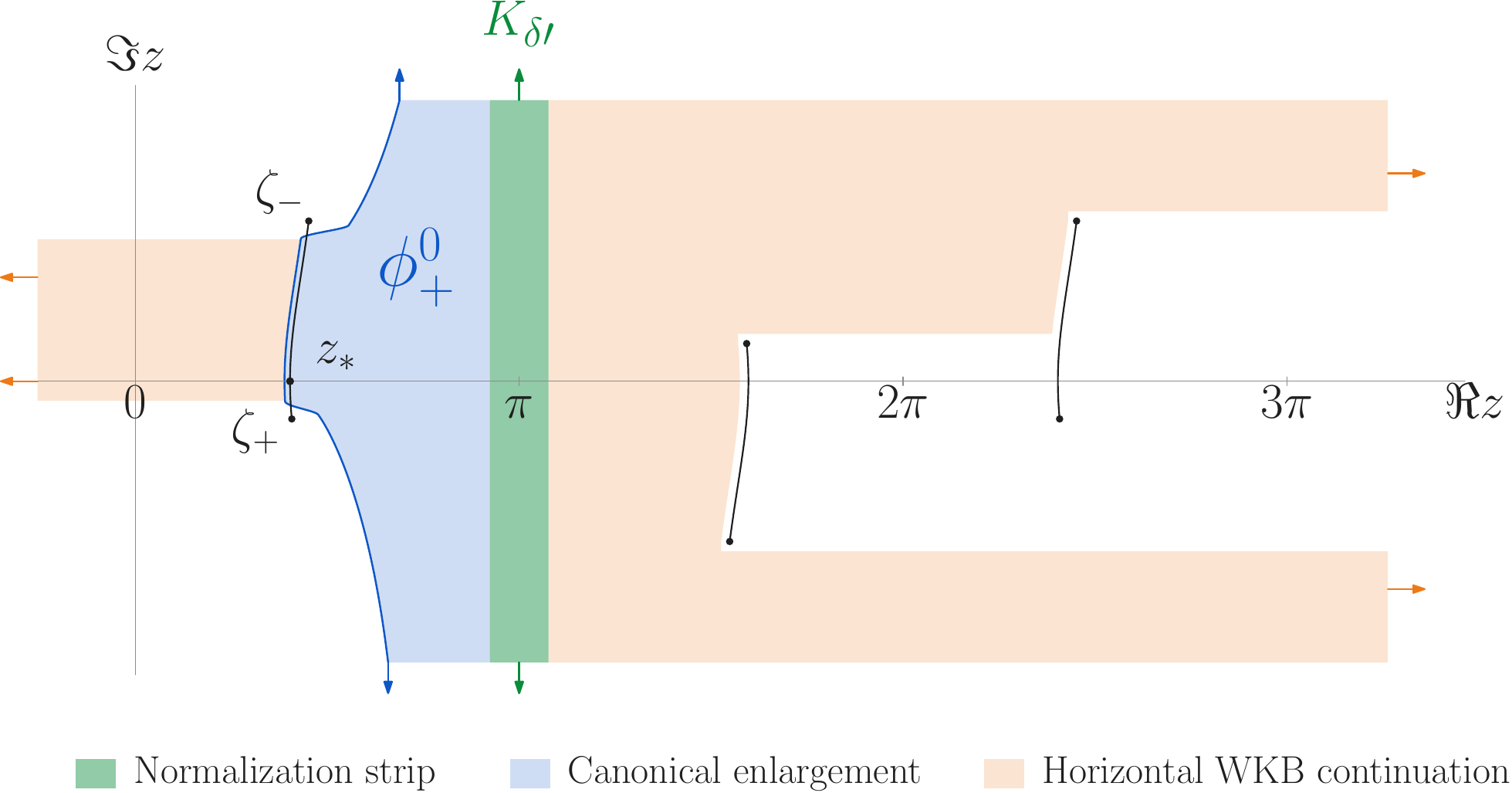}
  \caption{Continuation of the WKB estimate for $\phi_+^0$: the normalization
  strip (green), canonical enlargement across the real-momentum arc (blue),
  and horizontal continuation to a neighborhood of $0$ (orange).}
  \label{fig:interior-domain-construction}
\end{figure}
}

\newcommand{\exteriorwkbcomparisonfigure}{%
\begin{figure}[htbp]
  \centering
  \includegraphics[width=0.94\linewidth,height=0.72\textheight,keepaspectratio]{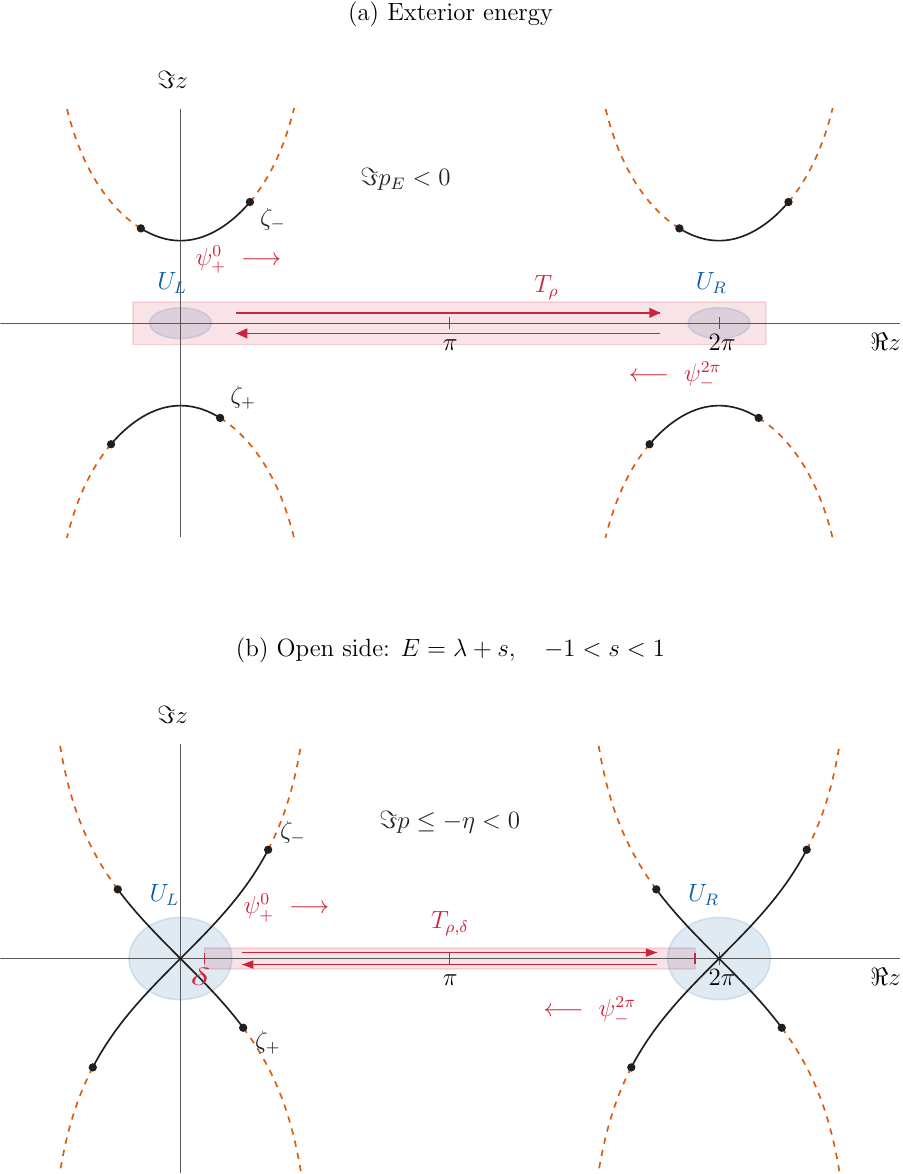}
  \caption{Horizontal WKB comparisons: (a) exterior energies;
  (b) the open side $E=\lambda+s$, $-1<s<1$.
  The red tube surrounds $[0,2\pi]$ in (a) and
  $[\delta,2\pi-\delta]$ in (b). In (b), the momentum is real
  at $0$ and $2\pi$.
  Blue marks the endpoint neighborhoods $U_L,U_R$.
  Arrows indicate rightward plus and leftward minus continuation.
  Black curves are real-momentum arcs; dashed orange extensions have
  real $\cos p$ but nonreal $p$. Shading is schematic.}
  \label{fig:exterior-wkb-comparison}
\end{figure}
}

\newcommand{\actioncontoursfigure}{%
\begin{figure}[htbp]
  \centering
  \includegraphics[width=0.96\linewidth,trim=50bp 30bp 50bp 40bp,clip]{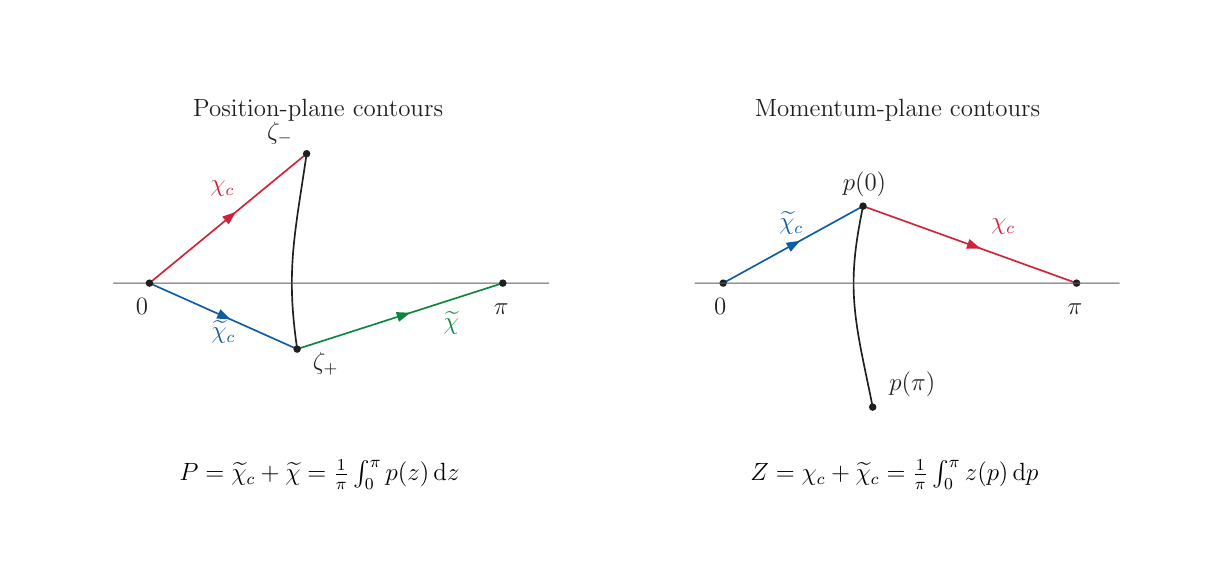}
  \caption{Action contours in the $z$-plane (left) and the corresponding
  $p$-plane paths for $Z=\chi_c+\widetilde\chi_c$ (right).
  Contours avoid the chosen cuts. Black curves are the real-momentum
  arc (left) and the real-position curve from $p(0)$ to $p(\pi)$ (right).}
  \label{fig:action-contours}
\end{figure}
}

\newcommand{\vanishingcyclefigure}{%
\begin{figure}[htbp]
  \centering
  \includegraphics[width=0.96\linewidth]{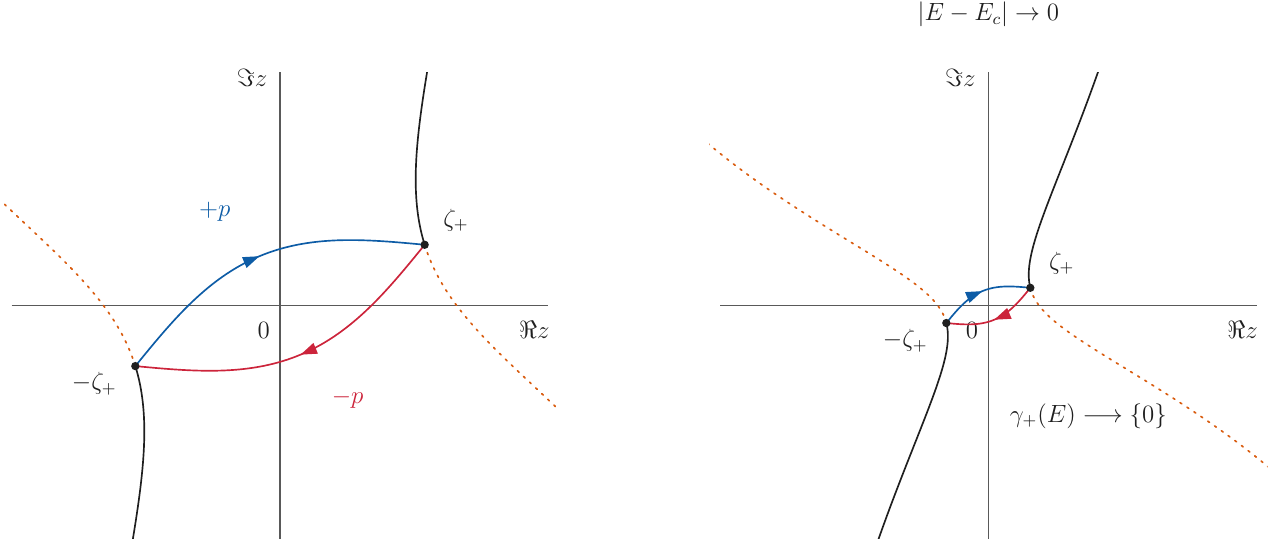}
  \caption{The vanishing cycle at two energies near $E_+=1+\lambda$,
  shown on the same $z$-scale. Blue and red indicate outward and
  return paths on opposite momentum sheets as $\pm\zeta_+$ coalesce
  at zero. Black curves are real-momentum arcs; dotted orange
  extensions have $\cos p\in\mathbb R\setminus[-1,1]$ and are not
  the cuts used for this cycle.}
  \label{fig:vanishing-cycle}
\end{figure}
}

\newcommand{\exteriorcontinuationpathsfigure}{%
\begin{figure}[htbp]
  \centering
  \includegraphics[width=0.96\linewidth]{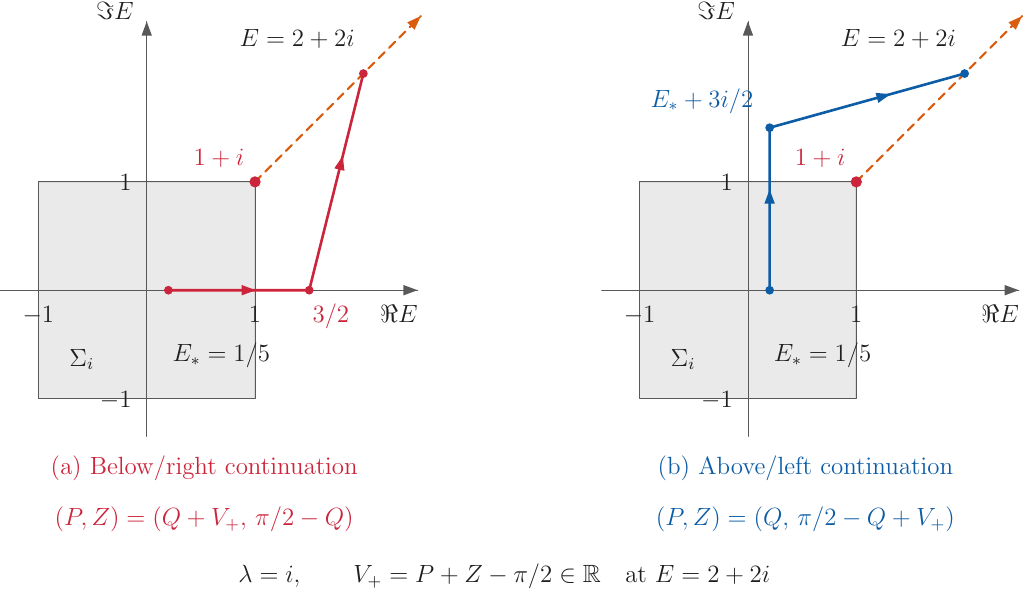}
  \caption{Two continuation paths for $\lambda=i$ from $E=1/5$ to
  $E=2+2i$ on the upper-right equal-growth cut.
  Red passes through $3/2$; blue through $1/5+3i/2$.
  They give different pairs $(P,Z)$ but the same vanishing action $V_+$.}
  \label{fig:exterior-continuation-paths}
\end{figure}
}

\newcommand{\stablecornerdomainfigure}{%
\begin{figure}[htbp]
  \centering
  \includegraphics[width=\linewidth]{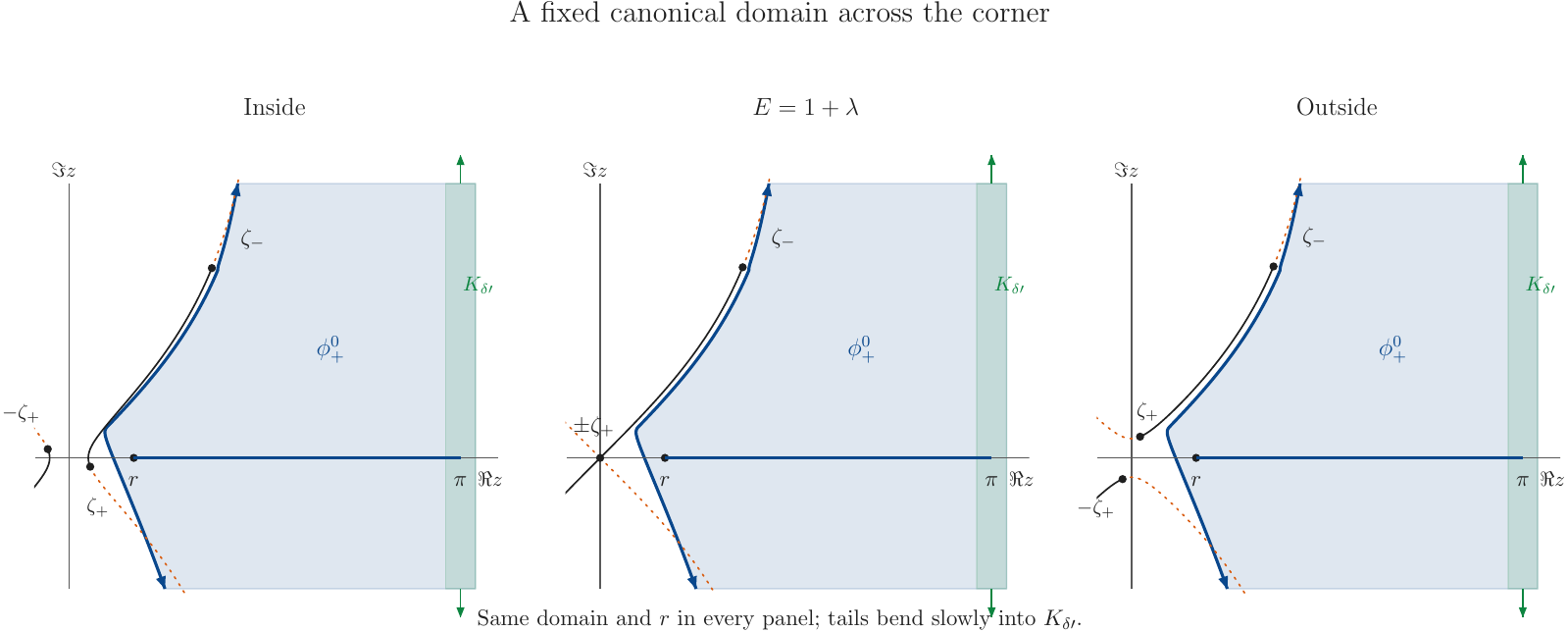}
  \caption{A fixed canonical domain for $\phi_+^0$ (blue) at nearby
  interior, corner, and exterior energies, on the same scale.
  The pair $\pm\zeta_+$ coalesces at zero outside the domain, while
  $\zeta_-$ stays separated. Tails join the green strip $K_{\delta'}$;
  the real interval $[r,\pi]$, for a fixed small $r>0$, lies inside the domain.
  Black curves are real-momentum arcs; dotted orange extensions have
  $\cos p\in\mathbb R\setminus[-1,1]$. Boundary joins and tail bends
  are schematic.}
  \label{fig:stable-corner-domain}
\end{figure}
}

\newcommand{\cornernormalizationfigure}{%
\begin{figure}[htbp]
  \centering
  \includegraphics[width=0.94\linewidth]{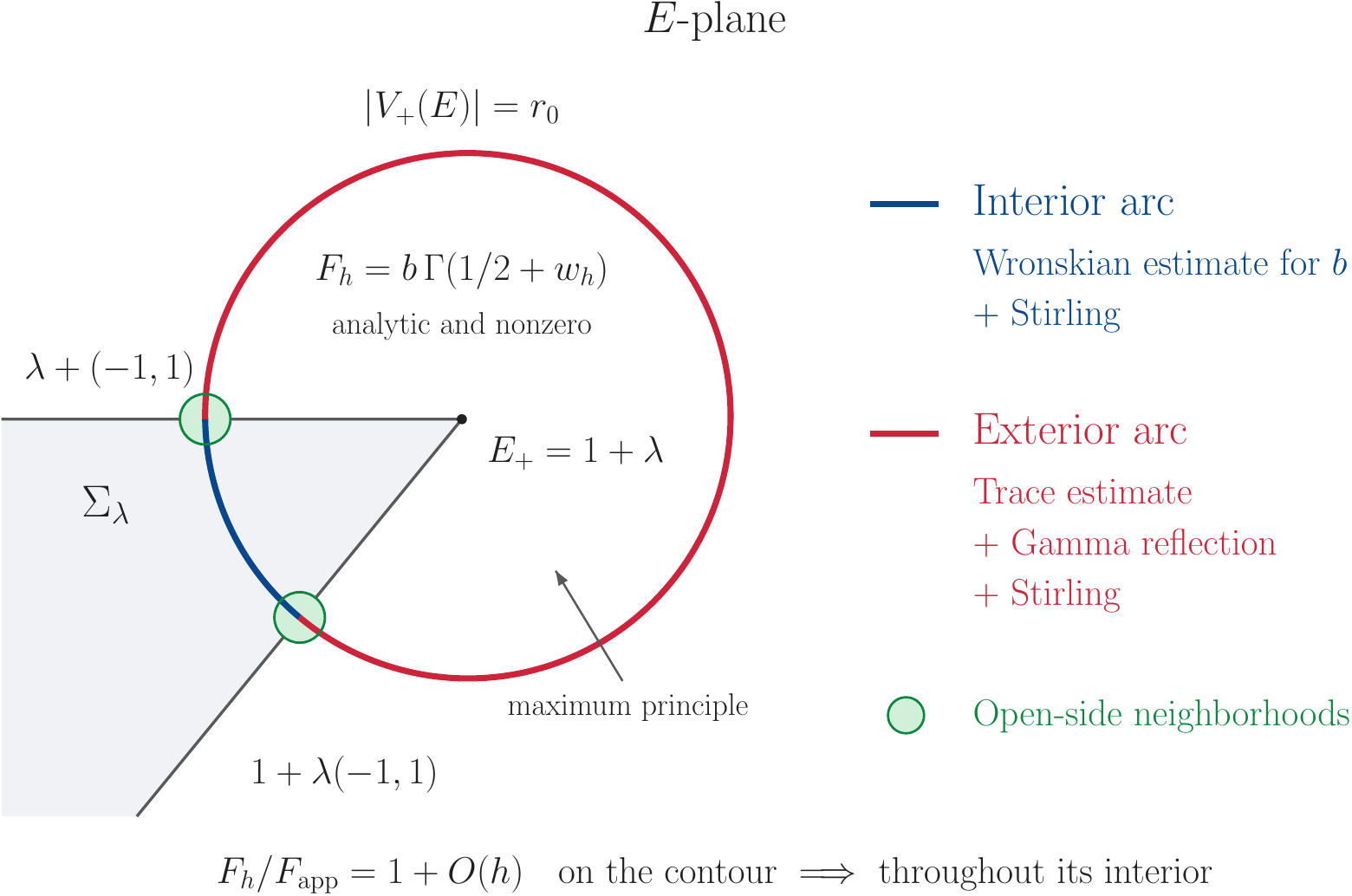}
  \caption{The fixed normalization contour $\lvert V_+(E)\rvert=r_0$
  around $E_+=1+\lambda$ in the $E$-plane. Blue and red mark its arcs
  inside and outside $\Sigma_\lambda$, respectively. Green marks
  neighborhoods of the open-side crossings.}
  \label{fig:corner-boundary-normalization}
\end{figure}
}

\newcommand{\negativelevelsfigure}{%
\begin{figure}[htbp]
  \centering
  \includegraphics[width=\linewidth]{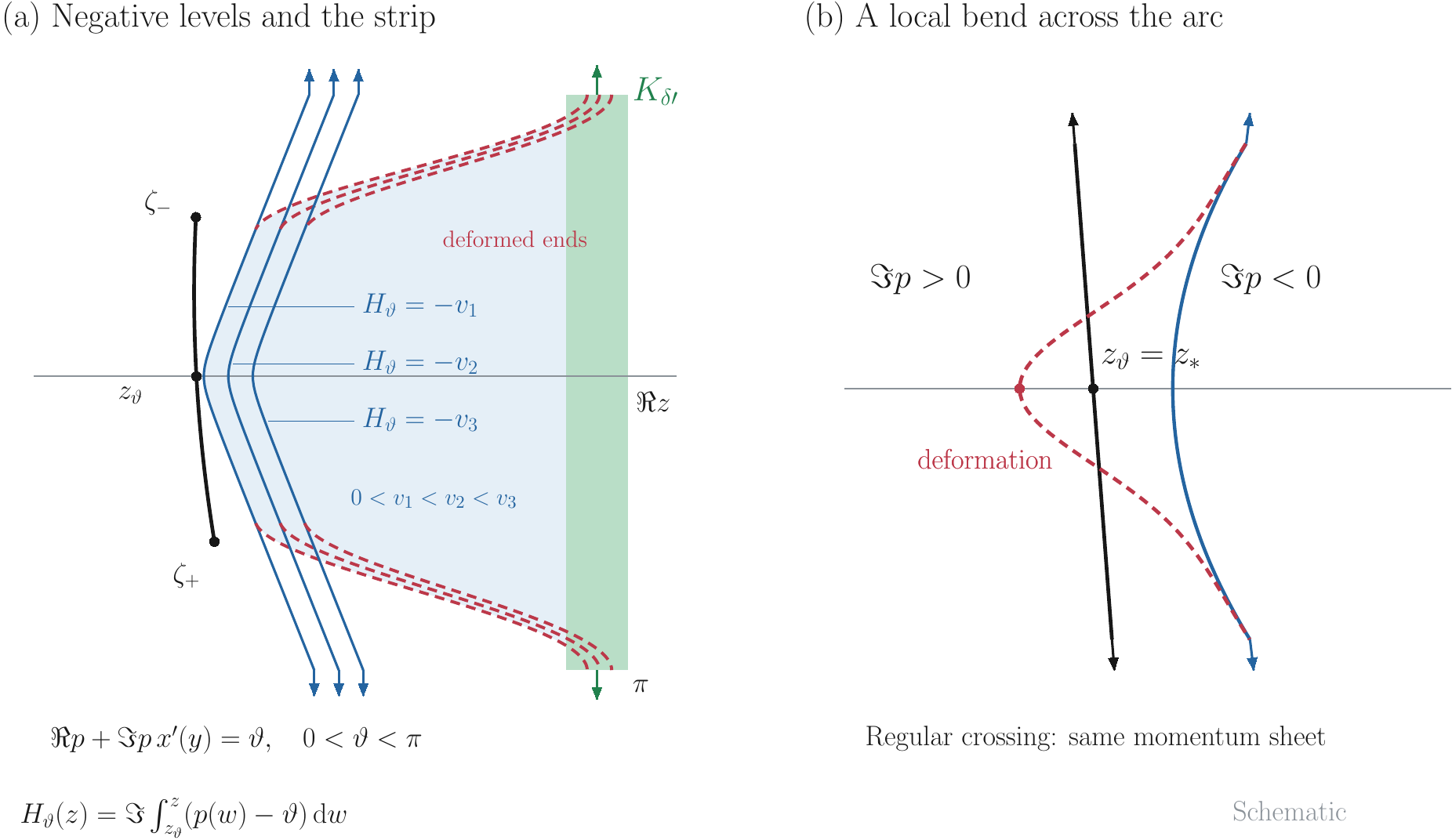}
  \caption{Canonical curves from negative levels of $H_\vartheta$.
  (a) Blue level curves and dashed red deformations joining the green
  normalization strip. (b) A local deformation across the regular
  real-momentum arc, followed by horizontal continuation toward zero.
  The red curves are not level curves; crossing the regular arc
  does not change sheets. The drawing is schematic.}
  \label{fig:negative-level-construction}
\end{figure}
}

\newcommand{\globaltwoactionfigure}{%
\begin{figure}[htbp]
  \centering
  \includegraphics[width=0.9\linewidth]{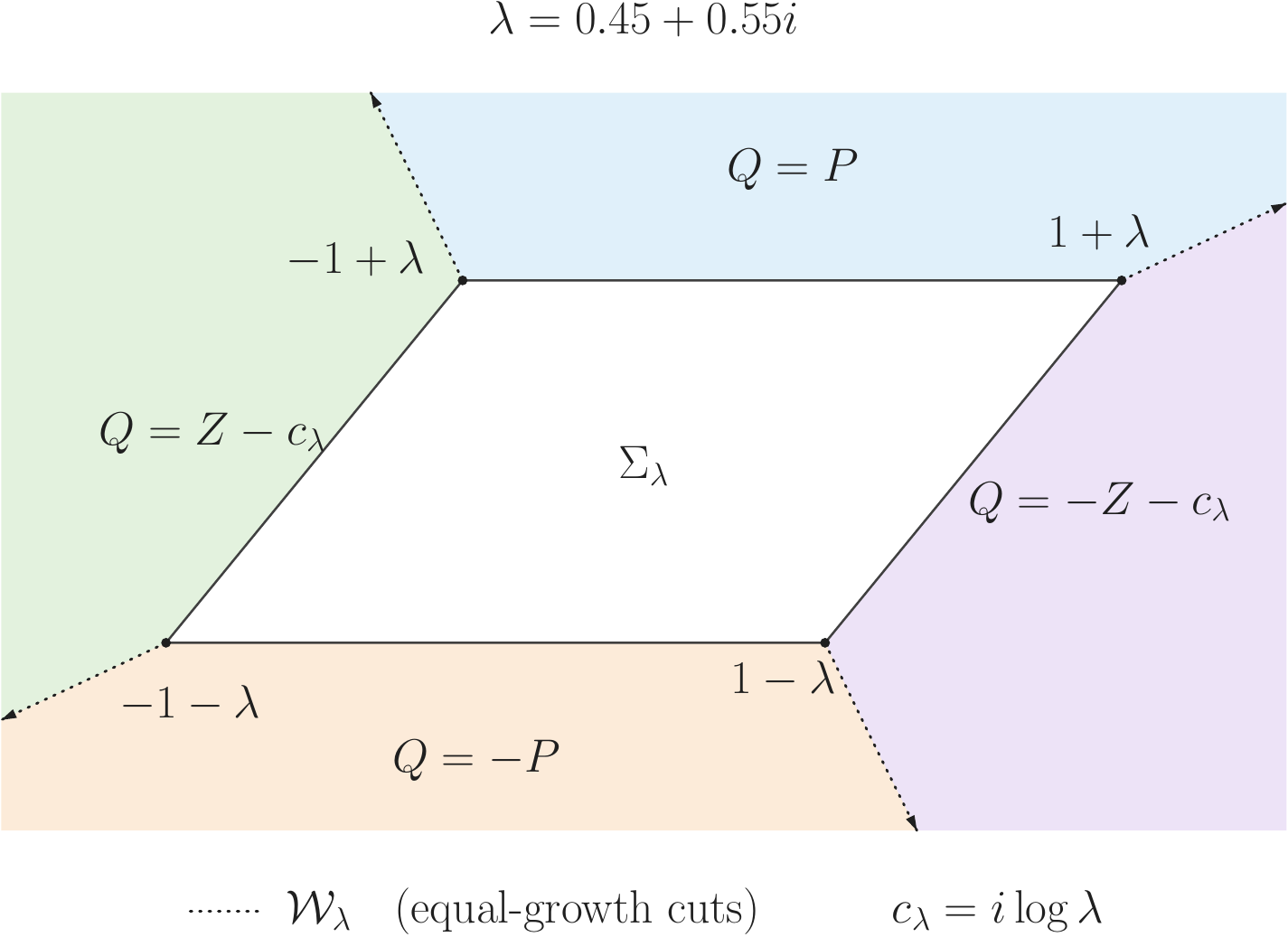}
  \caption{Exterior equal-growth cuts $\mathcal W_\lambda$ (dotted) and the
  dominant action $Q$ in each shaded exterior region. Action identities
  hold modulo $2\pi\mathbb Z$, with $c_\lambda=i\log\lambda$.}
  \label{fig:global-two-action-regions}
\end{figure}
}

\newcommand{\fourexteriorregionsfigure}{%
\begin{figure}[H]
  \centering
  \includegraphics[width=0.9\linewidth]{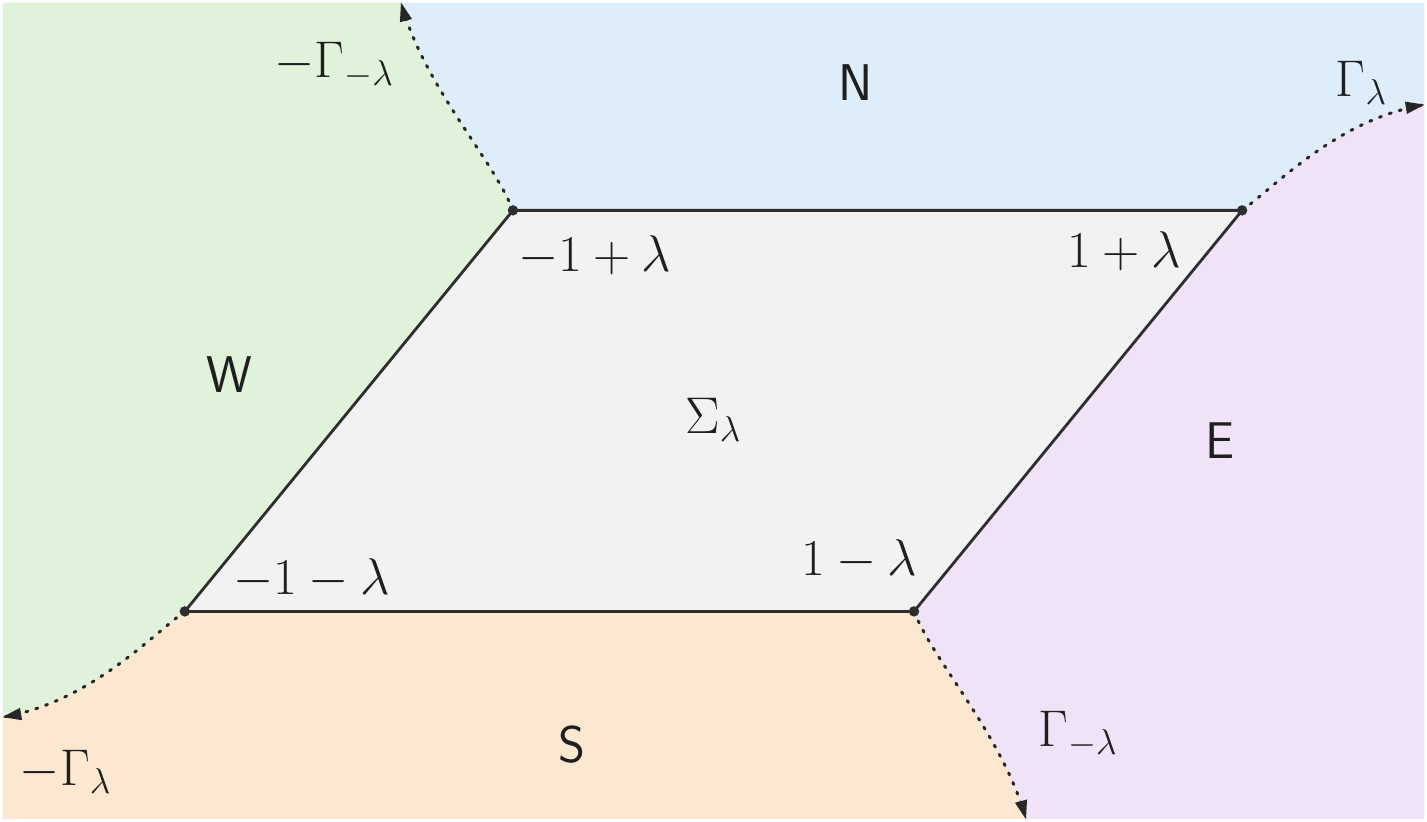}
  \caption{The four exterior regions separated by the cuts issuing from
  the corners of $\Sigma_\lambda$. The curves are schematic.}
  \label{fig:four-exterior-regions}
\end{figure}
}

\crefname{prop}{proposition}{propositions}
\Crefname{prop}{Proposition}{Propositions}
\crefname{defi}{definition}{definitions}
\Crefname{defi}{Definition}{Definitions}
\crefname{corr}{corollary}{corollaries}
\Crefname{corr}{Corollary}{Corollaries}
\crefname{rem}{remark}{remarks}
\Crefname{rem}{Remark}{Remarks}
\crefname{appendix}{appendix}{appendices}
\Crefname{appendix}{Appendix}{Appendices}

\newtheorem{theorem}{Theorem}
\newtheorem{prop}{Proposition}[section]
\newtheorem{defi}[prop]{Definition}
\newtheorem{lemma}[prop]{Lemma}
\newtheorem{corr}[prop]{Corollary}
\newtheorem{example}[prop]{Example}
\newtheorem{rem}{Remark}[section]
\numberwithin{equation}{section}

\DeclareMathOperator{\tr}{tr}

\DeclareMathOperator{\Spec}{Spec}
\newcommand*{\dd}{\mathop{}\!\mathrm{d}}
\renewcommand{\Re}{\mathrm{Re}}
\renewcommand{\Im}{\mathrm{Im}}

\newcommand{\abs}[1]{\left\lvert #1 \right\rvert}
\newcommand{\set}[1]{ \left \{ #1 \right  \}}

\newcommand{\N}{\mathbb{N}}
\newcommand{\Z}{\mathbb{Z}}
\newcommand{\R}{\mathbb{R}}
\newcommand{\C}{\mathbb{C}}
\renewcommand{\phi}{\varphi}

\usepackage{autonum}
\makeatletter
\AtBeginDocument{%
  \autonum@generatePatchedReference{eqref}%
  \autonum@generatePatchedReferenceCSL{Cref}%
}
\makeatother

\makeatletter
\renewcommand{\l@section}
  {\@tocline{1}{0pt}{0pt}{}{}}
\renewcommand{\l@subsection}
  {\@tocline{2}{0pt}{1.5em}{}{}}
\renewcommand{\l@subsubsection}
  {\@tocline{3}{0pt}{3em}{}{}}
\renewcommand{\l@paragraph}
  {\@tocline{4}{0pt}{4.5em}{}{}}
\renewcommand{\l@subparagraph}
  {\@tocline{5}{0pt}{6em}{}{}}
\makeatother

\title{The Spectrum of the Discrete Mathieu Operator with Nonreal Coupling}
\author{Fr{\'e}d{\'e}ric Klopp}
\address[Fr{\'e}d{\'e}ric Klopp]{Sorbonne Universit{\'e},
Universit{\'e} Paris Diderot, CNRS,
Institut de Math{\'e}matiques de Jussieu - Paris Rive Gauche,
F-75005, Paris, France}
\email{frederic.klopp@imj-prg.fr}
\author{Izak Oltman}
\address[Izak Oltman]{Department of Mathematics, Northwestern University, 2033 Sheridan Rd, Evanston, IL 60208}
\email{ioltman@northwestern.edu}
\date{October 19, 2026}

\begin{document}

\begin{abstract}
We study the discrete Mathieu operator with nonreal coupling in the large-period limit, using complex WKB methods for analytic finite-difference equations. We obtain quantization formulas in different regions of the complex energy plane. Away from the corners of the image of the classical symbol, these formulas involve cosines of classical action integrals with small analytic corrections to the actions. Near each corner, we derive a uniform quantization formula involving the reciprocal Gamma function. The quantization formulas constrain spectral accumulation to curves determined by the competing exponential terms. Numerical computations illustrate the resulting spectral geometry.
\end{abstract}

\maketitle
\raggedbottom

\section{Introduction}

For fixed $p_0,x_0\in\C$ and $\lambda\in\C\setminus\R$, define the
operator $H^N_{p_0,x_0,\lambda}$ (for $N\in \mathbb N$) on sequences $u=(u_n)_{n\in\Z}$ by
\begin{align}
 (H^N_{p_0,x_0,\lambda}u)_n
 &=\tfrac12\bigl(e^{ip_0}u_{n+1}+e^{-ip_0}u_{n-1}\bigr)
   +\lambda\cos(2\pi n/N+x_0)u_n.
 \label{eq:finite-mathieu-operator}
\end{align}
We study the eigenvalues of this operator on $N$-periodic sequences
as $N\to\infty$. More precisely, we seek asymptotic formulas for the
values $E\in\C$ for which there is a nonzero sequence $u$ satisfying
\begin{align}
H^N_{p_0,x_0,\lambda}u&=Eu,
& u_{n+N}&=u_n,
& n&\in\Z.
\end{align}
We note that the matrix representation of $H_{p_0,x_0,\lambda}^N$ on $\C^{N}$ is the semiclassical Weyl quantization
on $\R^2/(2\pi\Z)^2$ of the symbol $\cos(p+p_0)+\lambda\cos(x+x_0)$ with semiclassical parameter $h=2\pi/N$.
For background on torus quantization, see
\cite[\S2.3]{ChrZwo2010},
\cite[\S2.2]{DyatlovJezequel2024CatMaps},
and \cite{bouzouina1996equipartition}.

\subsection{Main results}
\label{sec:main-result}

The image of the unshifted classical symbol $\cos p+\lambda\cos x$
on the real torus is the parallelogram
\begin{align}
 \Sigma_\lambda&\coloneq[-1,1]+\lambda[-1,1]
 \label{eq:intro-classical-region}
\end{align}
shown in \Cref{fig:global-energy-regimes}.
A complex number is called an \emph{interior energy} if it belongs to $\operatorname{int}\Sigma_\lambda$, an \emph{exterior energy} if it belongs to $\C\setminus\Sigma_\lambda$,
and a \emph{boundary energy} if it belongs to the boundary of $\Sigma_\lambda$.  An \emph{open side} is a side of
$\Sigma_\lambda$ with its two endpoints removed.  

\globalregimefigure

For $E\in\operatorname{int}\Sigma_\lambda$, define the two classical actions
\begin{equation}\label{eq:intro-classical-actions}
  \begin{split}
     P(E)&\coloneq\frac1{2\pi}\int_0^{2\pi}\arccos(E-\lambda\cos z)\dd z,\\
 Z(E)&\coloneq\frac1{2\pi}\int_0^{2\pi}\arccos((E-\cos p)/\lambda)\dd p,
  \end{split}
\end{equation}
using the inverse cosines with real parts in $(0,\pi)$.\footnote{Thus $P(E) = (2\pi)^{-1}\int_{0}^{2\pi} p(z)\dd z$ and $Z(E) = (2\pi)^{-1}\int_0^{2\pi}z(p) \dd p$, where $p$ and $z$ are branches of analytic functions satisfying $\cos p + \lambda \cos (z) = E$.}
Outside $\operatorname{int}\Sigma_\lambda$, contours and branches are continued together.

Set
\begin{align}
 R_N(p_0,x_0,\lambda)
 &\coloneq\cos(Np_0)+\lambda^N\cos(Nx_0).
 \label{eq:intro-bloch-trace}
\end{align}
Away from the corners of $\Sigma_\lambda$, eigenvalues of $H^N_{p_0,x_0,\lambda}$ satisfy the following quantization condition.

\begin{theorem}[Quantization away from corners]
 \label{thm:main}
\normalfont
Fix $\lambda\in\C\setminus\R$ and
$E_*\in\C\setminus\{\pm1\pm\lambda\}$.

\begin{itemize}
\item[\textup{(i)}] \emph{Interior and open sides.}
If $E_*$ belongs to $\Sigma_\lambda$, there are a
complex neighborhood $U$ of $E_*$ and analytic functions
\begin{align}
 P_N&=P+O(N^{-2}),& Z_N&=Z+O(N^{-2})
 \label{eq:intro-corrected-actions}
\end{align}
(with $P$ and $Z$ defined in \eqref{eq:intro-classical-actions}) such that, for all sufficiently large $N$ and all $E\in U$,
\begin{align}
 E\in\Spec(H^N_{p_0,x_0,\lambda})
 \quad\Longleftrightarrow\quad
 \cos(NP_N)
 +\lambda^N\cos(NZ_N)
 &=R_N(p_0,x_0,\lambda).
 \label{eq:intro-spectral-condition}
\end{align}

\item[\textup{(ii)}] \emph{Outside $\Sigma_\lambda$.}
If $E_*\notin\Sigma_\lambda$, choose the periodic momentum branch
$p(x)$, locally analytic in $E$, by
\begin{align}
 \cos p(x)&=E-\lambda\cos x,
 &\Im p(x)&<0,
\end{align}
and put
\begin{align}
 Q(E)&\coloneq\frac1{2\pi}\int_0^{2\pi}p(x)\dd x.
 \label{eq:intro-exterior-action}
\end{align}
There are a neighborhood $U\subset\C\setminus\Sigma_\lambda$ of $E_*$
and an analytic function $Q_N=Q+O(N^{-2})$ such that, for all
sufficiently large $N$ and all $E\in U$,
\begin{align}
 E\in\Spec(H^N_{p_0,x_0,\lambda})
 \quad\Longleftrightarrow\quad
 \cos(NQ_N)
 &=R_N(p_0,x_0,\lambda).
 \label{eq:intro-exterior-spectral-condition}
\end{align}
\end{itemize}

\noindent The functions $P_N,Z_N$, and $Q_N$ (which we refer to as
\textit{analytic correction actions}) are independent of
$p_0,x_0$. All the error terms are uniform in $E$.
\end{theorem}

\Cref{thm:main}\textup{(ii)} is the simplest global description of the exterior.
Section~\ref{sec:global-two-action-continuation} explains how \Cref{thm:main}\textup{(i)} continues to the exterior energy region.

We next give a quantization condition on a fixed neighborhood of a corner of $\Sigma_\lambda$ (which we from now on just refer to as a corner), including
the corner itself. We state it near $1+\lambda$ for $\Im\lambda>0$.
Conjugation and sign symmetries give the other corners.

Define
\begin{align}
 w(E,N)&\coloneq
 \frac{N}{2\pi}
 \bigl(P(E)+Z(E)+i\log\lambda\bigr),
 \qquad 0<\arg\lambda<\pi.
 \label{eq:intro-classical-corner-coordinate}
\end{align}
The term $P+Z+i\log\lambda$ is the action that
vanishes at $1+\lambda$, as proved in
\Cref{lem:vanishing-action-two-periods}. The actions $P,Z$
are continued together from the interior of $\Sigma_\lambda$.

\begin{theorem}[Quantization near the corner]
\label{thm:main-corner}
\normalfont
Fix $\lambda$ with $\Im\lambda>0$. There is a complex
neighborhood $U$ of $1+\lambda$, independent of $N$, and
constants $C,\eta>0$ such that, for all sufficiently large $N$,
\begin{itemize}
\item there are functions $w_N,a_N,r_{+,N}$ analytic in $U$
satisfying, uniformly on $U$,
\begin{align}
 w_N&=w+O(N^{-1}),
 &a_N&=1+O(N^{-1}\log N),\\
 \lvert r_{+,N}(E)\rvert
 &\le C\max\{1,\lvert\lambda\rvert\}^{N}e^{-\eta N}.
\end{align}
\item for $E\in U$, $E\in\Spec(H^N_{p_0,x_0,\lambda})$ if and only if
\begin{align}
 \sqrt{\frac\pi2}\,\lambda^{N/2}e^{iN(P-Z)/2}
 \frac{e^{w(1-\log w)}}{\Gamma(\tfrac12-w_N)}\,a_N(E)
 +r_{+,N}(E)&=R_N(p_0,x_0,\lambda).
 \label{eq:intro-corner-spectral-condition}
\end{align}
\end{itemize}
Here, we choose $\log w$ to be real for $w$ positive.
The left-hand side of \eqref{eq:intro-corner-spectral-condition}
is understood by analytic continuation through $E=1+\lambda$.
All terms on the left-hand side of~\eqref{eq:intro-corner-spectral-condition}
are independent of $p_0$ and $x_0$.
\end{theorem}

For real phase shifts, the poles of
$\Gamma(\tfrac12-w_N)$ give the leading corner quantization levels.
Taylor expansion of the vanishing action $P+Z+i\log\lambda$ at $1+\lambda$
then gives an explicit formula for each fixed level.

\begin{corr}
\label{cor:real-phase-corner-quantization}
Assume $\Im\lambda>0$ and $p_0,x_0\in\R$.
There is $\varepsilon>0$ such that, for all sufficiently
large $N$ and every integer $0\le n\le\varepsilon N$, there
is a unique simple eigenvalue $E_{n,N}\in U$ satisfying
\begin{align}
 P(E_{n,N})+Z(E_{n,N})+i\log\lambda
 &=\frac{2\pi}{N}\left(n+\tfrac12\right)+O(N^{-2}),
\end{align}
uniformly in $n,x_0$, and $p_0$.
In particular, for each fixed $n\ge0$,
\begin{align}
 E_{n,N}
 &=1+\lambda-\frac{2\pi\sqrt{\lambda}}{N}
   \left(n+\tfrac12\right)+O(N^{-2}).
\end{align}
\end{corr}
The proof is given in Appendix~\ref{app:corner-B-zeros}.

In \Cref{fig:global-energy-regimes}, \Cref{thm:main}\textup{(i)} applies
to the green interior and open sides, and \Cref{thm:main}\textup{(ii)}
to the blue exterior. \Cref{thm:main-corner} and its symmetry reductions
apply to the amber corner neighborhoods. These overlap the regular
regions, and their drawn sizes are schematic.

\begin{rem}
\label{rem:three-regime-transitions}
The corresponding formulas at the other corners of $\Sigma_\lambda$
and for all corners when $\Im\lambda<0$ follow from
Theorem~\ref{thm:main-corner} and
\Cref{cor:corner-symmetry-reduction}.

The left-hand sides of the quantization conditions in \Cref{thm:main}
and \Cref{thm:main-corner} are the same phase-independent polynomial
$m_N(E,\lambda)/2$, where $m_N$ is constructed
as the trace of a monodromy matrix.
Its existence and the spectral identity $m_N=2R_N$ are proved in
\Cref{lem:trace-polynomial}. Thus, on overlapping neighborhoods away
from the corners,
\begin{align}
 \cos(NP_N)
 +\lambda^N\cos(NZ_N)
 &=\cos(NQ_N)
\end{align}
exactly. Reversing the sign of an action or adding $2\pi k$,
$k\in\Z$, leaves its cosine unchanged.
The corner expression (\Cref{thm:main-corner}) agrees with this formula wherever their
neighborhoods overlap. At any corner $E_c\in\{\pm1\pm\lambda\}$,
it also covers $N^{-1}\ll\lvert E-E_c\rvert \ll1$, where the large-argument Gamma
asymptotics recover the regular WKB formulas. At the corner itself
the leading coefficient changes, so the unchanged classical
two-cosine expression does not extend through it (see \Cref{prop:corner-obstruction} to see exactly how the two-term action fails at the corner).
\end{rem}

\begin{rem}
As elementary computations show, quantization of the symbol
$\cos(p+p_0)+\lambda\cos(x+x_0)$ for all $\lambda\notin\R$ and
$(p_0,x_0)\in\C^2$ is equivalent, up to multiplication
by a nonzero constant, to quantizing all the symbols
$ae^{ip}+be^{-ip}+ce^{ix}+de^{-ix}$ for $(a,b,c,d)\in\C^4$ such
that $abcd\ne0$ and $cd/(ab)\notin(0,\infty)$.
{
Indeed, after choosing square roots, set
$\lambda=\sqrt{cd}/\sqrt{ab}$,
$e^{ip_0}=a/\sqrt{ab}$, and $e^{ix_0}=c/\sqrt{cd}$. Then
\begin{align}
 ae^{ip}+be^{-ip}+ce^{ix}+de^{-ix}
 &=2\sqrt{ab}\bigl(\cos(p+p_0)+\lambda\cos(x+x_0)\bigr).
\end{align}
If $cd/(ab)\in(0,\infty)$, this reduction gives real coupling,
but the phase shifts may still be complex. In particular,
$a=\overline b$ and $c=\overline d$ give the self-adjoint case.
}
When $abcd=0$, the equation (possibly modulo a Fourier transform)
reduces to a first-order finite-difference equation or
to the spectral problem for a multiplication operator, both of which
can be treated explicitly.
\end{rem}

\subsection{Interpretation of the results}
\label{sec:interpretation}

The quantization formulas separate two questions: where can the
exponential terms cancel, and which energies satisfy their phase
condition?  We first describe the resulting energy diagrams and
their spectral consequences, then
specialize to the Scottish flag matrix.
The two-action comparisons take place on fixed compact subsets of
the neighborhoods in \Cref{thm:main}\textup{(i)}, or of their exterior
continuations in \Cref{cor:global-two-action}.  Near the corners we
use the reciprocal-Gamma formula of \Cref{thm:main-corner}.

\subsubsection{Energy diagrams}
\label{sec:energy-diagram}

The quantization conditions of \Cref{thm:main} describe both where eigenvalues can
accumulate and how individual eigenvalues are selected.  To see the
first point, expand \eqref{eq:intro-spectral-condition} as
\begin{align}
 e^{iNP_N}+e^{-iNP_N}
 +\lambda^N e^{iNZ_N}+\lambda^N e^{-iNZ_N}-2R_N&=0.
 \label{eq:diagram-five-term-equation}
\end{align}
If one term on the left-hand side is
exponentially larger than all the others, cancellation is impossible
for large $N$.
We can compare their size to show eigenvalues concentrate on curves and use their phase to give more precise estimates of where they lie.

Using the principal logarithm, set
\begin{align}
 c_\lambda&\coloneq i\log\lambda,
 &\Im c_\lambda&=\log\lvert\lambda\rvert.
 \label{eq:exterior-log-period}
\end{align}
The logarithm of the magnitude of each term in \eqref{eq:diagram-five-term-equation}, divided by $N$, is given by the right column below,
with an $O(N^{-2})$ error in the first four entries:
\begin{center}
\begin{tabular}{c|c}
Term & Growth rate\\ \hline
$e^{iNP_N}$ & $-\Im P$\\
$e^{-iNP_N}$ & $\Im P$\\
$\lambda^N e^{iNZ_N}$ & $\Im(c_\lambda-Z)$\\
$\lambda^N e^{-iNZ_N}$ & $\Im(c_\lambda+Z)$\\
$-2R_N$ & $r_N\coloneq N^{-1}\log\lvert 2R_N\rvert$
\end{tabular}
\end{center}
If $R_N=0$, set $r_N=-\infty$, omit the last row, and omit
the balances involving $R_N$ in \Cref{tab:dominant-balances}.

For \eqref{eq:diagram-five-term-equation} to hold, the two largest
growth rates must differ by at most $O(N^{-1})$.
The following proposition makes this bound precise.

\begin{prop}
\label{prop:spectral-growth-gap}
Let $g_{1,N}(E)\geq g_{2,N}(E)$ be the two largest growth rates
listed in the growth-rate table above.
Every eigenvalue in a fixed compact subset of a two-action
neighborhood satisfies
\begin{align}
 g_{1,N}(E)-g_{2,N}(E)
 &\leq\frac{\log4}{N}+O(N^{-2}),
 \label{eq:spectral-growth-gap}
\end{align}
uniformly there.  A region where one rate exceeds all the others by
a fixed positive amount therefore contains no eigenvalues for
sufficiently large $N$.
\end{prop}
\begin{proof}
At a zero of \eqref{eq:diagram-five-term-equation}, the largest term
in magnitude is at most the sum of the other four, hence at most
four times the second largest.  Taking logarithms gives the bound
for the corrected actions.  Replacing $P_N,Z_N$ by $P,Z$ changes
each rate, including the two largest, by $O(N^{-2})$.
\end{proof}

By considering all such cases, we can derive the following table listing the equalities between the largest
classical growth rates and the inequalities ensuring that
no other rate is larger.

\begin{table}[htbp]
\centering
\caption{Dominant balances in the two-action quantization condition.}
\label{tab:dominant-balances}
\setlength{\tabcolsep}{4pt}
\renewcommand{\arraystretch}{1.3}
\begin{tabular}{@{}lll@{}}
Competing terms & Balance equality & Condition for dominance\\ \hline
Two $P$-terms & $\Im P=0$
 & $\Im c_\lambda+\lvert\Im Z\rvert\le0,\quad r_N\le0$\\
Two $Z$-terms & $\Im Z=0$
 & $\lvert\Im P\rvert\le\Im c_\lambda,\quad r_N\le\Im c_\lambda$\\
$P$-term and $Z$-term
 & $\lvert\Im P\rvert=\Im c_\lambda+\lvert\Im Z\rvert$
 & $\lvert\Im P\rvert\ge r_N$\\
$P$-term and $2R_N$ & $\lvert\Im P\rvert=r_N$
 & $\Im c_\lambda+\lvert\Im Z\rvert\le r_N$\\
$Z$-term and $2R_N$ & $\Im c_\lambda+\lvert\Im Z\rvert=r_N$
 & $\lvert\Im P\rvert\le r_N$
\end{tabular}
\end{table}

\Cref{fig:energy-diagram} numerically illustrates this. We (numerically) compute eigenvalues of $H_{p_0,x_0,\lambda}^N$ for varying parameters and overlay the balancing equalities from \Cref{tab:dominant-balances}.
This gives a good approximation of the eigenvalue location, and partially explains why the spectrum accumulates to ``spectral arcs''.

\energydiagramfigure

\begin{example}[The Scottish flag matrix]
\label{sec:scottish-flag}
We now use our analysis to analyze the spectrum of the Scottish flag matrix\footnote{
This matrix has been studied in various contexts for its spectral instability and connection to geometric quantization, see for instance \cite{BorthwickUribe,TrefethenChapman2004,TrefethenEmbree}.
Numerically, its spectrum was observed to lie on crossed lines in the complex plane (a rough approximation of the flag of Scotland -- leading to its name).
Only recently has this spectral confinement been rigorously established, see \cite{BeckerOltman2026Scottish}.} where $\lambda=i$ and $p_0=x_0=0$.
In this case, the corresponding matrix (for $N >4$) is
\begin{align}
  H^N_{0,0,i}
=
\begin{pmatrix}
i & \frac12 & 0 & \cdots & \frac12 \\
\frac12 & i\cos(2\pi/N) & \frac12 & \ddots & 0 \\
0 & \frac12 & i\cos(4\pi/N) & \ddots & \vdots \\
\vdots & \ddots & \ddots & \ddots & \frac12 \\
\frac12 & 0 & \cdots & \frac12
& i\cos(2\pi(N-1)/N)
\end{pmatrix}.
\end{align}

Write $H^N_{0,0,i}=C_N+iD_N$ for its hopping and diagonal parts.
Both are self-adjoint contractions.  If
$H^N_{0,0,i}u=Eu$ and $\lVert u\rVert =1$, then
\begin{align}
 \Re E&=\langle C_Nu,u\rangle\in[-1,1],&
 \Im E&=\langle D_Nu,u\rangle\in[-1,1].
\end{align}
Thus the spectrum lies in $\Sigma_i = [-1,1]+i[-1,1]$.

We now restrict our analysis to $N=4m$ for simplicity.
In this case $R_N = 1 + (i)^N=2$ and $\Im c_\lambda = 0$.
Therefore, by \Cref{tab:dominant-balances} the spectrum concentrates on the sets
\begin{align}
\{ \Im Z =  \Im P\} \cup \{\Im Z =- \Im P \}.
\end{align}

Using
\begin{align}
 \cos a+\cos b=2\cos\frac{a+b}{2}\cos\frac{a-b}{2},
\end{align}
we rewrite the spectral equation in \Cref{thm:main}, away from the corners, as
\begin{align}
 \cos\left(\frac{N(P_N+Z_N)}2\right)
 \cos\left(\frac{N(P_N-Z_N)}2\right)&=1.
 \label{eq:scottish-balanced-factorization}
\end{align}
By \Cref{lem:scottish-bulk-actions}, on compact subsegments
of each diagonal arm away from the center and corners, one factor is
exponentially large.  At an eigenvalue near such a subsegment, the
other factor is therefore exponentially small.  This yields the
classical Bohr--Sommerfeld conditions
\begin{align}
 \frac1{2\pi}\left(\int_0^{2\pi}p(z)\dd z
 +\int_0^{2\pi}z(p)\dd p\right)
 &=\frac{2\pi}{N}\left(n+\frac12\right),\\
 \frac1{2\pi}\left(\int_0^{2\pi}z(p)\dd p
 -\int_0^{2\pi}p(z)\dd z\right)
 &=\frac{2\pi}{N}\left(n+\frac12\right),
 \label{eq:scottish-classical-diagonal-quantization}
\end{align}
on the diagonals toward $1+i$ and $-1+i$, respectively, with $n\in\mathbb Z$.
Here $p(z)=\arccos(E-i\cos z)$ and
$z(p)=\arccos((E-\cos p)/i)$ use the branches in
\eqref{eq:intro-classical-actions}.  These equalities define the
classical candidate locations: each candidate in a fixed compact
subsegment approximates a unique simple eigenvalue within $O(N^{-2})$,
and every eigenvalue near the subsegment satisfies the corresponding
condition with an $O(N^{-2})$ remainder.  The opposite diagonals follow by
$E\mapsto-E$.

Moreover, by \Cref{lem:scottish-explicit-action}, the two action combinations
have the same radial expression, for $0<t<\sqrt2$:
\begin{align}
 S(t)
 &\coloneq P(e^{i\pi/4}t)+Z(e^{i\pi/4}t)-\frac\pi2\\
 &=Z(e^{3i\pi/4}t)-P(e^{3i\pi/4}t)+\frac\pi2\\
 &=\frac2\pi\int_0^{\arccos(t^2/2)}
 \arccos\left(\frac{t}{\sqrt{2\cos u}}\right)\dd u.
 \label{eq:scottish-interior-action-explicit}
\end{align}
Appendix~\ref{app:scottish-bulk} proves that $S$
decreases strictly from $S(0)=\pi/2$ to $S(\sqrt2)=0$.
Since $\pi/2=(2\pi/N)m$, the shifts only relabel the quantization integer.
The Bohr--Sommerfeld candidates used in the numerics in \Cref{fig:scottish-flag-numerics} are therefore
\begin{align}
 S(t_{n,N})&=\frac{2\pi}{N}\left(n+\frac12\right),&
 n&=0,\ldots,m-2.
 \label{eq:scottish-interior-Bohr-Sommerfeld}
\end{align}
For the numerical comparison, we omit the innermost classical
candidate $n=m-1$, so rotation of $e^{i\pi/4}t_{n,N}$
gives $N-4$ candidates on the four diagonals.
At $E=0$, neither cosine factor in
\eqref{eq:scottish-balanced-factorization} is exponentially large,
so the reduction to a one-factor quantization condition does not
apply.

The actions remain regular at $E=0$, where $\Im P=\Im Z=0$.
The four classical exponential terms therefore have equal magnitude,
and the corrected terms remain comparable.

It turns out that $E=0$ is an eigenvalue of multiplicity $4$.
Indeed, by
\cite[Theorem~3.4]{BeckerOltman2026Scottish},
\begin{align}
 \det(EI_N-H^N_{0,0,i})
 &=E^4\prod_{r=1}^{m-1}(E^4+\tau_{r,N}),&
 \tau_{r,N}&>0.
 \label{eq:scottish-exact-factorization}
\end{align}

For the corners, we apply \Cref{cor:real-phase-corner-quantization} with
$\lambda=i$ and $p_0=x_0=0$.
We obtain, for each fixed $n\ge0$,
a unique simple eigenvalue near $1+i$ satisfying
\begin{align}
 E_{n,N}
 &=E_{n,N}^{\mathrm{BS}}+O(N^{-2}),&
 E_{n,N}^{\mathrm{BS}}
 &\coloneq
 1+i-\frac{\pi(2n+1)}{\sqrt2N}(1+i).
 \label{eq:scottish-corner-Bohr-Sommerfeld}
\end{align}
The corner symmetries give the corresponding formulas
at the other three corners.

\Cref{fig:scottish-flag-numerics} compares the computed eigenvalues
with the bulk predictions \eqref{eq:scottish-interior-Bohr-Sommerfeld}
in panels (a),(c) and the corner predictions
\eqref{eq:scottish-corner-Bohr-Sommerfeld} in panels (b),(d).

\begin{figure}[htbp]
 \centering
 \includegraphics[width=\linewidth]{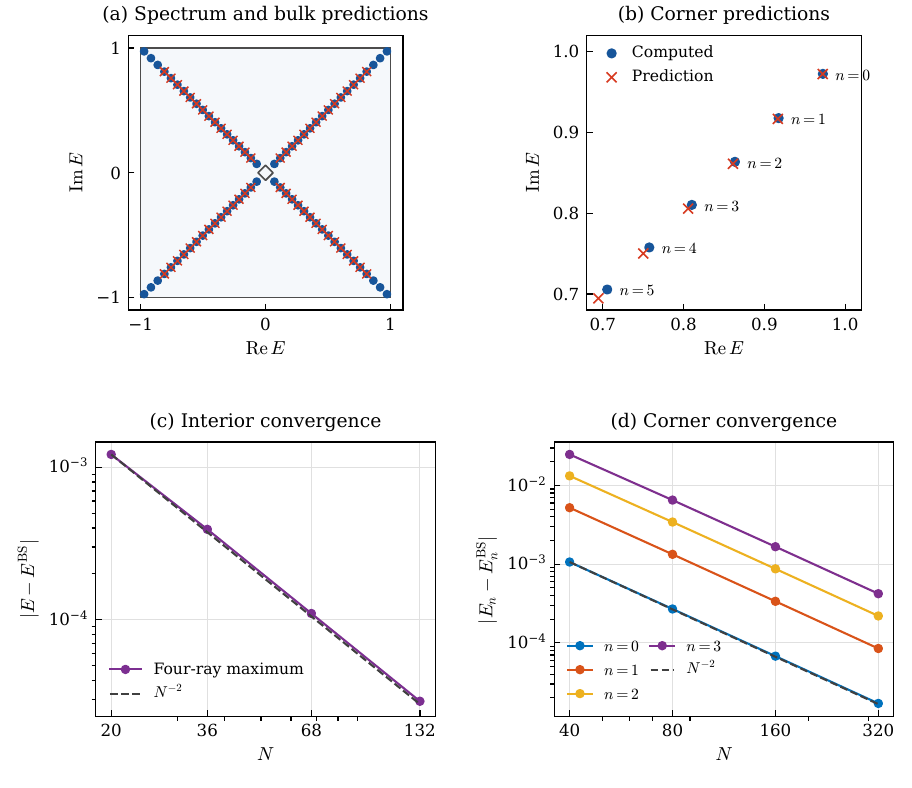}
 \caption{Scottish-flag numerics, $\lambda=i$, $p_0=x_0=0$.
 Blue dots are numerically computed eigenvalues; red crosses are
 Bohr--Sommerfeld predictions. At $N=80$: (a) eigenvalue comparison away from the origin and corners; (b) the first six upper-right corner modes; (c) maximum error (away from the origin and corners) over four rays at
 $(n+1/2)/N=1/8$; (d) corner errors for $n=0,1,2,3$.
 Dashed lines show the reference rate $N^{-2}$.}
 \label{fig:scottish-flag-numerics}
\end{figure}

\end{example}

\subsection{Relation to earlier work}

Trefethen and Chapman~\cite{TrefethenChapman2004} studied the
Scottish flag matrix as an example of a nonnormal matrix with
exponentially accurate, localized pseudomodes.  Their work explains
the large difference between its spectrum and its pseudospectra.
For small random perturbations of quantized tori, Christiansen and
Zworski~\cite{ChrZwo2010} proved a Weyl law for the expected eigenvalue
count; Vogel~\cite{Vogel2020Tori} established a Weyl law with high
probability.  Oltman~\cite{Oltman2023ProbabilisticWeyl} generalized the
probabilistic Weyl law to randomly perturbed Berezin--Toeplitz operators
on compact K\"ahler manifolds equipped with positive line bundles.
For the unperturbed almost Mathieu matrices, Becker and
Oltman~\cite{BeckerOltman2026Scottish} derive exact finite-dimensional
identities, prove confinement to perpendicular lines under specified
phase and parity conditions, and compute limiting eigenvalue measures.

Our approach builds on the complex WKB method for analytic
finite-difference equations developed by Buslaev and
Fedotov~\cite{B-F:94}.  
For a survey of the complex WKB method, see \cite{Fed23}.
We use the canonical-domain constructions of
Fedotov and Klopp~\cite{Klopp2004} and the difference-equation
estimates of Fedotov and Shchetka~\cite{Fed2017,FedotovShchetka2022Berry}.
Fedotov and Klopp~\cite{FedotovKlopp2019Airy} give uniform
approximations near simple turning points.  For two coalescing turning
points, Fedotov studies semiclassical solution
asymptotics~\cite{Fedotov2021Coalescent} and obtains uniform
transition-matrix asymptotics by comparing WKB bases away from the
turning points~\cite{Fedotov2022Transition}.
We combine WKB comparisons with exact Wronskian identities to obtain
a trace formula for arbitrary complex phases throughout the energy plane.

Analytic Berezin--Toeplitz methods also give quantization results.
Deleporte and Le Floch~\cite[Theorem~1]{DeleporteLeFloch2025BohrSommerfeld}
establish Bohr--Sommerfeld rules with exponentially small errors near
regular energy levels on compact real-analytic K\"ahler manifolds of
complex dimension one.  They treat small non-self-adjoint perturbations
of self-adjoint operators, with perturbation size independent of the
semiclassical parameter.  Near critical energies,
Deleporte~\cite{Deleporte2023WKB} constructs WKB quasimodes
near nondegenerate wells with exponentially small errors.
R\'eguer~\cite[Theorem~1.1]{Reguer2026EllipticToeplitz} obtains local
eigenvalue quantization near elliptic critical values for analytic
non-self-adjoint Toeplitz operators.  After identifying the Toeplitz
symbol through the correspondence of
Rouby~\cite{Rouby2017TorusQuantization}, one could apply this result
to local corner quantization in the real-phase problem.  We prefer
a direct Wronskian argument, which also determines the nonvanishing
trace prefactor needed for arbitrary complex phases.

\subsection{Outline of the proof}

The proof reduces the spectral equation to Wronskians of solutions to
analytic finite difference equations and estimates them
by WKB comparison.

\paragraph{\textbf{Exact reduction.}}
Section~\ref{sec:preliminary-reductions} reduces the spectral problem
to estimating the trace of a monodromy matrix that is $2\pi/N$-periodic. The trace is a first order
trigonometric polynomial and its Fourier coefficients of order $\pm1$
are computed explicitly. Thus, it suffices to compute its $0$-order
Fourier coefficient, say, $m_N$ that is polynomial in $E$ and $\lambda$.
This coefficient inherits a number of symmetries from those of the
symbol $\cos p+\lambda\cos x$. These are used to reduce the domain in
$(E,\lambda)$ where the coefficient has to be computed.

The components of the monodromy matrix in the WKB basis used below
are ratios of Wronskians of
solutions to an analytic finite difference equation which makes the
use of the complex WKB approximation method developed for these
equations very natural. The $(1,1)$-coefficient of this matrix is
constant, say, $b$; its $(1,2)$-coefficient is
the first-order trigonometric polynomial
$(c+ae^{-iNz})/(-2i)$ in the normalization used below.
Exact identities reduce the computation of $m_N$ to that of $b$ and $a/c$.
Aubry--Andr\'e duality relates these two quantities.

\paragraph{\textbf{Regular energies.}}
Section~\ref{sec:away-from-corners} proves \Cref{thm:main} by
evaluating $b$ throughout the interior of $\Sigma_\lambda$.
Aubry--Andr\'e duality then yields $a/c$. The large $N$ asymptotics of
these coefficients are expressed in terms of phases
whose principal terms are the actions $P$ and $Z$.
Outside the interior of $\Sigma_\lambda$, away from its corners,
i.e., on its open sides and in the exterior of $\Sigma_\lambda$,
horizontal continuation of the WKB asymptotics yields $m_N$ directly.
Comparing the actions recovers the two-action formula across the sides.

\paragraph{\textbf{Corners.}}
Section~\ref{sec:corner-quantization} proves \Cref{thm:main-corner}.
It continues one WKB comparison through the corner and uses the exact
trace identity to locate the zeros of the remaining Wronskian.
Multiplication by a Gamma factor cancels these zeros.
WKB estimates on a fixed energy contour around the corner and the
maximum principle determine the remaining nonzero factor inside
the contour. By the symmetries of $m_N$, the calculation needs only
to be performed for the corner $1+\lambda$ with $\Im\lambda>0$.

\paragraph{\textbf{Global continuation.}}
Section~\ref{sec:global-two-action-continuation} combines exterior WKB
with the corner actions. Four equal-growth curves for
$\lvert\Im P\rvert$ and
$\log\lvert\lambda\rvert+\lvert\Im Z\rvert$
provide natural cuts for the asymptotics of
$m_N$ outside which $P$ and $Z$ are continued starting from the
interior of $\Sigma_\lambda$ and the two-action formulas obtained in
the interior continue to hold locally.

Appendix~\ref{app:wkb-tools} contains the uniform WKB estimates,
normalization, and comparison domains used throughout the paper.
Appendix~\ref{app:fourier-normalization} develops the Fourier
transformation used in the Aubry--Andr\'e duality argument.
Appendix~\ref{app:exterior-cuts} studies the equal-growth curves, and
Appendix~\ref{app:uniform-corner-comparison} proves the corner estimates.
Finally, the action calculations for the Scottish flag in
\Cref{sec:scottish-flag} are performed in
Appendix~\ref{app:action-identity}.

\subsection{Notation}
\label{sec:notation}

The following table collects the notation used throughout the paper.
Later entries are listed here for reference; their definitions remain
at the cited places.

{\small
\renewcommand{\arraystretch}{1.12}
\begin{longtable}{|p{0.16\textwidth}|p{0.47\textwidth}|p{0.27\textwidth}|}
\hline
\textbf{Variable} & \textbf{Brief explanation} &
\textbf{Definition or first use}\\ \hline
\endfirsthead
\hline
\textbf{Variable} & \textbf{Brief explanation} &
\textbf{Definition or first use}\\ \hline
\endhead
$N$, $h=2\pi/N$ & Matrix size and analytic finite-difference step &
\eqref{eq:finite-mathieu-operator};
\Cref{sec:comp-asympt-m_n}\\ \hline

$p_0,x_0,\lambda$ & Momentum phase, position phase, and nonreal coupling &
\eqref{eq:finite-mathieu-operator}\\ \hline

$H^N_{p_0,x_0,\lambda}$ & Cyclic discrete Mathieu operator &
\eqref{eq:finite-mathieu-operator}\\ \hline

$\Sigma_\lambda$ & Classical parallelogram $[-1,1]+\lambda[-1,1]$ &
\eqref{eq:intro-classical-region}\\ \hline

$m_N$ & Phase-independent trace polynomial &
\Cref{lem:trace-polynomial}\\ \hline

$R_N$ & Phase-dependent side of the spectral equation &
\eqref{eq:intro-bloch-trace}\\ \hline

$c_\lambda$ &
Constant $c_\lambda=i\log\lambda$,
with $\Im c_\lambda=\log\lvert\lambda\rvert$ &
\eqref{eq:exterior-log-period}\\ \hline

$r_N$ & Growth rate
$r_N=N^{-1}\log\lvert 2R_N\rvert$, with $r_N=-\infty$ if $R_N=0$ &
\Cref{sec:energy-diagram}\\ \hline

$p(z)$, $p(z;E)$ & Momentum branches satisfying
$\cos p(z)+\lambda\cos z=E$ &
\eqref{eq:cosine-momentum};
\eqref{eq:exterior-momentum}\\ \hline

$\zeta_+,\zeta_-$ & Turning points with $\cos p=1$ and
$\cos p=-1$, respectively &
\eqref{eq:turning-points}\\ \hline

$K_\delta$ & Initial canonical strip $\lvert\Re z-\pi\rvert<\delta$ &
\eqref{eq:unbounded-pi-strip}\\ \hline

$\phi_\pm^0,\phi_\pm^{2\pi}$ & Strip-normalized WKB solutions and
their one-period translates &
\Cref{prop:entire-WKB-normalization};
\Cref{sec:unified-period-map}\\ \hline

$W[\phi,\psi]$ & Discrete Wronskian &
\eqref{eq:discrete-wronskian}\\ \hline

$T(E,h;z_0)$ & Transfer matrix carrying the solution basis backward
by one period on $z_0+h\Z$ &
\eqref{eq:unified-connection-matrix};
\eqref{eq:unified-wronskian-matrix}\\ \hline

$\mathcal A,\mathcal B,\mathcal C,\mathcal D$ &
Four Wronskians encoding the one-period connection matrix &
\eqref{eq:four-wronskians}\\ \hline

$a,b,c$ & $b=\mathcal B/\mathcal A$ and the two coefficients in
$\mathcal C(z)=c+a e^{2\pi iz/h}$ &
\eqref{eq:scalar-coefficient-notation};
\eqref{eq:exact-C-two-modes}\\ \hline

$E^\dagger,\lambda^\dagger$ & Conjugated dual energy and coupling &
\eqref{eq:coefficient-dual-parameters}\\ \hline

$\chi_c,\widetilde\chi_c,\widetilde\chi$ &
Action integrals between $0$, the turning points, and $\pi$, normalized by $1/\pi$ &
\eqref{eq:action-definitions}\\ \hline

$P,Z$ & The two classical actions, normalized by $1/(2\pi)$ &
\eqref{eq:intro-classical-actions}\\ \hline

$P_N,Z_N$\newline $P_h,Z_h$ &
Analytically corrected actions, with $P_N=P_h$ and $Z_N=Z_h$
when $h=2\pi/N$ &
\eqref{eq:intro-corrected-actions};
\eqref{eq:exact-quantum-actions}\\ \hline

$Q,Q_N,Q_h$ &
Classical and corrected exterior actions, with $Q_N=Q_h$
when $h=2\pi/N$ &
\eqref{eq:intro-exterior-action};
\eqref{eq:exterior-quantum-action-Wronskian}\\ \hline

$E_c,V_c,V_{c,h}$\newline $E_+,V_+,V_h$ &
Corner energy, vanishing action, and corrected vanishing action.
The second line denotes the corner $1+\lambda$ &
\eqref{eq:corner-energy-definition};
\eqref{eq:corner-action-definition};
\eqref{eq:upper-right-vanishing-action-normalization};
\eqref{eq:corner-corrected-action}\\ \hline

$w,w_N,w_h$ &
Scaled corner actions $w=V_+/h$ and
$w_N=w_h=V_h/h$ when $h=2\pi/N$ &
\eqref{eq:intro-classical-corner-coordinate};
\Cref{thm:main-corner};
\eqref{eq:corner-corrected-action}\\ \hline

$a_N,a_h$ & Corner amplitude correction, with $a_N=a_h$ when
$h=2\pi/N$ (distinct from the Fourier coefficient $a$) &
\eqref{eq:corner-intro-correction}\\ \hline

$r_{c,N},r_{+,N}$ &
Exponentially small remainders in the corner formulas &
\Cref{thm:main-corner};
\eqref{eq:all-corner-WKB-matching}\\ \hline

$\Gamma_\lambda,\mathcal W_\lambda$ &
Outward equal-growth curve from $1+\lambda$ and the union
of the four corner cuts &
\eqref{eq:outward-equal-growth-curve};
\eqref{eq:global-cut-system}\\ \hline
\end{longtable}
}

\paragraph{\textbf{Acknowledgments.}} We thank Maciej Zworski for initially suggesting this problem.
Partial support for this project was provided by Maciej Zworski's Simons Foundation's Targeted Grant Award No. 896630 ``Moir\'e Materials Magic''.
The second author was supported by NSF grant DMS-2136217.

\paragraph{\textbf{AI Usage.}} This work began in 2022, several years before competent large language models were available.
In June of 2026, we began using Anthropic and ChatGPT models (including Claude Opus, Claude Fable, ChatGPT Sol 5.6, and ChatGPT 6 Astra) for filling in gaps, finding simplifications, proofreading, and reorganizing.
The vast majority of the ideas of this paper are our own.
But we must acknowledge that many details, simplifications, and ideas were suggested or implemented by AI, including (but not limited to) the details of the corner quantization condition and the idea to use canonical strips around $\pi$.
We have checked all AI generated arguments and take full responsibility for any errors or omissions.

\tableofcontents

\section{Preliminary reductions}
\label{sec:preliminary-reductions}

This section expresses the spectral condition through the coefficients
$b$ and $a/c$.
First we remove the phase variables $p_0$ and $x_0$ and reduce the quantization problem to estimating a polynomial
$m_N(E,\lambda)$.  
To estimate $m_N(E,\lambda)$ we study (the seemingly more complicated) analytic finite difference equation generalizing \eqref{eq:finite-mathieu-operator}. 
$m_N(E,\lambda)$ is expressed as Wronskians of solutions to this analytic finite-difference equation.
We apply analytic finite-difference WKB methods to estimate these Wronskians.

\subsection{The trace polynomial and Aubry--Andr\'e duality}
\label{sec:finite-reduction}

The phase variables $p_0,x_0$ can be separated from the energy by an
elementary trace calculation.

\begin{lemma}[The trace polynomial]
\label{lem:trace-polynomial}
There is a polynomial $m_N(E,\lambda)$ of degree $N$ in $E$ (with leading coefficient $2^N$) and at most
$N$ in $\lambda$ such that
\begin{align}
 E\in\Spec(H^N_{p_0,x_0,\lambda})
 \quad\Longleftrightarrow\quad
  m_N(E,\lambda)=2R_N(p_0,x_0,\lambda).
 \label{eq:finite-spectral-condition}
\end{align}
More precisely,
\begin{align}
 2^N\det\bigl(EI_N-H^N_{p_0,x_0,\lambda}\bigr)
  &=m_N(E,\lambda)-2R_N(p_0,x_0,\lambda).
 \label{eq:finite-characteristic-determinant}
\end{align}
Here $I_N$ is the $N\times N$ identity matrix.
\end{lemma}
\begin{proof}
\textbf{Floquet reduction.}
Suppose $H_{p_0,x_0,\lambda}^N u = E u$.
Setting $v_n \coloneq e^{inp_0}u_n$, we get that
\begin{align}
  \begin{pmatrix}
    v_{N+1}\\v_N
  \end{pmatrix} = M_N(0,x_0,E,\lambda) \begin{pmatrix}
    v_1\\v_0
  \end{pmatrix}
\end{align}
where the transfer matrix $M_N(0,x_0,E,\lambda)$ is
\begin{align}
  M_N(0,x_0,E,\lambda) &\coloneq G_N\cdots G_1,&
  G_n &\coloneq \begin{pmatrix}
    2(E-\lambda\cos(x_0+2\pi n/N))&-1\\1&0  \end{pmatrix}.
  \label{eq:finite-transfer-matrix}
\end{align}
Because $v_{n+N}=e^{iNp_0}v_n$, $H_{p_0,x_0,\lambda}^N u = E u$ has a nonzero periodic solution $u$ exactly when $e^{iNp_0}$ is an eigenvalue of $M_N$.
Because $\det G_j = 1$, the other eigenvalue of $M_N$ is $e^{-iNp_0}$, thus $\tr M_N = 2\cos(Np_0)$.

For fixed $x_0$ and $\lambda$, the polynomial $E\mapsto\tr M_N$ has degree
$N$ and leading coefficient $2^N$.  Unless $2\cos(Np_0)$ is one of
its finitely many critical values, $\tr M_N-2\cos(Np_0)$ has $N$
distinct roots.  By the equivalence just proved, these are exactly the
roots of the characteristic polynomial.  Comparing leading coefficients
therefore gives
\begin{align}
 2^N\det\bigl(EI_N-H^N_{p_0,x_0,\lambda}\bigr)
 &=\tr M_N(0,x_0,E,\lambda)-2\cos(Np_0).
\end{align}
Both sides depend continuously on $p_0$. Since the identity holds whenever the roots are distinct, approximation by such phases proves it for every $p_0$.

\textbf{Dependence on the position phase.}
The upper-left entry of $G_n$ is
\begin{align}
    2E-\lambda e^{ix_0}e^{2\pi in/N}
       -\lambda e^{-ix_0}e^{-2\pi in/N}.
\end{align}
Thus $\tr(G_N\cdots G_1)$ is a linear combination of
$e^{ikx_0}$ with $-N\le k\le N$.
Replacing $x_0$ by $x_0+2\pi/N$ cyclically permutes the
factors, leaving the trace unchanged. Consequently, only
the terms with $k=0,\pm N$ can occur.

The coefficient of $e^{\pm iNx_0}$ comes from choosing
$-\lambda e^{\pm ix_0}e^{\pm2\pi in/N}$ in every upper-left
entry. It is therefore
\begin{align}
    (-\lambda)^N
    \exp\left(\pm\frac{2\pi i}{N}\sum_{n=1}^N n\right)
    =(-\lambda)^N(-1)^{N+1}
    =-\lambda^N.
\end{align}
Writing $m_N(E,\lambda)$ for the term independent of $x_0$ gives
\begin{align}
    \tr M_N(0,x_0,E,\lambda)
    &=m_N(E,\lambda)-2\lambda^N\cos(Nx_0).
    \label{eq:trace-polynomial-decomposition}
\end{align}
Each factor has degree at most one in $E$ and in $\lambda$,
so $m_N$ has degree at most $N$ in each variable.
Choosing $2E$ in every upper-left entry gives the unique
term of degree $N$ in $E$, namely $2^NE^N$.
\end{proof}

Thus, to determine the spectrum of $H_{p_0,x_0,\lambda}^N$, it suffices to analyze $m_N(E,\lambda)$ which is independent of $x_0$ and $p_0$.
In computing $m_N(E,\lambda)$, we may therefore take $p_0 = x_0 = 0$.
 
We next observe the following symmetry of $m_N(E,\lambda)$ allowing us to interchange the role of the position and momentum.

\begin{lemma}[Aubry--Andr\'e duality]
\label{lem:position-momentum-duality}
For $\lambda\ne0$,
\begin{align}
 m_N(E,\lambda)&=\lambda^N m_N(E/\lambda,\lambda^{-1}).
 \label{eq:mN-duality}
\end{align}
\end{lemma}
\begin{proof}
Let $u=(u_n)$ be an $N$-periodic sequence, and define its discrete
Fourier transform by
\begin{align}
    \widehat u_k
    &\coloneq \frac{1}{\sqrt N}
    \sum_{n=0}^{N-1}u_n e^{2\pi ink/N}.
\end{align}
Writing $H=H^N_{p_0,x_0,\lambda}$, we compute
\begin{align}
    \widehat{Hu}_k
    &=\frac{1}{\sqrt N}\sum_{n=0}^{N-1}
    \left[
        \frac12\left(e^{ip_0}u_{n+1}+e^{-ip_0}u_{n-1}\right)
        +\lambda\cos(2\pi n/N+x_0)u_n
    \right]e^{2\pi ink/N}\\
    &=\frac12\left(
        e^{ip_0-2\pi ik/N}+e^{-ip_0+2\pi ik/N}
    \right)\widehat u_k
    +\frac{\lambda}{2}\left(
        e^{ix_0}\widehat u_{k+1}
        +e^{-ix_0}\widehat u_{k-1}
    \right)\\
    &=\lambda\left[
        \frac12\left(
            e^{ix_0}\widehat u_{k+1}
            +e^{-ix_0}\widehat u_{k-1}
        \right)
        +\lambda^{-1}\cos(2\pi k/N-p_0)\widehat u_k
    \right]\\
    &=\lambda
    \bigl(H^N_{x_0,-p_0,\lambda^{-1}}\widehat u\bigr)_k.
\end{align}
In the second equality, we reindex the hopping terms using the
$N$-periodicity of $u$ and expand the cosine into exponentials.
In particular,
\begin{align}
    H^N_{p_0,x_0,\lambda}u=Eu
    \quad\Longleftrightarrow\quad
    H^N_{x_0,-p_0,\lambda^{-1}}\widehat u
    =\frac{E}{\lambda}\widehat u.
\end{align}

Since the Fourier transform is invertible, the operator identity
above and \eqref{eq:finite-characteristic-determinant} give
\begin{align}
    &m_N(E,\lambda)-2\cos(Np_0)-2\lambda^N\cos(Nx_0)\\
    &\qquad=2^N\det\bigl(EI_N-H^N_{p_0,x_0,\lambda}\bigr)\\
    &\qquad=2^N\lambda^N
    \det\bigl((E/\lambda)I_N
        -H^N_{x_0,-p_0,\lambda^{-1}}\bigr)\\
    &\qquad=\lambda^N m_N(E/\lambda,\lambda^{-1})
        -2\lambda^N\cos(Nx_0)-2\cos(Np_0).
\end{align}
The second equality uses the Fourier conjugation and factors
$\lambda$ from each of the $N$ rows. The last applies
\eqref{eq:finite-characteristic-determinant} to the transformed
operator and uses $\cos(-Np_0)=\cos(Np_0)$.
Canceling the cosine terms proves the claim.
\end{proof}

Together with \Cref{lem:position-momentum-duality}, the following elementary symmetries will let us
prove several geometric cases only once.  Complex conjugation and shifts
of $p_0$ or $x_0$ by $\pi$ give
\begin{align}
 m_N(\overline E,\overline\lambda)
 &=\overline{m_N(E,\lambda)},
 \label{eq:mN-conjugation}\\
 m_N(E,-\lambda)&=m_N(E,\lambda),&
 m_N(-E,\lambda)&=(-1)^N m_N(E,\lambda).
 \label{eq:mN-sign-symmetries}
\end{align}
When studying energies in the interior of $\Sigma_\lambda$, we use complex conjugation to assume $\Im\lambda>0$. No restriction on $\lvert \lambda\rvert $ is needed.

\subsection{The analytic difference equation}
\label{sec:comp-asympt-m_n}

We now consider the analogous (and seemingly much more complicated) analytic finite difference equation to \eqref{eq:finite-mathieu-operator}.
Set $h\coloneq 2\pi/N$ and consider analytic functions $\psi$ solving
\begin{align}
 \frac12\bigl(\psi(z+h)+\psi(z-h)\bigr)
       +\lambda\cos z\,\psi(z)&=E\psi(z).
 \label{eq:main-cos}
\end{align}

We stress here that it is not at all obvious why a nonzero analytic
solution to \eqref{eq:main-cos} even exists, or how studying this
much more rigid equation (compared to
\eqref{eq:finite-mathieu-operator}) gives necessary and sufficient
conditions for eigenvalues.
Existence results and WKB estimates for such equations were
developed by Buslaev and Fedotov \cite{B-F:94} and by Fedotov and
Shchetka \cite{FedotovShchetka2022Berry}.

Observe that an analytic solution $\psi$ of \eqref{eq:main-cos}
restricts to a solution $u_n=\psi(nh)$ of the discrete equation
\eqref{eq:finite-mathieu-operator} on the lattice $h\Z$.
If this sequence is nonzero and $N$-periodic, it gives an eigenvector
of $H^N_{0,0,\lambda}$.

Conversely, for sufficiently small $h$, the WKB construction below
provides two analytic solutions from which we can form a linear
combination agreeing with any prescribed discrete solution at two
consecutive lattice points.\footnote{The required condition is that
the discrete Wronskian of the two analytic solutions be nonzero at
a lattice point. This makes the two equations determining the
coefficients of the linear combination uniquely solvable. The WKB
estimates give this nonvanishing on a regular domain; see
\Cref{thm:uniform-wkb}.}
The recurrence then gives agreement on the whole lattice.
Thus every periodic discrete solution extends to an analytic
solution of \eqref{eq:main-cos}, although this extension need not
be periodic away from the lattice.

\subsection{WKB solutions and canonical domains}
\label{sec:wkb-method-second}

We now describe the asymptotic behavior of solutions to \eqref{eq:main-cos}.
We first define the complex momentum as a local branch $p(z)$ satisfying
$\cos p(z)=E-\lambda\cos z$.
\Cref{thm:uniform-wkb} gives solutions $\psi_\pm$ to
\eqref{eq:main-cos} with the following asymptotics:
\begin{align}
 \psi_\pm(z)
 &=\frac{1}{\sqrt{\sin p(z)}}\,
   \exp\left(\pm\frac{i}{h}\int_{z_0}^z p(w)\dd w\right)(1+O(h)).
 \label{eq:wkb-form}
\end{align}
These asymptotics are valid in so-called \textit{canonical domains} (\Cref{def:canonical_domain}) whose size is independent of $h$.

\subsubsection{The geometric condition}
\label{sec:canonical-domains}

A turning point of the complex momentum $p$ is where $\sin p = 0$.
A region is \emph{regular} if it contains no turning point.
A regular curve is \emph{vertical} if its tangent is never parallel
to the abscissa axis. We parameterize it as $z=x(y)+iy$.
A \emph{straight vertical} line is parallel to the ordinate axis.
\begin{defi}[Canonical curve]
\label{def:2}
A regular vertical curve $z=x(y)+iy$ is canonical for the chosen branch of $p$ if
\begin{align}
 0<\Re p(x(y)+iy)+\Im p(x(y)+iy)\,x'(y)<\pi.
 \label{eq:canonical-curve-criterion}
\end{align}
\end{defi}
The lower inequality says that $\Im\int^z p(w)\dd w$ increases
along the curve; the upper says that
$\Im\int^z(p(w)-\pi)\dd w$ decreases.
These are the canonical-curve conditions of Fedotov and
Klopp \cite[Section~3.3.1]{MR2097997}; for difference equations, see
Fedotov and Shchetka \cite[(2.10)]{FedotovShchetka2022Berry}.
For a straight vertical line, the condition reduces to $0<\Re p<\pi$.

\begin{defi}[Canonical domain]\label{def:canonical_domain}
A {canonical domain} is a regular domain whose horizontal
slices are intervals and is the union of canonical curves.
A bounded canonical domain is composed of canonical curves all joining two common boundary points.
An unbounded canonical domain is composed of canonical curves extending to both imaginary infinities.
\end{defi}

\Cref{thm:uniform-wkb} gives solutions to \eqref{eq:main-cos}
of the form \eqref{eq:wkb-form} on each fixed subdomain
$\widetilde U\Subset U$ of a bounded canonical domain $U$,
uniformly for nearby energies, provided its hypotheses hold
uniformly as $E$ varies.

Once a WKB estimate is established on an open neighborhood,
it extends along a fixed regular horizontal segment starting
in that neighborhood, in the direction in which its exponential
grows, provided $\Im p$ has a strict sign.
Thus the plus estimate propagates rightward when $\Im p<0$ and
leftward when $\Im p>0$; the directions reverse for the minus estimate.
This is the horizontal continuation principle, or rectangle lemma,
stated in \Cref{thm:finite-segment-continuation}.

\subsubsection{One initial strip}

Define $f_E(z)\coloneq E-\lambda\cos z$ so that $\cos(p)=f_E(z)$.  For $\Im\lambda>0$,
\begin{align}
 \Im f_E(\pi+iy)&=\Im E+\Im\lambda\cosh y>0
 \qquad\text{if }\Im E>-\Im\lambda.
 \label{eq:straight-canonical-line}
\end{align}
We fix the branch of the complex momentum such that $0<\Re p<\pi$ and $\Im p<0$.
Thus the line $\Re z=\pi$ is canonical at every interior energy.

\begin{prop}
\label{prop:unbounded-canonical-strip}
Let $\Im\lambda>0$ and let $\mathcal K$ be a compact subset of
$\{\Im E>-\Im\lambda\}$.  There are $\delta>0$ and a neighborhood
$U$ of $\mathcal K$ such that
\begin{align}
 K_\delta&\coloneq\{\pi+x+iy:\lvert x\rvert <\delta,\ y\in\R\}
 \label{eq:unbounded-pi-strip}
\end{align}
is regular for every $E\in U$.  On this strip the branch has
$0<\Re p(z;E)<\pi$, $\Im p(z;E)<0$, and every straight vertical line
is canonical.
\end{prop}
\begin{proof}
Choose $\delta$ so small that
\begin{align}
 c_\delta=(\Im\lambda)\cos\delta-\lvert \Re\lambda\rvert \sin\delta
 >\max\bigl(0,-\min_{E\in\mathcal K}\Im E\bigr).
\end{align}
Shrink a bounded neighborhood $U$ so that the same inequality holds
with $e_*=\inf_{E\in U}\Im E$.  For $\lvert x\rvert \le\delta$,
\begin{align}
 \Im f_E(\pi+x+iy)
 &\ge e_*+c_\delta\cosh y
 \ge\bigl(c_\delta+\min(e_*,0)\bigr)\cosh y>0.
\end{align}
Since $f_E$ takes values in the upper half-plane, we may choose
the analytic inverse-cosine branch satisfying
$0<\Re p(z;E)<\pi$ and $\Im p(z;E)<0$ throughout $K_\delta$.
In particular, $f_E\ne\pm1$, so the strip contains no turning points.
Moreover, $0<\Re p(z;E)<\pi$ is precisely the canonical condition
\eqref{eq:canonical-curve-criterion} for a straight vertical line.
\end{proof}

\begin{prop}[Strip-normalized entire WKB solution]
\label{prop:entire-WKB-normalization}
On every fixed narrower strip $K_{\delta'}$, $0<\delta'<\delta$,
there is an exact solution
\begin{align}
 \phi_+^0(z;E,\lambda,h)
 &=\frac{1}{\sqrt{\sin p(z;E)}}\,
   \exp\left(\frac{i}{h}S_E(z)\right)(1+O(h)),
 &
 S_E(z)&\coloneq\int_\pi^z p(w;E)\dd w,
 \label{eq:entire-strip-WKB}
\end{align}
uniformly for all $z\in K_{\delta'}$ and $E\in\mathcal K$.
It is entire in $z$ and analytic in $E$ near $\mathcal K$.
Its reflection $\phi_-^0(z)\coloneq\phi_+^0(2\pi-z)$ has the minus
WKB estimate on the same strip.
\end{prop}

The unbounded WKB estimates of
\cite[Theorems~2.2--2.3]{FedotovShchetka2022Berry}
give uniform control as $\Im z\to\pm\infty$ within the strip.
The recurrence relation from \eqref{eq:main-cos} then extends $\phi_+^0$ horizontally to the entire
plane.
The detailed verification is in Appendix~\ref{app:growth}.
Multiplication of $\phi_+^0$ and $\phi_-^0$ by a common scalar
$1+O(h)$ preserves these estimates; we fix that scalar in the exact
trace reduction below.

\subsection{The one-period Wronskian matrix}
\label{sec:unified-period-map}

For analytic solutions $\phi,\psi$ of \eqref{eq:main-cos}, the Wronskian is defined by
\begin{align}
 W[\phi,\psi](z)
 &\coloneq\phi(z+h)\psi(z)-\phi(z)\psi(z+h).
 \label{eq:discrete-wronskian}
\end{align}
The recurrence relation in \eqref{eq:main-cos} makes the Wronskian $h$-periodic.  
If $W[\phi,\psi]$ does not vanish, we say $\phi$ and $\psi$ are linearly independent and any other 
solution $\xi$ satisfies
\begin{align}
 W[\phi,\psi]\,\xi
 &=W[\xi,\psi]\,\phi+W[\phi,\xi]\,\psi.
 \label{eq:wronskian-expansion}
\end{align}

Let $\phi_+^0,\phi_-^0$ be linearly independent solutions to \eqref{eq:main-cos}, and define their translates as
$\phi_\pm^{2\pi}(z)\coloneq\phi_\pm^0(z-2\pi)$.
Because $\phi_+^0,\phi_-^0$ are linearly independent, for every
$z_0\in\C$ there exists a unique $2\times2$ matrix $T(E,h;z_0)$
such that on the lattice $z_0+h\Z$,
\begin{align}
 (\phi_+^{2\pi},\phi_-^{2\pi})
 &=(\phi_+^0,\phi_-^0)\,T(E,h;z_0).
 \label{eq:unified-connection-matrix}
\end{align}
The Wronskian identity \eqref{eq:wronskian-expansion} gives
\begin{align}
 T(E,h;z_0)
 =\frac1{W[\phi_+^0,\phi_-^0](z_0)}
 \begin{pmatrix}
  W[\phi_+^{2\pi},\phi_-^0]&
  W[\phi_-^{2\pi},\phi_-^0]\\
  W[\phi_+^0,\phi_+^{2\pi}]&
  W[\phi_+^0,\phi_-^{2\pi}]
 \end{pmatrix}_{z=z_0}.
 \label{eq:unified-wronskian-matrix}
\end{align}
The matrix $T$ represents the backward shift by one period
in the solution basis $(\phi_+^0,\phi_-^0)$, so it is similar
to $M_N(0,z_0,E,\lambda)^{-1}$ (defined in \eqref{eq:finite-transfer-matrix}).
Since $M_N$ is a $2\times2$ matrix with determinant one,
$\tr T=\tr(M_N^{-1})=\tr M_N$.
Thus, by \eqref{eq:trace-polynomial-decomposition},
\begin{align}
 \tr T(E,h;z_0)
 &=m_N(E,\lambda)-2\lambda^N\cos(Nz_0).
 \label{eq:unified-wronskian-trace}
\end{align}
In particular, choosing $z_0=h/4$ gives
\begin{align}
 \tr T(E,h;h/4)&=m_N(E,\lambda).
 \label{eq:quarter-phase-trace}
\end{align}
Equations \eqref{eq:unified-wronskian-matrix}--\eqref{eq:quarter-phase-trace}
are the common starting point for the interior, exterior, and
regular-boundary calculations.  Only the bases used to estimate the four
Wronskians change.  The corner argument uses the same trace polynomial
through the exact coefficient identities derived next.

\subsection{Exact reduction of the trace}
\label{sec:exact-reduction}

We now specialize the Wronskian formula of
\Cref{sec:unified-period-map} to the solutions normalized on the
unbounded strip, and set
\begin{align}
 \mathcal{A}&\coloneq W[\phi_-^0,\phi_+^0],
 &\mathcal{B}&\coloneq W[\phi_-^0,\phi_+^{2\pi}],\\
 \mathcal{C}&\coloneq W[\phi_+^0,\phi_+^{2\pi}],
 &\mathcal{D}&\coloneq W[\phi_-^{2\pi},\phi_+^0].
 \label{eq:four-wronskians}
\end{align}
We write $\mathcal C(z)=\mathcal C(z,E)$, displaying the energy
argument only when needed.
Recall that $\phi_-^0(z)=\phi_+^0(2\pi-z)$, so that
\begin{align}
 W[\phi_-^0,\phi_-^{2\pi}](z)
 &=\mathcal C(4\pi-z-h)=\mathcal C(-z),
\end{align}
where the last equality uses the $h$-periodicity of $\mathcal C$
and $2\pi=Nh$.
Thus
\begin{align}\label{eq:T_components}
  T(E,h;z_0) = \frac{1}{  \mathcal{A}(z_0)}\begin{pmatrix}
  \mathcal{B}(z_0) & \mathcal{C}(-z_0) \\
  -\mathcal{C}(z_0) &  \mathcal{D}(z_0)
 \end{pmatrix}.
\end{align}

Our main objective is 
to estimate
\begin{align}
 m_N(E,\lambda)
 &=\tr T(E,h;h/4)= \mathcal{A}(h/4)^{-1}(\mathcal{B}(h/4)+\mathcal{D}(h/4)).
\end{align}

Two exact structural facts make this reduction useful.  The first fixes
the normalization of the global basis.

\begin{prop}
\label{prop:constant-AB}
The functions $\mathcal{A}$ and $\mathcal{B}$ are independent of $z$, and
$\mathcal{A}=-2i(1+O(h))$.  Multiplying $\phi_+^0$ and $\phi_-^0$ by a common
scalar $1+O(h)$ makes $\mathcal{A}=-2i$ exactly.
\end{prop}

Use this normalization from now on, and put
\begin{align}
 G(z)&\coloneq e^{2\pi iz/h},& b&\coloneq \mathcal{B}/\mathcal{A}.
 \label{eq:scalar-coefficient-notation}
\end{align}

The second fact is the exact Fourier truncation that makes the later WKB
calculation finite.

\begin{prop}
\label{prop:exact-C-modes}
For sufficiently small $h=2\pi/N$,
\begin{align}
 \mathcal{C}(z,E)&=c(E,h)+a(E,h)G(z),
 \label{eq:exact-C-two-modes}
\end{align}
where $a,c$ are analytic in $E$.
\end{prop}

\begin{prop}
\label{cor:entire-normalization-unique}
Any two exact solutions with the uniform plus estimate
\eqref{eq:entire-strip-WKB} differ by a scalar $1+O(h)$.
If both have $\mathcal{A}=-2i$, they coincide.
\end{prop}

The proofs are in Appendix~\ref{app:wronskian-growth}.  The key observation is simple:
a nonconstant $h$-periodic (in the real direction) Fourier term grows at least as
$e^{2\pi\lvert \Im z\rvert /h}$ at one imaginary end.  Bounds on the solutions
exclude those terms in $\mathcal{A},\mathcal{B}$ and all but two terms in $\mathcal{C}$.
The same bounds give uniqueness.

\begin{prop}[Exact trace identity]
\label{prop:exact-two-mode-reduction}
The coefficients $a$, $c$ (from \Cref{prop:exact-C-modes}), and $b$ (from \eqref{eq:scalar-coefficient-notation}) satisfy
\begin{align}
 ac&=\mathcal{A}\mathcal{B}\lambda^N,&
 \mathcal{A}\mathcal{B}m_N&=\mathcal{A}^2+\mathcal{B}^2-a^2-c^2.
 \label{eq:division-free-reduction}
\end{align}
If $\mathcal{B}\ne0$, then $a,c\ne0$ and
\begin{align}
 m_N&=b+b^{-1}-\lambda^N(a/c+c/a).
 \label{eq:normalized-exact-reduction}
\end{align}
In this normalization, $ac=-4b\lambda^N$.
\end{prop}
\begin{proof}
Because $\det T = 1$, we have (by \eqref{eq:T_components}),
\begin{align}
\mathcal{B} \mathcal{D}(z) + \mathcal{C}(z)\mathcal{C}(-z) = \mathcal{A}^2.\label{eq:1493}
\end{align}
By \Cref{prop:exact-C-modes}, $\mathcal{C}(z)\mathcal{C}(-z)=a^2+c^2+ac(G+G^{-1})$.
Next, the trace \eqref{eq:unified-wronskian-trace} can be rewritten as
\begin{align}
  \frac{\mathcal{B} + \mathcal{D}(z )}{\mathcal{A}} = m_N - \lambda^N (G(z) + G(z)^{-1}),\label{eq:1497}
\end{align}
so that
$\mathcal{D}=\mathcal{A}m_N-\mathcal{B}-\mathcal{A}\lambda^N(G+G^{-1})$.
Substituting this expression into \eqref{eq:1493}, and using the expression for $\mathcal{C}(z)\mathcal{C}(-z)$, we obtain
\begin{align}
\mathcal{B} (\mathcal{A}m_N-\mathcal{B}-\mathcal{A}\lambda^N(G+G^{-1})) + a^2+c^2+ac(G+G^{-1}) = \mathcal{A}^2.
\end{align}    
Comparing constant and first Fourier coefficients gives the two identities in \eqref{eq:division-free-reduction}.  When $\mathcal{B}\ne0$, we can divide by $\mathcal{A}\mathcal{B}$ and use $ac=\mathcal{A}\mathcal{B}\lambda^N$ to obtain \eqref{eq:normalized-exact-reduction}.   
\end{proof}

The operator duality also determines the individual coefficients.
Unlike trace duality alone, the following identity will give the
second interior estimate from the first.

\begin{prop}[Aubry--Andr\'e duality for the coefficients]
\label{prop:coefficient-duality}
Fix $\Im\lambda>0$ and a compact set
$\mathcal K\Subset\operatorname{int}\Sigma_\lambda$, and put
\begin{align}
 E^\dagger&\coloneq\overline{E/\lambda},&
 \lambda^\dagger&\coloneq1/\overline\lambda.
 \label{eq:coefficient-dual-parameters}
\end{align}
For all sufficiently large $N$, the normalized coefficients satisfy
\begin{align}
 a(E,\lambda)+c(E,\lambda)
       \overline{b(E^\dagger,\lambda^\dagger)}&=0
 \label{eq:coefficient-duality}
\end{align}
on a neighborhood of $\mathcal K$.
\end{prop}
\begin{proof}
The transformed coupling has positive imaginary part. Moreover,
$E=s+\lambda t$, $-1<s,t<1$, gives
$E^\dagger=t+\lambda^\dagger s$, so the dual energy is also interior.
Let
\begin{align}
 \widehat\phi_+^0(\rho)
 &\coloneq\overline{\phi_+^0(\overline\rho;E^\dagger,\lambda^\dagger)},&
 \widehat\phi_-^0(\rho)&\coloneq\widehat\phi_+^0(2\pi-\rho).
\end{align}
Here we display the energy and coupling and suppress the unchanged
parameter $h$.
\Cref{lem:fourier-transform,lem:fourier-normalization} construct a
Fourier transform $\mathcal F_h$, with kernel $e^{-i\rho z/h}$, for which
\begin{align}
 \mathcal F_h\phi_+^0&=\kappa(E,h)\widehat\phi_+^0,
 &\kappa(E,h)&\ne0.
\end{align}
The transform respects translation and reflection, including when
the defining integrals initially converge on different strips.

By \eqref{eq:T_components}, the first connection relation is
\begin{align}
 \phi_+^0(z-2\pi)
 &=b\phi_+^0(z)+\frac{c+a e^{2\pi iz/h}}{2 i }\phi_-^0(z). \label{eq:exact-connection-relation2}
\end{align}

We then apply $\mathcal F_h$ to both sides.
By \Cref{lem:fourier-transform},
\begin{align}
 \mathcal F_h[\phi_+^0(\,\cdot-2\pi)](\rho)
&=\kappa e^{-2\pi i\rho/h}\widehat\phi_+^0(\rho), \qquad \mathcal F_h[\phi_-^0](\rho)
 =\kappa e^{-2\pi i\rho/h}\widehat\phi_+^0(-\rho),\\
 \mathcal F_h[e^{2\pi i(\cdot)/h}\phi_-^0](\rho)
 &=(\mathcal F_h\phi_-^0)(\rho-2\pi)\\
 &=\kappa e^{-2\pi i(\rho-2\pi)/h}
       \widehat\phi_+^0(2\pi-\rho)\\
 &=\kappa e^{-2\pi i\rho/h}\widehat\phi_-^0(\rho).
\end{align}
Here the last equality uses
$\widehat\phi_-^0(\rho)=\widehat\phi_+^0(2\pi-\rho)$ and
$e^{4\pi^2i/h}=e^{2\pi iN}=1$.

Substituting these identities into the Fourier transform of \eqref{eq:exact-connection-relation2} and cancelling the nonzero factor
$\kappa e^{-2\pi i\rho/h}$ gives
\begin{align}
 \widehat\phi_+^0(\rho)
 &=b e^{2\pi i\rho/h}\widehat\phi_+^0(\rho)
   +\frac{c\widehat\phi_+^0(-\rho)
          +a\widehat\phi_-^0(\rho)}{2i}.
\end{align}
Rearranging, we obtain
\begin{align}
 c\widehat\phi_+^0(-\rho)+a\widehat\phi_-^0(\rho)
 &=2i\bigl(1-b e^{2\pi i\rho/h}\bigr)
       \widehat\phi_+^0(\rho).
 \label{eq:fourier-connection-comparison}
\end{align}
Applying $W[\widehat\phi_+^0,\cdot]$ to both sides (using that $W[\widehat\phi_+^0,\widehat\phi_-^0]=-2i$ and $e^{2\pi i\rho/h}$ is $h$-periodic) gives
\begin{align}
cW[\widehat\phi_+^0,\widehat\phi_+^0(-\cdot)] -2i a =0 .\label{eq:1598}
\end{align}

On the other hand, replacing $z$ by $2\pi-z$ in
\eqref{eq:exact-connection-relation2} gives
\begin{align}
 \phi_+^0(-z)
 &=b\phi_+^0(2\pi-z)
   +\frac{\mathcal C(2\pi-z)}{2i}\phi_-^0(2\pi-z)\\
 &=b\phi_-^0(z)
   +\frac{\mathcal C(2\pi-z)}{2i}\phi_+^0(z).
\end{align}
Applying this identity at $(E^\dagger,\lambda^\dagger)$ and
conjugating, we obtain
\begin{align}
 \widehat\phi_+^0(-\rho)
 &=\overline{b(E^\dagger,\lambda^\dagger)}
       \widehat\phi_-^0(\rho)
   +\beta(\rho)\widehat\phi_+^0(\rho),
\end{align}
where $\beta$ is $h$-periodic because $\mathcal C$ is.

Applying $W[\widehat\phi_+^0,\cdot]$ to both sides gives
\begin{align}
 W[\widehat\phi_+^0,\widehat\phi_+^0(-\cdot)]
 &=-2i\,\overline{b(E^\dagger,\lambda^\dagger)}.
\end{align}
Substituting this into \eqref{eq:1598} gives
\begin{align}
 -2i\left(c\,\overline{b(E^\dagger,\lambda^\dagger)}+a\right)
 &=0.
\end{align}
\end{proof}

On a simply connected energy neighborhood where $\mathcal{B}\ne0$, we define the following analytic functions
\begin{align}
 P_h&\coloneq\frac{ih}{2\pi}\log b,&
 Z_h&\coloneq-\frac{ih}{2\pi}\log(-a/c).
 \label{eq:exact-quantum-actions}
\end{align}
Substitution into \eqref{eq:normalized-exact-reduction} gives
\begin{align}
 m_N
 &=2\cos(2\pi P_h/h)+2\lambda^N\cos(2\pi Z_h/h).
 \label{eq:exact-action-factorization}
\end{align}
Changing the logarithm branches shifts $P_h$ and $Z_h$ by
integer multiples of $h$, leaving
\eqref{eq:exact-action-factorization} unchanged.
In the interior, the WKB estimates below give
$P_h=P+O(h^2)$ and $Z_h=Z+O(h^2)$, where $P,Z$ are the
classical actions. Across an open side, a horizontal WKB comparison
gives a two-action formula with locally constructed corrections of
the same order. See \Cref{thm:interior-side-quantization}.

\section{Quantization away from the corners}
\label{sec:away-from-corners}

We now estimate the trace polynomial $m_N$ using the exact identities of the
preceding section. For interior energies, we continue the WKB estimate
of $\phi_+^0$ from the normalization strip to a neighborhood of $z=0$.
This allows us to evaluate $\mathcal B$ and hence estimate $b$.
Applying the same calculation to the dual equation then gives $a/c$
by \Cref{prop:coefficient-duality}.
For open-side and exterior energies, we instead compare WKB solutions
across one period to estimate the trace directly.
We first construct the required comparison regions, then carry out
these calculations and express the resulting contour integrals in
terms of the classical actions.

\subsection{Momentum geometry}
\label{sec:momentum-geometry}
Assume $\Im\lambda>0$ and $E\in\operatorname{int}\Sigma_\lambda$.
Recall the complex momentum $p(z)$ satisfies
\begin{align}
 \cos p(z)&=E-\lambda\cos z.
 \label{eq:cosine-momentum}
\end{align}
For real $z$, we choose a branch with $0<\Re p<\pi$.
Indeed, the segment $E-\lambda[-1,1]$ meets the real axis
at exactly one point in $(-1,1)$, so it avoids the cuts
$(-\infty,-1]\cup[1,\infty)$ of this inverse-cosine branch.
For each real $z$, this choice is analytic in $E$ throughout
$\operatorname{int}\Sigma_\lambda$.

Thus, there is a unique $z_*\in(0,\pi)$ such that
\begin{align}
 E-\lambda\cos z_*&\in(-1,1).
\end{align}
Put $p_*\coloneq p(z_*)\in(0,\pi)$.  The quantity
$E-\lambda\cos z$ can be real for real $z\in(0,\pi)$ only at $z_*$,
where it equals $\cos p_*$.  
The same argument, with $p,z$ interchanged, fixes the inverse branch
$z(p)=\arccos((E-\cos p)/\lambda)$ on the real $p$-axis.

The branch points of $p$, called \emph{turning points}, occur where
$\cos p=\pm1$.  In one period they are
$\zeta_+,\zeta_-,2\pi-\zeta_+,2\pi-\zeta_-$, where
\begin{align}
 \cos\zeta_+&=(E-1)/\lambda,&
 \cos\zeta_-&=(E+1)/\lambda,&
 0<\Re\zeta_\pm&<\pi.
 \label{eq:turning-points}
\end{align}
The curve
\begin{align}
 z(u)&=\arccos\left(\frac{E-u}{\lambda}\right),
 \qquad -1\le u\le1,
 \label{eq:action-arc}
\end{align}
is chosen by requiring $z(\cos p_*)=z_*$.  It joins $\zeta_-$ to
$\zeta_+$ and has real momentum $p=\arccos u$.  We call it the
 \emph{real-momentum arc}.  Its interior is regular and its endpoints are
two simple turning points.

The real-momentum arc, its reflected and periodic copies, the continuations
of the locus $\cos p\in\mathbb R$ beyond its endpoints, and the initial
strip $K_\delta$ are shown in
\Cref{fig:real-momentum-strip}.

\realmomentumstripfigure

At an exterior energy, the segment $E-\lambda[-1,1]$ is disjoint from
$[-1,1]$.  On a sufficiently small simply connected energy
neighborhood there is therefore a $2\pi$-periodic branch $p$, analytic
near the real $z$-axis, such that
\begin{align}
 \cos p(z)&=E-\lambda\cos z,&
 \Im p(z)&<0.
 \label{eq:exterior-momentum}
\end{align}
The path $z\mapsto E-\lambda\cos z$ traverses the segment and retraces
it, so the chosen inverse-cosine branch returns to its initial value.

\subsection{Comparison regions}
\label{sec:interior-comparison}

\subsubsection{The adjacent canonical-domain principle}
\label{sec:same-solution-enlargement}

Extending the exact solution and extending its WKB estimate are
different steps.  The recurrence relation from \eqref{eq:main-cos} gives the first.  For the second,
we either continue in a growing direction (rectangle lemma, see \Cref{thm:finite-segment-continuation}) or enlarge the canonical
domain while keeping the original infinite strip (\Cref{lem:same-solution-enlargement}).

The following is the finite-difference analogue of the adjacent canonical-domain
principle of Fedotov and Klopp \cite[Proposition~5.1]{MR2097997}.

\begin{lemma}[Adjacent canonical-domain principle]
\label{lem:same-solution-enlargement}
Let $K$ be an unbounded canonical domain containing the strip
$K_{\delta'}$ on which $\phi_+^0$ was constructed, with the same
momentum branch and phase based at $\pi$.
Suppose the unbounded WKB results
\cite[Theorems~2.2--2.3]{FedotovShchetka2022Berry}
apply on $K$ with uniform parameter bounds.
Then the solution $\phi_+^0$ already chosen has its
plus WKB estimate on every compact subset of $K$.
By reflection, $\phi_-^0(z)=\phi_+^0(2\pi-z)$ has its minus WKB
estimate on every compact subset of $2\pi-K$, with the reflected
momentum branch.
\end{lemma}
\begin{proof}
Construct a plus solution on $K$ and extend it by the recurrence.
On the original strip it has the uniform estimate at both imaginary
ends.  By \Cref{cor:entire-normalization-unique}, it is a scalar
$1+O(h)$ times $\phi_+^0$.  Normalizing its reflected Wronskian to
$-2i$ makes it exactly $\phi_+^0$, without changing the estimate on
$K$.
\end{proof}

\subsubsection{WKB comparison regions for the Wronskians and trace}
\label{sec:general-interior-domains}

In the interior we continue the globally normalized solution
$\phi_+^0$ to $z=0$ to estimate $\mathcal B$. Across an open side and
in the exterior we use local bases to estimate the trace.

\begin{prop}[WKB comparison regions]
\label{prop:uniform-WKB-comparison-geometry}
Assume that $\Im\lambda>0$.  Each construction below is uniform on a
complex neighborhood of the indicated energy $E_*$.

\begin{enumerate}
\item[\textup{(i)}] \emph{Interior comparisons.}
At every $E_*\in\operatorname{int}\Sigma_\lambda$,
the WKB description of $\phi_+^0$ extends to a fixed neighborhood of
$z=0$.

\item[\textup{(ii)}] \emph{Exterior and open-side comparisons.}
Suppose either $E_*\notin\Sigma_\lambda$ or
$E_*=s+\lambda$ with $-1<s<1$.  There are a fixed $\delta>0$,
local exact WKB bases $\psi_\pm^0$ near $[-2\delta,2\delta]$, and
their translates
$\psi_\pm^{2\pi}(z)=\psi_\pm^0(z-2\pi)$ near $2\pi$.
The momentum is even and periodic near the real period, with
\begin{align}
 \Im p(x)&<0 &&(\delta\leq x\leq2\pi-\delta),&
 \Im\int_0^\delta p(x)\dd x&<0.
 \label{eq:buffered-horizontal-signs}
\end{align}
The plus estimate for $\psi_+^0$ continues rightward, and the minus
estimate for $\psi_-^{2\pi}$ continues leftward.
In particular, there are fixed overlaps near $\delta$ and
$2\pi-\delta$ on which the corresponding local and continued estimates
hold. 
In the exterior, use the branch $p$ defined by \eqref{eq:exterior-momentum}. Near the open side $E_*=s+\lambda$, use the analytic continuation in $E$ of $p$ from the interior.

\end{enumerate}

All comparison regions are independent of sufficiently small $h$ and
remain a fixed positive distance from the turning points.  The
canonical inequalities, and the growing-direction inequalities on the
middle interval in \textup{(ii)}, have uniform strict margins.
\end{prop}

\begin{proof}[Sketch of proof]

\textbf{Interior energies.}
Choose a regular point $z_\vartheta$ of the real-momentum arc
\eqref{eq:action-arc} with $p(z_\vartheta)=\vartheta\in(0,\pi)$,
and set
\begin{align}
 H_\vartheta(z)&\coloneq
 \Im\int_{z_\vartheta}^z(p(w)-\vartheta)\dd w.
\end{align}
Along a regular vertical level curve $z=x(y)+iy$,
\begin{align}
 0&=\frac{\dd}{\dd y}H_\vartheta(x(y)+iy)
   =\Im p\,x'(y)+\Re p-\vartheta,\\
 &\Rightarrow \Re p+\Im p\,x'(y)=\vartheta\in(0,\pi).
\end{align}
Thus these curves satisfy the canonical inequalities
\eqref{eq:canonical-curve-criterion} with a strict margin.
Small negative level curves pass just to the right of $z_\vartheta$.
Since $H_\vartheta\geq0$ on the real-momentum arc, they cannot
cross the arc or meet its endpoints.
Each end of a level curve either reaches $K_{\delta'}$ or extends to imaginary infinity. In the latter case, we bend it into $K_{\delta'}$ while preserving the canonical inequalities.
The resulting family of curves fills an
unbounded canonical domain containing the normalization strip.

Now take $z_\vartheta=z_*$. A small local deformation carries
part of this domain across the arc near $z_*$, while preserving
the strict canonical inequalities and leaving the normalization
strip unchanged. The adjacent canonical-domain principle
(\Cref{lem:same-solution-enlargement}) identifies its plus solution
with $\phi_+^0$. Immediately to the left of the crossing,
$\Im p>0$; starting there, the rectangle lemma
(\Cref{thm:finite-segment-continuation}) carries the WKB estimate
leftward to zero.

\Cref{fig:interior-domain-construction} illustrates the canonical
enlargement and horizontal continuation used to reach zero.

\interiorconstructionfigure

\textbf{Exterior and open-side energies.}
For part \textup{(ii)}, continue the plus estimate from the left endpoint
rightward and the minus estimate from the right endpoint leftward,
using $\Im p<0$ on the middle interval; see
\Cref{fig:exterior-wkb-comparison}.
At an open-side energy $E=s+\lambda$, $-1<s<1$, we have
$p(0)=\arccos s\in(0,\pi)$, so zero is regular despite
$\Im p(0)=0$.
Horizontal continuation requires $\Im p(x)<0$ only on
$[\delta,2\pi-\delta]$.
The other sides will be obtained by symmetry and duality at the level
of the trace.

\exteriorwkbcomparisonfigure

\end{proof}

The complete constructions, including the bends at infinity and
uniformity for nearby energies, are proved in Appendix~\ref{app:geometry}.

\subsection{Wronskian calculations}
\label{sec:raw-Wronskian-calculations}

\subsubsection{Interior energies}

\begin{prop}
\label{prop:general-interior-B}
For $\Im\lambda>0$, locally near every interior energy,
\begin{align}
 \mathcal{B}
 &=-2i\exp\left(-\frac{i}{h}\int_0^{2\pi}p(z)\dd z\right)
   (1+O(h)).
 \label{eq:general-interior-B}
\end{align}
The error is analytic in $E$ and locally uniform.  In
particular, $\mathcal{B}\ne0$ for sufficiently small $h$.
\end{prop}
\begin{proof}
By \eqref{eq:discrete-wronskian} and \eqref{eq:four-wronskians},
\begin{align}
 \mathcal{B}
 &=W[\phi_-^0,\phi_+^{2\pi}](2\pi)\\
 &=\phi_-^0(2\pi+h)\phi_+^{2\pi}(2\pi)
   -\phi_-^0(2\pi)\phi_+^{2\pi}(2\pi+h)\\
 &=\phi_+^0(0)\bigl(\phi_+^0(-h)-\phi_+^0(h)\bigr).
 \label{eq:1388}
\end{align}
The last equality uses $\phi_-^0(z)=\phi_+^0(2\pi-z)$ and
$\phi_+^{2\pi}(z)=\phi_+^0(z-2\pi)$.
By \Cref{prop:uniform-WKB-comparison-geometry}\textup{(i)},
the WKB estimate of $\phi_+^0$ extends to a fixed neighborhood
of $z=0$.
Taylor expansion at $z=0$ gives
\begin{align}
 \phi_+^0(\pm h)
 &=\frac{1}{\sqrt{\sin p(\pm h)}}
   \exp\left(\frac{i}{h}\int_\pi^{\pm h}p(w)\dd w\right)
   (1+O(h))\\
 &=\frac{e^{\pm ip(0)}}{\sqrt{\sin p(0)}}
   \exp\left(\frac{i}{h}\int_\pi^0p(w)\dd w\right)
   (1+O(h)).
\end{align}
Substitution into \eqref{eq:1388} yields
\begin{align}
 \mathcal{B}
 &=\frac{\exp\left(\frac{2i}{h}\int_\pi^0p(w)\dd w\right)}
         {\sin p(0)}
   \left[e^{-ip(0)}(1+O(h))-e^{ip(0)}(1+O(h))\right]\\
 &=-2i\exp\left(\frac{2i}{h}\int_\pi^0p(w)\dd w\right)(1+O(h)).
\end{align}
Here $\sin p(0)$ is bounded away from zero on compact interior
energy sets: its vanishing would require $E=\lambda\pm1$, which
are corners.

Finally, $p(2\pi-x)=p(x)$ gives
$2\int_\pi^0p=-\int_0^{2\pi}p$, proving the formula.
The error is analytic in $E$, since both $\mathcal B$
and its nonvanishing leading exponential depend analytically on $E$.
Its $O(h)$ bound is locally uniform by the uniform WKB estimates.
\end{proof}

\subsubsection{A horizontal comparison across an open side and in the exterior}
\label{sec:boundary-transfer}
\label{sec:exterior-wronskian}

We next estimate $m_N$ across the open side $E-\lambda\in(-1,1)$
and in the exterior using the Wronskian matrix
\eqref{eq:unified-wronskian-matrix} in a local basis.

Use the branch of $p$ and endpoint bases $\psi_{\pm}^{0,2\pi}$ in part \textup{(ii)} of
\Cref{prop:uniform-WKB-comparison-geometry}.  Recall
the normalized action
\begin{align}
 Q(E)&\coloneq\frac{1}{2\pi}\int_0^{2\pi}p(x)\dd x.
 \label{eq:exterior-action}
\end{align}
Write the change of basis matrix $T$ as
\begin{align}
 (\psi_+^{2\pi},\psi_-^{2\pi})
 &=(\psi_+^0,\psi_-^0)T.
\end{align}
We evaluate $T$ on the lattice $h/4+h\mathbb Z$, so its trace
is exactly $m_N$ by \eqref{eq:quarter-phase-trace}.  The action $Q$ has
strictly negative imaginary part, and hence $\lvert e^{2\pi iQ/h}\rvert \geq e^{c/h}$ for
some $c>0$ on a smaller energy neighborhood.

\begin{prop}
\label{prop:exterior-Wronskian-estimates}
Locally near every exterior energy and every point of
$E-\lambda\in(-1,1)$ with $\Im\lambda>0$,
\begin{align}
 T_{22}&=e^{2\pi iQ/h}(1+O(h)),\\
 T_{12},T_{21}&=O\!\left(he^{2\pi iQ/h}\right),
 \label{eq:exterior-Wronskian-entries}\\
 T_{11}&=e^{-2\pi iQ/h}(1+O(h))+O\!\left(h^2e^{2\pi iQ/h}\right).
 \label{eq:exterior-fourth-entry}
\end{align}
Consequently $m_N=e^{2\pi iQ/h}(1+O(h))$.  The errors are analytic in $E$ and
locally uniform.
\end{prop}
\begin{proof}
Choose $x_L,x_R\in h/4+h\Z$ in fixed compact intervals
contained in the overlaps near $\delta$ and $2\pi-\delta$,
respectively, supplied by
\Cref{prop:uniform-WKB-comparison-geometry}\textup{(ii)}.
Choose them so that $x_L+h$ and $x_R+h$ remain in the
respective intervals.
Such points exist for all small $h$.  Since the
connection coefficients are $h$-periodic, comparisons at these two
points estimate the same matrix.

Normalize the local basis by
\begin{align}
 \psi_\pm^0(z)
 &=\frac{1}{\sqrt{\sin p(z)}}
   \exp\left(\pm\frac{i}{h}\int_0^z p(w)\dd w\right)(1+O(h)).
 \label{eq:exterior-left-WKB-basis}
\end{align}
Translation gives the normalization at $2\pi$.

\textbf{The dominant entry.}
By a similar computation as in the proof of \Cref{prop:general-interior-B},
\begin{align}
 W[\psi_+^0,\psi_-^{2\pi}](x_L)
 &=\frac{e^{2\pi iQ/h}}{\sin p(x_L)}
   \bigl(e^{ip(x_L)}-e^{-ip(x_L)}+O(h)\bigr)\\
 &=2ie^{2\pi iQ/h}(1+O(h)).
 \label{eq:exterior-dominant-Wronskian}
\end{align}
The relative error is uniform because the overlap stays away from
turning points.  
Because $W[\psi_+^0,\psi_-^0](x_L)=2i(1+O(h))$, we get that $ T_{22}=e^{2\pi iQ/h}(1+O(h))$.

\textbf{The off-diagonal entries.}
The off diagonal entries contain Wronskians of solutions with the same sign of their phase.  We then get
\begin{align}
 W[\psi_+^0,\psi_+^{2\pi}](x_R)
 &=O\!\left(he^{2\pi iQ/h}\right)\exp\left(\frac{2i}{h}
       \int_{2\pi}^{x_R}p(w)\dd w\right),
 \label{eq:exterior-plus-Wronskian}\\
 W[\psi_-^{2\pi},\psi_-^0](x_L)
 &=O\!\left(he^{2\pi iQ/h}\right)\exp\left(-\frac{2i}{h}
       \int_0^{x_L}p(w)\dd w\right).
 \label{eq:exterior-minus-Wronskian}
\end{align}
The sign conditions \eqref{eq:buffered-horizontal-signs}
and $p(2\pi-x)=p(x)$ imply
\begin{align}
 \Im\int_{2\pi}^{x_R}p(w)\dd w&>0,&
 \Im\int_0^{x_L}p(w)\dd w&<0.
\end{align}
Thus both exponential factors have modulus at most one.
Dividing by the local Wronskian gives
$T_{21},T_{12}=O\!\left(he^{2\pi iQ/h}\right)$.

\textbf{The remaining entry.}
Since $\det T=1$, we obtain
\begin{align}
 T_{11}
 =\frac{1+T_{12}T_{21}}{T_{22}}=e^{-2\pi iQ/h}(1+O(h))
   +O\!\left(h^2e^{2\pi iQ/h}\right).
\end{align}
Using that $\lvert e^{2\pi iQ/h}\rvert \ge e^{c/h}$, we get
$m_N = \tr T =e^{2\pi iQ/h}(1+O(h))$.
The fixed regular overlaps and strict phase inequalities make all
errors analytic in $E$ and locally uniform.
\end{proof}

\subsection{Identification of the actions}
\label{sec:classical-actions}

We now identify the contour integrals in these estimates with the
classical actions $P,Z,Q$ and substitute them into the trace identity
\eqref{eq:normalized-exact-reduction}.

The action $P$ is the normalized integral of the complex momentum over one period:
\begin{align}
 P(E)&\coloneq\frac{1}{2\pi}\int_0^{2\pi}p(z)\dd z,
 \qquad E\in\operatorname{int}\Sigma_\lambda.
 \label{eq:ordinary-action-from-Wronskian}
\end{align}
Here the contour is the real interval and the branch is the one fixed in
\Cref{sec:momentum-geometry}.  Under continuation, the contour and its
momentum branch are transported together.

The dual action $Z$ is the normalized integral of the position variable over one period:
\begin{align}
 Z(E)&\coloneq\frac{1}{2\pi}\int_0^{2\pi}z(p)\dd p.
 \label{eq:dual-action}
\end{align}
For $0\le p\le\pi$, the branch of $z(p)$ traces the real-momentum
arc from $\zeta_+$ to $\zeta_-$. Under continuation, the momentum contour and this
branch are transported together.

Observe that position--momentum duality simply interchanges the two closed
periods.
\begin{lemma}
\label{lem:duality-swaps-actions}
Let $\Im\lambda>0$ and $E\in\operatorname{int}\Sigma_\lambda$,
and use the dual variables in \eqref{eq:coefficient-dual-parameters}.
With corresponding contours and branches, the actions satisfy
\begin{align}
 \overline{P(E^\dagger,\lambda^\dagger)}
 &=Z(E,\lambda),&
 \overline{Z(E^\dagger,\lambda^\dagger)}
 &=P(E,\lambda).
 \label{eq:duality-swaps-actions}
\end{align}
\end{lemma}
\begin{proof}
Dividing the momentum equation by $\lambda$ and conjugating gives
\begin{align}
 \cos\overline z+\lambda^\dagger\cos\overline p
 &=E^\dagger.
\end{align}
Thus $(z,p)\mapsto(\overline p,\overline z)$ takes the original
equation to the dual equation. Choose the corresponding contours
and branches by this transformation, preserving their
parametrizations. The definitions of the actions then give
\begin{align}
 P(E^\dagger,\lambda^\dagger)
 &=\frac1{2\pi}\int_0^{2\pi}\overline{z(p)}\dd p
 =\overline{Z(E,\lambda)},\\
 Z(E^\dagger,\lambda^\dagger)
 &=\frac1{2\pi}\int_0^{2\pi}\overline{p(z)}\dd z
 =\overline{P(E,\lambda)}.
\end{align}
\end{proof}

\begin{corr}
\label{prop:general-interior-match}
For $\Im\lambda>0$, locally near every interior energy,
\begin{align}
 b&=\exp(-2\pi iP/h+h r_P),&
 a/c&=-\exp(2\pi iZ/h+h r_Z),
 \label{eq:target-relative-normalization}
\end{align}
where $r_P,r_Z$ are analytic and locally bounded in $E$.
\end{corr}
\begin{proof}
By \Cref{prop:general-interior-B}, we have
\begin{align}
 b = \mathcal{B}/\mathcal{A}=e^{-2\pi iP/h}(1+O(h)).
\end{align}
In particular, $b\ne0$ for sufficiently small $h$.
The identity $ac=-4b\lambda^N$ (\Cref{prop:exact-two-mode-reduction}) then gives $a,c\ne0$.
For the parameters in \eqref{eq:coefficient-dual-parameters},
\Cref{prop:coefficient-duality} then gives us
\begin{align}
 -\frac ac
 &=\overline{b(E^\dagger,\lambda^\dagger)}
 =\exp\left(\frac{2\pi i}{h}\overline{P(E^\dagger,\lambda^\dagger)}\right)
   (1+O(h))
 =e^{2\pi iZ/h}(1+O(h)).
\end{align}
The last equality follows from \Cref{lem:duality-swaps-actions}.
Taking logarithms of the relative errors near one proves the result,
locally uniformly in the interior.
\end{proof}

For the exterior, open-side, and corner calculations, it is helpful to define the following normalized contour integrals:
\begin{equation}\label{eq:action-definitions}
 \begin{array}{@{}c@{\qquad}c@{\qquad}c@{}}
  \chi_c\coloneq\frac1\pi\int_0^{\zeta_-}(\pi-p(z))\dd z,
  &\widetilde\chi_c\coloneq\frac1\pi\int_0^{\zeta_+}p(z)\dd z,
  &\widetilde\chi\coloneq\frac1\pi\int_{\zeta_+}^{\pi}p(z)\dd z.
 \end{array}
\end{equation}
Their contours (which avoid cuts of $p$) are shown in \Cref{fig:action-contours}.
In \Cref{prop:half-period-action-identity} we prove that $\widetilde\chi+\chi_c=-i\log\lambda$ (which we use below).

\actioncontoursfigure

\begin{lemma}
\label{lem:dual-action-contour}
For $\widetilde\chi_c$ and $\chi_c$ defined in \eqref{eq:action-definitions} and $Z$ defined in \eqref{eq:dual-action},
\begin{align}
 Z&=\chi_c+\widetilde\chi_c.
 \label{eq:second-action-from-Wronskians}
\end{align}
\end{lemma}

\begin{proof}
Integration by parts gives
\begin{align}
 \widetilde\chi_c
 &=\frac1\pi\bigl[zp(z)\bigr]_0^{\zeta_+}
   -\frac1\pi\int_{p(0)}^0 z(p)\dd p
 =\frac1\pi\int_0^{p(0)}z(p)\dd p,\\
 \chi_c
 &=\frac1\pi\bigl[z(\pi-p(z))\bigr]_0^{\zeta_-}
   +\frac1\pi\int_{p(0)}^\pi z(p)\dd p
 =\frac1\pi\int_{p(0)}^\pi z(p)\dd p.
\end{align}
The endpoint terms vanish because $p(\zeta_+)=0$ and
$p(\zeta_-)=\pi$.  The concatenated momentum contour deforms to
$[0,\pi]$ on the chosen inverse branch; see \Cref{fig:action-contours}.
Using \eqref{eq:dual-action} and $z(2\pi-p)=z(p)$,
we obtain
\begin{align}
 Z=\chi_c+\widetilde\chi_c
  =\frac{1}{\pi}\int_0^\pi z(p)\dd p
  =\frac{1}{2\pi}\int_0^{2\pi}z(p)\dd p.
\end{align}
\end{proof}

We next relate the position period to the half-period integrals
needed in the corner calculation.

\begin{lemma}
\label{lem:action-identity}
With $\arg\lambda\in(0,\pi)$,
\begin{align}
 P&=-\chi_c+\widetilde\chi_c-i\log\lambda.
 \label{eq:action-identity}
\end{align}
\end{lemma}
\begin{proof}
The reflection $z\mapsto2\pi-z$ preserves the chosen momentum branch
on the real axis.  Hence
\begin{align}
 P
 &=\frac{1}{\pi}\int_0^\pi p(z)\dd z
  =\frac{1}{\pi}\left(\int_0^{\zeta_+}p(z)\dd z
          +\int_{\zeta_+}^{\pi}p(z)\dd z\right)=\widetilde\chi_c+\widetilde\chi.
\end{align}
Here the concatenation of the contours from $0$ to $\zeta_+$
and from $\zeta_+$ to $\pi$ deforms to $[0,\pi]$ without crossing
the cuts.  By \Cref{prop:half-period-action-identity},
$\widetilde\chi+\chi_c=-i\log\lambda$, which gives
\eqref{eq:action-identity}.  The continued identity follows by
transporting all paths and branches together.
\end{proof}

For the open sides we also use the exact trace duality
\begin{align}
 m_N(E,\lambda)
 &=\lambda^N\overline{m_N(E^\dagger,\lambda^\dagger)},
 \label{eq:conjugated-dual-trace}
\end{align}
which follows from \eqref{eq:mN-duality} and \eqref{eq:mN-conjugation}.

On an open side, write $E=s+\lambda t$, where either
$s=\pm1$ and $-1<t<1$, or $t=\pm1$ and $-1<s<1$.
The transformation
$(E,\lambda)\mapsto(E^\dagger, \lambda ^\dagger)$
exchanges $s,t$, while $E\mapsto-E$ reverses both signs.
Together with \eqref{eq:conjugated-dual-trace} and
\eqref{eq:mN-sign-symmetries}, these operations reduce the proof on
every open side to $t=1$ (i.e., $E - \lambda \in (-1,1)$) with $\Im\lambda>0$.

Near $E_*=s+\lambda$, $-1<s<1$, the action $Q$ agrees with
the continued interior action $P$. Thus
\Cref{prop:exterior-Wronskian-estimates} gives
\begin{align}
 m_N(E,\lambda)&=e^{2\pi iP(E)/h}(1+O(h)).
\end{align}
To recover the two-action formula, we first show that
$\lambda^N e^{2\pi iZ/h}$ and $\lambda^N e^{-2\pi iZ/h}$
are exponentially smaller than this leading term.

\begin{lemma}
\label{lem:open-side-dominant-period}
At $E_*=s+\lambda$, $-1<s<1$, $\Im\lambda>0$, the continued
interior actions satisfy
\begin{align}
 -\Im P&>0,&
 -\Im P&>\log\lvert \lambda\rvert +\lvert \Im Z\rvert .
 \label{eq:open-side-dominance}
\end{align}
Both inequalities remain strict on a complex neighborhood of $E_*$.
\end{lemma}
\begin{proof}
For $0<x<2\pi$, we have
$\Im\cos p(x)=(\Im\lambda)(1-\cos x)>0$.
Since $0<\Re p(x)<\pi$, this implies $\Im p(x)<0$.
Integrating gives $\Im P<0$.

Here $p(0)=p_*=\arccos s$. The inverse branch is
\begin{align}
 z(p)&=\arccos\left(1+\frac{s-\cos p}{\lambda}\right),
 \qquad 0\leq\Re z(p)\leq\pi.
\end{align}
It has negative imaginary part for $0<p<p_*$ and positive imaginary
part for $p_*<p<\pi$.  Indeed, the imaginary part of its cosine has
the opposite sign to $s-\cos p$, and
$\Im\cos z=-\sin(\Re z)\sinh(\Im z)$.
Integration by parts along the two arms from zero gives
\begin{align}
 \widetilde\chi_c&=\frac1\pi\int_0^{p_*}z(p)\dd p,&
 \chi_c&=\frac1\pi\int_{p_*}^\pi z(p)\dd p.
\end{align}
Thus $\Im\widetilde\chi_c<0<\Im\chi_c$.  The half-period relation
\eqref{eq:action-identity} and the definition
\eqref{eq:second-action-from-Wronskians} of $Z$ give
\begin{align}
 -\Im P-\log\lvert \lambda\rvert -\Im Z
 &=-2\Im\widetilde\chi_c>0,\\
 -\Im P-\log\lvert \lambda\rvert +\Im Z
 &=2\Im\chi_c>0.
\end{align}
The actions are analytic at an open side: no turning points coalesce,
and their contours can be transported locally.  The strict inequalities
therefore persist on a smaller neighborhood.
\end{proof}

\subsection{Quantization in the three regular regimes}
\label{sec:regular-quantization}

To obtain exact cosine formulas from the WKB estimates, we use an
elementary observation: a relative $O(h)$ error in a dominant exponential
can be absorbed into an $O(h^2)$ correction of its phase.

\begin{lemma}[Correcting a dominant exponential]
\label{lem:dominant-exponential-correction}
Let $Q$ and $F_h$ be analytic on a fixed energy neighborhood.  Suppose
that $\Im Q\leq-c<0$ and
\begin{align}
 F_h(E)&=e^{2\pi iQ(E)/h}(1+O(h))
\end{align}
uniformly there.  For all sufficiently small $h$, there is an analytic
correction $Q_h$ such that, exactly,
\begin{align}
 F_h(E)&=2\cos\bigl(2\pi Q_h(E)/h\bigr),&
 Q_h&=Q+O(h^2).
 \label{eq:dominant-exponential-correction}
\end{align}
The correction estimate is uniform on the same neighborhood.
\end{lemma}
\begin{proof}
We aim to construct $Q_h$ satisfying
\begin{align}
 e^{2\pi iQ_h/h}+e^{-2\pi iQ_h/h}&=F_h.\label{eq:2445}
\end{align}
Setting $e^{2\pi iQ_h/h}=\Lambda_h$, \eqref{eq:2445} can be rewritten
\begin{align}
 \Lambda_h^2-F_h\Lambda_h+1&=0.
\end{align}
Since $F_h=e^{2\pi iQ/h}(1+O(h))$ is exponentially large, choose the root
close to $F_h$:
\begin{align}
 \Lambda_h
 &=\frac{F_h}{2}\left(1+\sqrt{1-\frac{4}{F_h^2}}\right)\\
 &=F_h\bigl(1+O(e^{-4\pi c/h})\bigr)
  =e^{2\pi iQ/h}(1+O(h)).
\end{align}
The square root is the analytic branch near one.

Consequently, $e^{-2\pi iQ/h}\Lambda_h$ is close to one and has an analytic
logarithm near zero.  Define
\begin{align}
 Q_h&\coloneq Q-\frac{ih}{2\pi}\log\left(e^{-2\pi iQ/h}\Lambda_h\right).
\end{align}
The logarithm is $O(h)$, so $Q_h-Q=O(h^2)$.  Moreover,
\begin{align}
 e^{2\pi iQ_h/h}&=\Lambda_h,&
 2\cos(2\pi Q_h/h)&=\Lambda_h+\Lambda_h^{-1}=F_h,
\end{align}
as required.
\end{proof}

\begin{theorem}[Quantization in the interior and across open sides]
\label{thm:interior-side-quantization}
Fix $\lambda\in\C\setminus\R$, and let $E_*$ lie in the interior or on
an open side of $\Sigma_\lambda$.  There is a complex neighborhood $U$
of $E_*$ on which the continued classical actions $P,Z$, defined in
\eqref{eq:ordinary-action-from-Wronskian} and
\eqref{eq:second-action-from-Wronskians}, have analytic corrections
\begin{align}
 P_h&=P+O(h^2),&Z_h&=Z+O(h^2)
 \label{eq:regular-side-corrected-actions}
\end{align}
such that, exactly for all sufficiently small $h$,
\begin{align}
 m_N(E,\lambda)
 &=2\cos(2\pi P_h(E)/h)+2\lambda^N\cos(2\pi Z_h(E)/h).
 \label{eq:regular-side-two-cosine}
\end{align}
The estimates are locally uniform for $E\in U$.
By \eqref{eq:finite-spectral-condition}, setting $m_N=2R_N$
gives the desired quantization condition \eqref{eq:intro-spectral-condition}, proving
\Cref{thm:main}\textup{(i)}.
\end{theorem}
\begin{proof}
\textbf{Interior energies.}
Suppose $\Im\lambda>0$. By \Cref{prop:general-interior-match},
choose the logarithms in \eqref{eq:exact-quantum-actions} so that
$P_h=P+O(h^2)$ and $Z_h=Z+O(h^2)$.
The exact factorization \eqref{eq:exact-action-factorization} gives
the result throughout the interior.

\textbf{The representative open side.}
Let $E_*=s+\lambda$, $-1<s<1$.  The horizontal Wronskian calculation
in \Cref{prop:exterior-Wronskian-estimates} gives
\begin{align}
 m_N&=e^{2\pi iP/h}(1+O(h))
\end{align}
on a full complex neighborhood of $E_*$.
The strict inequalities in \Cref{lem:open-side-dominant-period} show that
$\lambda^N e^{2\pi iZ/h}$ and $\lambda^N e^{-2\pi iZ/h}$ are exponentially
smaller than $e^{2\pi iP/h}$.  Consequently
\begin{align}
 T_h&\coloneq m_N-2\lambda^N\cos(2\pi Z/h)
       =e^{2\pi iP/h}(1+O(h)).
\end{align}
Apply \Cref{lem:dominant-exponential-correction} to $T_h$ with action $P$.
It gives an analytic $P_h$ satisfying
\begin{align}
 T_h&=2\cos(2\pi P_h/h),&P_h&=P+O(h^2).
\end{align}
Taking $Z_h=Z$ gives the claimed two-action formula.

\textbf{Symmetry reductions.}
Conjugated duality and the energy-sign symmetry give the other three
sides. For the energy-sign symmetry, the classical periods satisfy
\begin{align}
 P(E)&=\pi-P(-E),&
 Z(E)&=\pi-Z(-E).
\end{align}
Here the branches are continued along corresponding paths $E(t)$
and $-E(t)$. Define the corrected periods near the reflected energy
by the same formulas. The action shift
$\pi$ changes the cosine arguments by $2\pi^2/h=N\pi$,
multiplying both cosines by $(-1)^N$.
This is exactly the factor in the trace symmetry
$m_N(E,\lambda)=(-1)^Nm_N(-E,\lambda)$.

Complex conjugation gives the cases $\Im\lambda<0$.

The spectral statement follows from \Cref{lem:trace-polynomial}.
For the notation of \Cref{thm:main}, set $P_N\coloneq P_h$, $Z_N\coloneq Z_h$;
the action errors are then $O(N^{-2})$.
\end{proof}

The corrected actions are local and are not uniquely determined by the
trace.

\begin{theorem}[Exterior quantization]
\label{thm:exterior-quantization}
Fix $\lambda\in\C\setminus\R$, a simply connected exterior energy
domain $U\Subset\C\setminus\Sigma_\lambda$, and a compact set
$\mathcal K\Subset U$.
Let $Q$ be the exterior action in \eqref{eq:exterior-action}.
There is a locally analytic action $Q_h=Q+O(h^2)$, uniformly on
the compact set, for which
\begin{align}
 m_N(E,\lambda)
 &=2\cos(2\pi Q_h(E)/h)
 \label{eq:exterior-exact-cosine}
\end{align}
exactly for all sufficiently small $h$.  In particular,
$e^{-2\pi iQ/h}m_N=1+O(h)$.  The exterior spectral condition is obtained
by setting $m_N=2R_N$.
\end{theorem}
\begin{proof}
Since $\Im p\le-\eta<0$, one has $\Im Q\le-\eta<0$.
By \Cref{prop:exterior-Wronskian-estimates},
\eqref{eq:quarter-phase-trace}, and
\eqref{eq:exterior-Wronskian-entries}--
\eqref{eq:exterior-fourth-entry},
\begin{align}
 m_N(E,\lambda)&=e^{2\pi iQ/h}(1+O(h)).
 \label{eq:exterior-trace-from-Wronskians}
\end{align}

Apply \Cref{lem:dominant-exponential-correction} to $m_N$ with action $Q$.
It gives an analytic correction satisfying
\begin{align}
 m_N&=2\cos(2\pi Q_h/h),&Q_h&=Q+O(h^2).
 \label{eq:exterior-quantum-action-Wronskian}
\end{align}
The growing multiplier is therefore
\begin{align}
 \Lambda_h&\coloneq e^{2\pi iQ_h/h},&\Lambda_h+\Lambda_h^{-1}&=m_N.
 \label{eq:exterior-exact-multiplier}
\end{align}
It is the growing eigenvalue of $T$, since
$\det T=1$.  The normalized asymptotic is
\eqref{eq:exterior-trace-from-Wronskians}.
\end{proof}

For the notation of \Cref{thm:main}, set $Q_N\coloneq Q_h$.  Then
$Q_N-Q=O(N^{-2})$ and $2\pi Q_h/h=NQ_N$.  Together with
\Cref{lem:trace-polynomial}, this proves part \textup{(ii)} of
\Cref{thm:main}.

\section{Corner quantization}
\label{sec:corner-quantization}

We now estimate the trace on a fixed neighborhood of a corner.  A
canonical domain avoiding the coalescing turning points determines $c$
uniformly, but gives only an upper bound for $\mathcal{B}$.  We
recover $\mathcal{B}$ by locating its zeros and determining the
remaining nonzero factor.  Ordinary WKB expansion on a fixed energy
contour and the maximum principle replace the local connection
calculation.

\subsection{A uniform WKB comparison near the corner}
\label{sec:unif-wkb-comp}

First assume $\Im\lambda>0$, and put
\begin{align}
  E_+&\coloneq1+\lambda,
  &\alpha^2&=\lambda,
 &\arg\alpha&=\tfrac12\arg\lambda\in(0,\pi/2).
\end{align}
We remove this restriction at the end of the section.  Throughout,
$h=2\pi/N$ and the solutions retain their global strip normalization.
In particular $\mathcal{A}=-2i$ and $b=\mathcal{B}/\mathcal{A}$.

The coalescing turning points are $-\zeta_+(E)$ and
$\zeta_+(E)$, where $p=0$ and
\begin{align}
  \cos\zeta_+(E)&=\frac{E-1}{\lambda}, &\text{ one has }
  \zeta_+(E)&\longrightarrow0\text{ as }E\longrightarrow E_+.
\end{align}
The turning points where $p=\pi$ stay separated.  Let $\gamma_+(E)$ go from
$-\zeta_+(E)$ to $\zeta_+(E)$ on one momentum sheet and return on the
other.  Orient it so that the vanishing action
\begin{align}
 V_+(E)&\coloneq\frac1{2\pi}\int_{\gamma_+(E)}p\dd z,
 &V_+(E_+)&=0,
 &V_+'(E_+)&=-\frac1\alpha
 \label{eq:upper-right-vanishing-action-normalization}
\end{align}
is analytic at $E_+$.  In particular,
\begin{align}
 V_+(E)&=-\frac1\alpha(E-E_+)+O((E-E_+)^2).
 \label{eq:corner-vanishing-action-linearization}
\end{align}
The contraction and orientation are shown in
\Cref{fig:vanishing-cycle}.  An explicit formula for its derivative
appears in \Cref{lem:vanishing-period-derivative}.

\vanishingcyclefigure

The vanishing action $V_+$ can be expressed in terms of the $P$ and $Z$ action integrals.

\begin{lemma}
\label{lem:vanishing-action-two-periods}
For $\Im\lambda>0$ and $E$ in the interior sufficiently near $1+\lambda$,
with the branches fixed above,
\begin{align}
 V_+(E)&=P(E)+Z(E)+i\log\lambda.
 \label{eq:upper-right-action-difference}
\end{align}
The identity persists under simultaneous analytic continuation of the
actions and their defining contours.
\end{lemma}
\begin{proof}
We compute
\begin{align}
 V_+(E)=\frac1{2\pi}\left(\int_{-\zeta_+}^{\zeta_+}p(z)\dd z
 +\int_{\zeta_+}^{-\zeta_+}(-p(z))\dd z\right) =\frac2\pi\int_0^{\zeta_+}p(z)\dd z
 =2\widetilde\chi_c,
\end{align}
recalling the definition of $\widetilde\chi_c$ in \eqref{eq:action-definitions}.
Here $p(-z)=p(z)$ on the outgoing sheet, while the return sheet
has momentum $-p$. The orientation agrees with
\eqref{eq:upper-right-vanishing-action-normalization}, since
$2\widetilde\chi_c=-(E-E_+)/\alpha+O((E-E_+)^2)$.
Next, \eqref{eq:action-identity} and
\eqref{eq:second-action-from-Wronskians} give
\begin{align}
 P+Z+i\log\lambda
 =-\chi_c+\widetilde\chi_c-i\log\lambda
   +\chi_c+\widetilde\chi_c+i\log\lambda=2\widetilde\chi_c=V_+.
\end{align}
Continuing all contours and branches together preserves the identity.
\end{proof}

We construct a canonical domain for $\phi_+^0$ around a line that runs
close to the upper part of the real-momentum arc, passes to the right
of its endpoints, and extends to both imaginary infinities.  A local
deformation across a regular part of the arc gives an upper overlap
with the domain for $\phi_+^{2\pi}$.  On this overlap, the momentum
branches are opposite, so the leading terms in their Wronskian do not
cancel.  \Cref{fig:stable-corner-domain} shows the domain and the
positive real point used to bound $\mathcal{B}$.
\stablecornerdomainfigure

\begin{prop}[Uniform comparison at a corner]
\label{prop:stable-corner-comparison}
There are a fixed complex neighborhood $U$ of $E_+$ and
$\gamma>2\pi\log\max(1,\lvert \lambda\rvert )$
such that, uniformly on $U$,
\begin{align}
 c(E,h)&=2e^{-i\pi V_h(E)/h},
 &V_h(E)&=V_+(E)+O(h^2),
 &\lvert b(E,h)\rvert &\le Ce^{-\gamma/h}.
 \label{eq:corner-corrected-action}
\end{align}
The corrected action $E\mapsto V_h(E)$ is analytic
.
\end{prop}

Appendix~\ref{app:corner-stable-domain} constructs the domain and proves
these estimates.  The logarithm of the relative error in $c$ defines
$V_h$; the bound for $b$ uses WKB at a fixed positive real point and
the exact recurrence back to zero.  Neither estimate uses WKB at the
double turning point.

Since $c\ne0$, \eqref{eq:division-free-reduction} gives
$a=-4b\lambda^N/c$ and hence
\begin{align}
 b\,m_N&=1+b^2+\frac{a^2+c^2}{4}=1+\frac{c^2}{4}
       +\left(1+\frac{4\lambda^{2N}}{c^2}\right)b^2.
 \label{eq:corner-division-free-trace}
\end{align}
We do not divide by $b$: its zeros are the next ingredient of the
proof.

\subsection{Determining the corner factor}

Put
\begin{align}
 w_h&\coloneq\frac{V_h}h.
\end{align}
Thus $hw_h=V_++O(h^2)$. All neighborhoods below are fixed independently
of $h$ and can be made smaller a finite number of times.

\begin{lemma}
\label{lem:corner-B-zeros}
On a sufficiently small neighborhood of $E_+$, the zeros of $b$ are
simple and are exactly the points
\begin{align}
 w_h(E,h)&=-n-\tfrac12,\qquad n=0,1,2,\ldots,
 \label{eq:corner-B-zero-levels}
\end{align}
that belong to the neighborhood.
\end{lemma}

The proof in Appendix~\ref{app:corner-B-zeros} has two ingredients.  A
zero of $b$ forces $\phi_+^0(-z)=\pm\phi_+^0(z)$.  Its decay on the
imaginary axis then places the energy in the outward cone
\begin{align}
 E_++(0,\infty)+\lambda(0,\infty);
\end{align}
see \eqref{eq:corner-outward-zero-identity}.  The identity
\eqref{eq:corner-division-free-trace} restricts these zeros to
$\{E;\; w_h\in\Z+\tfrac12\}$. Rouch{\'e}'s theorem applied to
\eqref{eq:corner-reference-scalar} gives one simple zero of
$b(m_N-2R_N)$ near each energy in $\{E:\; w_h\in\Z+\tfrac12\}$.  The
spectrum of $H_{0,0,\lambda}^N$ lies in $\Sigma_\lambda$ and consists
of the zeros of $m_N-2R_N$. So every zero of $b(m_N-2R_N)$ outside
$\Sigma_\lambda$ is a zero of $b$.  The
identity~\eqref{eq:corner-division-free-trace} then places it at the
half-integer level itself.

Consequently, the function
\begin{align}
  E\mapsto F_h(E)&\coloneq b(E,h)\Gamma\bigl(\tfrac12+w_h(E,h)\bigr)
 \label{eq:corner-Gamma-factorization}
\end{align}
is analytic and nowhere zero: the poles of the Gamma factor cancel
exactly the simple zeros of $b$.  To normalize this factor, we first
remove the logarithmic singularity of $P$.

Choose $\log V_+$ real on the inward curve $V_+(E)>0$ and set
\begin{align}
 J(E)&\coloneq-2\pi iP(E)+V_+(E)\bigl(\log V_+(E)-1\bigr).
 \label{eq:corner-regularized-action}
\end{align}

\begin{lemma}
\label{lem:corner-regularized-action}
The function $J$ extends analytically through $E_+$.  In the outward
cone $E_++(0,\infty)+\lambda(0,\infty)$, let $Q$ be the real-period
action with $\Im p<0$.  Then
\begin{align}
 J&=-2\pi iQ+V_+\bigl(\log(-V_+)-1-i\pi\bigr),
 \label{eq:corner-exterior-regularization}
\end{align}
where $\log(-V_+)$ is real on the outward curve $V_+(E)<0$.
\end{lemma}

The contour and local calculation are given in
Appendix~\ref{app:corner-regularized-action}.  The change of $-2\pi iP$
around the corner cancels that of $V_+\log V_+$, leaving a removable
singularity.
When the coupling is displayed, we write $Q_+(E;\lambda)=Q(E)$.

\begin{prop}[Uniform corner trace]
\label{prop:uniform-corner-trace}
On a sufficiently small fixed neighborhood of $E_+$, say,
$\mathcal{V}_+$, there exists an analytic function
$E\in\mathcal{V}_+\mapsto J_h(E)=J(E)+O(h^2)$ with uniform error such
that, for all $N$ sufficiently large and $h=2\pi/N$, one has
\begin{align}
  m_N(E,\lambda)
  &=\frac{\sqrt{2\pi}\,h^{w_h} e^{-J_h/h-\pi iw_h}}
    {\Gamma(\tfrac12-w_h)} +
    \frac{\sqrt{2\pi}\,(1+\lambda^{2N}e^{2\pi iw_h})
    h^{-w_h}e^{J_h/h}}
    {\Gamma(\tfrac12+w_h)}
 \label{eq:corner-two-Gamma-trace}
\end{align}
with $w_h=V_h/h$ from \eqref{eq:corner-corrected-action} 
.  Both terms in the right-hand side
of~\eqref{eq:corner-two-Gamma-trace} are analytic through their
reciprocal-Gamma zeros.  On $\mathcal{V}_+$, the second term is
bounded by $\eta^{-1}\max(1,\lvert \lambda\rvert )^N e^{-\eta/h}$ for
some $\eta>0$.
\end{prop}

\begin{proof}
  The detailed WKB and Stirling calculation is given in
  \Cref{lem:corner-boundary-normalization} in
  Appendix~\ref{app:corner-boundary-normalization}.  It gives, on a
  fixed small contour around $E_+$,
  \begin{align}
    \frac{F_h(E)}{\sqrt{2\pi}\,e^{J(E)/h}h^{-w_h(E,h)}}
    &=1+O(h).
      \label{eq:corner-boundary-ratio}
  \end{align}

  \cornernormalizationfigure

  Here the interior (i.e., inside $\Sigma_\lambda$) estimate for $b$
  and exterior (i.e., outside $\Sigma_\lambda$) estimate for $m_N$ give
  the same normalization, by Gamma reflection and
  \eqref{eq:corner-exterior-regularization}.  Open-side estimates
  cover the transitions. \Cref{fig:corner-boundary-normalization}
  shows the contour in the $E$-plane. The map $E\mapsto V_+(E)$ is
  biholomorphic from a neighborhood of $E_+$ onto a neighborhood of $0$.

  We now complete the proof of
  Proposition~\ref{prop:uniform-corner-trace} by the maximum principle
  and the exact Wronskian identity. The quotient in
  \eqref{eq:corner-boundary-ratio} is analytic inside the contour.
  The maximum principle therefore gives the same estimate throughout
  its interior.  Its analytic logarithm tending to 0 defines
  $J_h=J+O(h^2)$ so that one has
  \begin{align}
    F_h&=\sqrt{2\pi}\,e^{J_h/h}h^{-w_h}.
  \end{align}
  Finally, the Gamma reflection formula and
  \eqref{eq:corner-division-free-trace} give
  \begin{align}
    F_h\bigl(m_N-(1+\lambda^{2N}e^{2\pi iw_h})b\bigr)
    &=\frac{2\pi e^{-\pi iw_h}}{\Gamma(\tfrac12-w_h)}.
  \end{align}
  Substitute the expression for $F_h$ and $b=F_h/\Gamma(\tfrac12+w_h)$
  to obtain \eqref{eq:corner-two-Gamma-trace}.  Its second term is
  exactly $(1+\lambda^{2N}e^{2\pi iw_h})b$.  The strict inequality for
  $\gamma$ in \Cref{prop:stable-corner-comparison}, and possibly shrinking
  the neighborhood of $E_+$, give the stated bound.
\end{proof}

The corrected $w_h$ in $h^{-w_h}$ avoids an $O(h\lvert \log h\rvert )$ relative error.
No small reciprocal-Gamma factor has been divided out, so the formula
remains informative at its zeros and at the corner.

\subsection{Quantization and the other corners}
\label{sec:quant-other-corn}

For the symmetry reductions, write the leading term using
$\nu_{+,N}=-w_h$:
\begin{align}
 \nu_{+,N}(E)&\coloneq-\frac{V_h(E)}h,
 \\
 K_{+,N}(E)
 &\coloneq\sqrt{\frac\pi2}\,h^{-\nu_{+,N}(E)}
       e^{-J_h(E)/h+i\pi\nu_{+,N}(E)}.
 \label{eq:corner-leading-prefactor}
\end{align}
The prefactor $K_{+,N}$ never vanishes.  { Define the remainder in
  the half-trace formula by
\begin{align}
 r_{+,N}(E)
 &\coloneq\frac12\bigl(1+\lambda^{2N}e^{2\pi iw_h(E,h)}\bigr)b(E,h).
 \label{eq:corner-half-trace-remainder}
\end{align}
}
In this notation,
\begin{align}
 \frac12m_N(E,\lambda)
 &=\frac{K_{+,N}(E)}{\Gamma(\nu_{+,N}(E)+1/2)}
     +r_{+,N}(E),
 &\lvert r_{+,N}\rvert &\le C\max(1,\lvert \lambda\rvert )^N e^{-\eta/h}.
 \label{eq:corner-leading-trace}
\end{align}

To obtain the form in \Cref{thm:main-corner}, put
$w\coloneq V_+/h$.  By
\eqref{eq:upper-right-action-difference} and
\eqref{eq:corner-regularized-action},
\begin{align}
 \frac\pi h\log\lambda+\frac{i\pi(P-Z)}h+w(1-\log w)
 &=-\frac Jh+w(\log h-i\pi).
\end{align}
Consequently
\begin{align}
 a_h(E)&\coloneq\exp\left(
 -\frac{J_h-J}{h}+(w_h-w)(\log h-i\pi)\right)
 =1+O(h\lvert \log h\rvert ).
 \label{eq:corner-intro-correction}
\end{align}
With $h=2\pi/N$, set $a_N\coloneq a_h$ and $w_N\coloneq w_h$.
This factor is analytic on the corner neighborhood, since
$J,J_h,w,w_h$ are analytic there.  Substitution in
\eqref{eq:corner-leading-trace} gives the left-hand side of
\eqref{eq:intro-corner-spectral-condition} as an exact expression for
$m_N/2$. The same calculation shows that its combined classical
prefactor is $\sqrt{\pi/2}\exp(-J/h+w(\log h-i\pi))$, analytic and
nonzero at the corner.


We next remove the half-plane restriction on $\lambda$.  Complex
conjugation gives
\begin{align}
 m_N(\overline E,\overline\lambda)
 &=\overline{m_N(E,\lambda)},
 \\
 \overline{H^N_{p_0,x_0,\lambda}}
 &=H^N_{-\overline{p_0},\overline{x_0},\overline\lambda}.
 \label{eq:corner-mN-conjugation}
\end{align}
For $\Im\lambda<0$, define $V_+,V_h,K_{+,N}$ and $r_{+,N}$
by conjugating their values at $(\overline E,\overline\lambda)$.
We conjugate the prefactor itself, not just $J_h$; this also conjugates
the factor $e^{i\pi\nu_{+,N}}$ in
\eqref{eq:corner-leading-prefactor}.

The upper-right formula now holds for every nonreal coupling, so the
sign symmetries may be applied even when they change its half-plane.
For $\epsilon_p,\epsilon_z\in\{1,-1\}$, put
\begin{align}
 E_c&\coloneq\epsilon_p+\lambda\epsilon_z,
 \label{eq:corner-energy-definition}\\
 \mu_c&\coloneq\epsilon_p\epsilon_z\lambda,
 &\widetilde E&\coloneq\epsilon_pE,
 \\
 \alpha_c^2&=\mu_c,
 &\arg\alpha_c&=\tfrac12\arg\mu_c\in(-\pi/2,\pi/2),\\
 V_c(E;\lambda)&\coloneq V_+(\widetilde E;\mu_c),
 &V_{c,h}(E;\lambda)&\coloneq V_h(\widetilde E;\mu_c).
 \label{eq:corner-action-definition}
\end{align}
When an exterior action is needed, retain the convention
$Q_c(E;\lambda)\coloneq Q_+(\widetilde E;\mu_c)+\frac\pi2(1-\epsilon_p)$
modulo $2\pi$; for the lower half-plane,
$Q_+(E;\lambda)=-\overline{Q_+(\overline E;\overline\lambda)}$.
These choices keep $\Im p<0$ on the exterior real period.

\begin{corr}
\label{cor:corner-symmetry-reduction}
For every nonreal $\lambda$ and each corner $E_c$,
\begin{align}
 V_c(E_c;\lambda)&=0,
 &\partial_EV_c(E_c;\lambda)&=-\frac{\epsilon_p}{\alpha_c}.
 \label{eq:corner-vanishing-action-normalization}
\end{align}
On a fixed neighborhood of $E_c$, define
\begin{align}
 \nu_{c,N}&\coloneq-\frac{V_{c,h}}h,
 &K_{c,N}(E;\lambda)
 &\coloneq\epsilon_p^NK_{+,N}(\widetilde E;\mu_c),\\
 r_{c,N}(E;\lambda)
 &\coloneq\epsilon_p^Nr_{+,N}(\widetilde E;\mu_c).
\end{align}
Then $K_{c,N}$ is nonzero and
$\nu_{c,N}=-V_c/h+O(h)$ uniformly on that neighborhood, and
\begin{align}
 \frac12m_N(E,\lambda)
 &=\frac{K_{c,N}(E;\lambda)}
         {\Gamma(\nu_{c,N}(E)+1/2)}+r_{c,N}(E;\lambda),
   \label{eq:all-corner-WKB-matching}\\
  \lvert r_{c,N}\rvert &\le \eta^{-1}\max(1,\lvert \lambda\rvert )^N e^{-\eta N}.
\end{align}
Setting this expression equal to $R_N$ extends \Cref{thm:main-corner}
to every corner and every nonreal coupling, for arbitrary complex phases.
\end{corr}
\begin{proof}
Apply $m_N(\widetilde E,\mu_c)=\epsilon_p^Nm_N(E,\lambda)$ to the
upper-right formula.  The chain rule gives the normalization of $V_c$;
conjugation and sign changes preserve the remainder bound, with a
change of $\eta$ when $h$ is replaced by $N$.
\end{proof}

For real phases the same argument gives a familiar consequence.  Near
$1+\lambda$, the inward eigenvalues satisfy
\begin{align}
 V_h(E_n)&=h(n+\tfrac12)+O(he^{-\eta/h}),
 \label{eq:corner-uniform-Bohr-Sommerfeld}
\end{align}
uniformly for the positive action levels in a fixed smaller interval.
For each fixed $n$, this reduces to
$E_n=1+\lambda-\alpha h(n+\tfrac12)+O(h^2)$.
The action statement, unlike this linearization, also covers
$n$ of order $h^{-1}$ within the fixed neighborhood.  The zero count
is part of Appendix~\ref{app:corner-B-zeros}.

\begin{rem}
\label{sec:corner-obstruction}
\label{prop:corner-obstruction}
In the upper-half-plane calculation, at $E=E_+$ write
\begin{align}
 S&\coloneq\frac1{2\pi}\int_0^{2\pi}
       \operatorname{arcosh}(1+\lambda(1-\cos x))\dd x
\end{align}
with positive real part of the inverse hyperbolic cosine.  Then
$J(E_+)=-2\pi S$ and $w_h(E_+)=O(h)$.  Hence
\begin{align}
 e^{-2\pi S/h}m_N(E_+,\lambda)&\longrightarrow\sqrt2.
 \label{eq:corner-trace-limit}
\end{align}
The continued classical two-cosine expression instead has normalized
limit $2$, rather than $\sqrt2$.  Corrections of size $O(h^2)$ to the classical actions
cannot change this leading coefficient.  The reciprocal-Gamma
factor is therefore essential.
\end{rem}

\section{Global continuation of the two-action formula}
\label{sec:global-two-action-continuation}

In this section we use $N=2\pi/h$ and write $P_N,Z_N,Q_N$ for the
corrected actions previously denoted by $P_h,Z_h,Q_h$.
So far, we have obtained the two-action formula
\begin{align}
 m_N&=2\cos(NP_N)+2\lambda^N\cos(NZ_N)
 \label{eq:2981}
\end{align}
in the interior and across the open sides
(\Cref{thm:interior-side-quantization}), and the one-action formula
\begin{align}
 m_N&=2\cos(NQ_N)
\end{align}
in the exterior (\Cref{thm:exterior-quantization}).
The corner analysis in \Cref{sec:corner-quantization} supplies a uniform formula near the remaining
energies. We now ask where the exterior formula can also be expressed
using the continued interior actions $P,Z$.

For any $\lambda \in \C \setminus \R$, define
\begin{align}
 \Gamma_\lambda
 &\coloneq
 \set{E=1+\lambda+t+\lambda s :
       t,s\geq0,\ \Im V_+(E;\lambda)=0},
 \label{eq:outward-equal-growth-curve}
\end{align}
where $V_+$, defined in
\eqref{eq:upper-right-vanishing-action-normalization},
is analytically continued from $1+\lambda$ into the outward cone
$\R_{\geq1}+\lambda\R_{\geq1}$.
We show in Appendix~\ref{app:exterior-cuts} that $\Gamma_\lambda$ is a curve originating from $1+\lambda$ and going off to infinity.
Set
\begin{align}
 \mathcal W_\lambda
 &\coloneq\Gamma_\lambda\cup(-\Gamma_\lambda)
       \cup\Gamma_{-\lambda}\cup(-\Gamma_{-\lambda})
 \label{eq:global-cut-system}
\end{align}
as the union of the four curves emanating from the four corners of $\Sigma_\lambda$.

The aim of this section is to prove the following theorem.

\begin{theorem}[The global two-action formula]
\label{cor:global-two-action}
Fix $\lambda\in\C\setminus\R$.
The actions $P,Z$ extend analytically from
$\operatorname{int}\Sigma_\lambda$ to
$\C\setminus\mathcal W_\lambda$. Every $E_*\in\C\setminus\mathcal W_\lambda$ has a neighborhood
$U\subset\C\setminus\mathcal W_\lambda$ on which, for all sufficiently
large $N$, there are analytic functions $P_N,Z_N$ satisfying,
uniformly for $E\in U$,
\begin{align}
 P_N-P&=O(N^{-2}),&Z_N-Z&=O(N^{-2}),\\
 m_N(E,\lambda)
 &=2\cos\left(NP_N(E)\right)
   +2\lambda^N\cos\left(NZ_N(E)\right).
\end{align}
\end{theorem}

Here we stress that the corrected actions and their error bounds are local.

\subsection{Recovering two cosines from one dominant exponential}

Recall $c_\lambda=i\log\lambda$ from \eqref{eq:exterior-log-period},
with the logarithm used in \eqref{eq:action-identity}.
Since $\lambda^N=e^{-iNc_\lambda}$, the classical expression
associated with \eqref{eq:2981} is the sum of exponentials with the four actions
\begin{align}
 P,\qquad -P,\qquad Z-c_\lambda,\qquad -Z-c_\lambda.
 \label{eq:four-classical-exterior-actions}
\end{align}
Note that the term with the strictly smallest imaginary part will dominate the others. If this term is $Q$ modulo
$2\pi\Z$, the following proposition gives the two-action formula in
the exterior of $\Sigma_\lambda$.

\begin{prop}
\label{cor:exterior-two-action}
On an open set $U\subset\C\setminus\Sigma_\lambda$, fix continued actions $P,Z$ from $\operatorname{int}\Sigma _\lambda$
and the exterior action $Q$.  Suppose one action in
\eqref{eq:four-classical-exterior-actions} equals $Q$ modulo
$2\pi\Z$, and the other three have strictly larger imaginary part.
Then there are analytic corrections $P_N-P,Z_N-Z=O(N^{-2})$ such that
\begin{align}
 m_N&=2\cos(NP_N)+2\lambda^N\cos(NZ_N)
\end{align}
exactly.  The estimates are uniform on compact subsets where the three
strict inequalities hold.
\end{prop}
\begin{proof}
First suppose the selected action is $P$, so that
$e^{iNQ}=e^{iNP}$. Then, by
\Cref{thm:exterior-quantization},
\begin{align}
 m_N-2\lambda^N\cos(NZ)&=e^{iNP}(1+O(N^{-1})).
\end{align}
By \Cref{lem:dominant-exponential-correction}, the right-hand side of this equality can be written as
$2\cos(NP_N)$ for an analytic $P_N=P+O(N^{-2})$.
Taking $Z_N=Z$ proves the claim.

If the selected action is $-P$, replace $P$ by $-P$ in this calculation;
the cosine is unchanged.  If it is $\sigma Z-c_\lambda$, with
$\sigma\in\{1,-1\}$, keep $P_N=P$ and instead use
\begin{align}
 \frac{m_N-2\cos(NP)}{\lambda^N}
 &=e^{iN\sigma Z}(1+O(N^{-1})).
\end{align}
Strict dominance over the term with action $-\sigma Z-c_\lambda$
gives $\Im(\sigma Z)<0$.  The same lemma corrects $\sigma Z$ by
$O(N^{-2})$; multiplying this corrected action by $\sigma$ defines $Z_N$.
Since cosine is even, this gives the required formula.
\end{proof}

\subsection{The four cuts and the proof of the theorem}
\label{sec:exterior-walls}

To see why these cuts matter, first consider the corner $1+\lambda$ (the others follow by symmetry).
For $\Im\lambda>0$, \Cref{lem:vanishing-action-two-periods} gives $V_+ = P + Z + c_\lambda$.
Thus $V_+$ is the difference between the actions $P$ and $-Z - c_\lambda$.
The corresponding exponentials have equal magnitude precisely when $\Im V_+=0$, which defines $\Gamma_\lambda$ within the outward cone.
Complex conjugation gives the case $\Im\lambda<0$.

To prove \Cref{cor:global-two-action}, we use the following proposition, which is proved in Appendix~\ref{app:exterior-cuts}.

\begin{prop}[Exterior cuts and the dominant action]
\label{prop:general-exterior-cuts}
For every nonreal $\lambda$, the four curves in
\eqref{eq:global-cut-system} are disjoint and extend from their corner
endpoints to infinity. On their complement, $P,Z$ can be continued
together from the interior.
At every exterior energy outside $\mathcal W_\lambda$, exactly one
action in \eqref{eq:four-classical-exterior-actions} agrees with $Q$
modulo $2\pi\Z$; the other three have strictly larger imaginary part.
\end{prop}

The dotted curves in \Cref{fig:global-two-action-regions} show
$\mathcal W_\lambda$, and the labels identify the dominant action in
each exterior region. The two-action formula extends into the exterior
away from these cuts. The one-action formula remains valid on the cuts
as well.
Appendix~\ref{app:exterior-cuts} constructs the outward curve from
the signs of the derivatives of $\Im V_+$ and proves the ordering of all four terms.
It also gives the action-difference identities at the other corners.

\globaltwoactionfigure

At fixed $\lambda$, the cuts have the central symmetry $E\mapsto-E$;
the reduction $\lambda\mapsto-\lambda$ is a change of coupling, not
in general a rotation in the same energy plane.
For $\lambda=i\rho$, $\rho>0$, conjugation gives reflections in both
coordinate axes. In particular,
\begin{align}
 \Gamma_i&=\{t(1+i):t\geq1\},
\end{align}
so the four cuts for the Scottish flag are the outward diagonal rays.
The more general imaginary-coupling curves are described in
\Cref{cor:imaginary-exterior-walls}.

\begin{proof}[Proof of \Cref{cor:global-two-action}]
Use \Cref{thm:interior-side-quantization} in the interior and on the
open sides. In the exterior, \Cref{prop:general-exterior-cuts} supplies
the dominant action, and \Cref{cor:exterior-two-action} gives the
two-action formula.
\end{proof}

The exclusion of the cuts is essential in general:
\Cref{prop:exterior-obstruction} shows that, for $\lambda=i$, the
continued two-action expression has an extra leading term on a cut
that $O(N^{-2})$ corrections to $P,Z$ cannot remove.
On the exterior cuts away from the corner, the one-action formula of
\Cref{thm:exterior-quantization} remains valid. At their corner
endpoints, \Cref{cor:corner-symmetry-reduction} supplies the
reciprocal-Gamma formula on a fixed neighborhood.

\appendix

\section{WKB estimates and continuation}
\label{app:wkb-tools}

We establish uniform WKB estimates, normalize the entire solutions,
and construct the comparison regions used in the proof.
We write the finite difference equation as
\begin{align}
 \psi(z+h)+\psi(z-h)+v(z)\psi(z)&=0,
 \label{eq:main}\\
 2\cos p(z)+v(z)&=0.
 \label{eq:momentum}
\end{align}
Our potential is $v(z)=2\lambda\cos z-2E$.
A regular domain avoids $v=\pm2$.
The notation $K_0\Subset K_1$ means that the closure of $K_0$ is
compact and lies inside $K_1$; this leaves a positive distance from
the boundary for the estimates.

\subsection{Uniform estimates on bounded domains}
\label{sec:uniform-wkb-bounded}
We use the following version of the WKB theorem on nested domains.
The parameter is the energy $E$ in our application.
\begin{theorem}[Uniform WKB estimates]
  \label{thm:uniform-wkb}
  Let $\mathcal Q\subset\C^d$ be compact, let $D\subset\C$ be a domain,
  and let $v_\alpha(z)$ be analytic in $(\alpha,z)$ near
  $\mathcal Q\times D$.  Suppose that
  $K_0\Subset K_1\Subset D$ are bounded simply connected domains,
  $K_1$ is horizontally connected, and the following data are fixed:
  \begin{enumerate}
  \item a branch $p_\alpha$ of
    $2\cos p_\alpha+v_\alpha=0$, analytic in $(\alpha,z)$ near
    $\mathcal Q\times\overline K_1$, together with one branch of
    $\sqrt{\sin p_\alpha}$, and a constant $c_0>0$ such that
    $\abs{\sin p_\alpha}\ge c_0$ there;
  \item a compact parameter set $\mathcal S$ and a $C^1$ map
    $\Gamma\colon\mathcal S\times[0,1]\to\overline K_1$, with uniformly
    bounded $C^1$ norm, such that the vertical curves
    $\gamma_s=\Gamma(s,\cdot)$ have the same two boundary endpoints
    and sweep out $K_1$, together with a number $\eta>0$ such that, on every
    $\gamma_s$ and for every $\alpha\in\mathcal Q$,
    \begin{align}
       \Im\bigl(p_\alpha(\gamma_s)\dot\gamma_s\bigr)&\ge\eta,
       &\Im\bigl((p_\alpha(\gamma_s)-\pi)\dot\gamma_s\bigr)&\le-\eta.
       \label{eq:uniform-strict-canonical-margin}
    \end{align}
  \end{enumerate}
  Fix $z_0\in K_0$.  Then there are an open neighborhood
  $\mathcal U\supset\mathcal Q$, numbers $h_0,C>0$, and two exact solutions
  $\psi_\pm(\alpha,h;z)$ of \eqref{eq:main} in $K_1$, analytic on
  $\mathcal U\times K_1$,
  such that, for $0<h<h_0$ and $z\in K_0$,
  \begin{align}
    \label{eq:asymp:wkb}
    \psi_\pm(\alpha,h;z)
      =\frac1{\sqrt{\sin p_\alpha(z)}}
       \exp\left(\pm\frac{i}{h}\int_{z_0}^z p_\alpha(\zeta)\dd\zeta
                 +h r_\pm(\alpha,h;z)\right).
  \end{align}
  Here $\abs{r_\pm}\le C$ uniformly on $K_0$.  Moreover, on $K_0$,
  \begin{align}
    W[\psi_+,\psi_-](z)=2i\bigl(1+O(h)\bigr),
    \label{eq:uniform-wkb-wronskian}
  \end{align}
  and hence this Wronskian does not vanish when $h$ is small.
\end{theorem}
\begin{proof}
\textbf{Fixed-parameter theorem and uniformity.}
For one fixed value of $\alpha$, this is the scalar specialization
\cite[Theorem~2.2 and (2.14)]{FedotovShchetka2022Berry} of the bounded-domain
WKB theorem.  We explain why its estimates are uniform here.

The cited construction diagonalizes the leading first-order system and
obtains the correction by a Neumann series on the canonical curves; see
\cite[Sections~4.1--4.2]{FedotovShchetka2022Berry} and
\cite[Sections~4--6]{MR3501077}.  The constants used there
come from Cauchy bounds for $v_\alpha,p_\alpha$, the lower bound on
$\lvert \sin p_\alpha\rvert $, and the canonical inequalities.  These are uniform
on $\mathcal Q\times\overline K_1$ and the fixed compact curve family.
On a small open neighborhood $\mathcal U$ of $\mathcal Q$, the strict
margin remains at least $\eta/2$.  The same series therefore converges
locally uniformly for one $h_0$; its analytic terms give analytic
parameter dependence.  Restriction to $K_0\Subset K_1$
gives separate relative errors $1+O(h)$ for the two modes.  Their
logarithms near one give the bounded remainders $r_\pm$.

\textbf{Wronskian estimate.}
Apply the formula on a compact subset containing $K_0$ and its
$h$-translate.  Taylor expansion gives the one-step factors
$e^{\pm ip_\alpha}(1+O(h))$.  The two errors are separate, but
$\lvert \sin p_\alpha\rvert \ge c_0$ prevents cancellation of the leading terms:
\begin{align}
  W[\psi_+,\psi_-]
    &=\frac{e^{ip_\alpha}(1+O(h))-e^{-ip_\alpha}(1+O(h))}
            {\sin p_\alpha}
      =2i\bigl(1+O(h)\bigr),
\end{align}
which proves \eqref{eq:uniform-wkb-wronskian}.
\end{proof}
\noindent For background, see \cite{B-F:94,MR3501077,Fed2017}.
The restriction $K_0\Subset K_1$ is essential: the estimates are
uniform on compact subsets, not up to the boundary.
\subsection{Horizontal continuation}
\label{sec:horizontal-continuation}

The following result carries an existing WKB estimate along a fixed
horizontal segment.  The sign condition is exactly the condition
that its exponential grows in the chosen direction.

\begin{theorem}[Uniform finite-segment continuation]
  \label{thm:finite-segment-continuation}
  Let $\mathcal Q\subset\C^d$ be compact and let
  $I=[w_0,w_1]$ be a horizontal segment with $\Re w_1>\Re w_0$.
  For $\rho>0$, put
  $I_\rho\coloneq\{z:\operatorname{dist}(z,I)<\rho\}$.
  Suppose $v_\alpha$ is analytic in $(\alpha,z)$ near
  $\mathcal Q\times\overline{I_\rho}$.
  Suppose that on this neighborhood the
  branches $p_\alpha$ and $\sqrt{\sin p_\alpha}$ are analytic in
  $(\alpha,z)$ and, for some $\eta>0$,
  \begin{align}
     \Im p_\alpha(z)\le-\eta,
     \qquad (\alpha,z)\in\mathcal Q\times I_\rho.
     \label{eq:horizontal-growth-margin}
  \end{align}
  Let $V_0\Subset I_\rho$ be a fixed neighborhood of $w_0$.  Suppose
  $\psi_+(\alpha,h;z)$ is an exact solution on a neighborhood of
  $\overline V_0$, analytic in $(\alpha,z)$, with the plus asymptotic in
  \eqref{eq:asymp:wkb} there, uniformly in $\alpha$, using one action
  $S_\alpha$ that satisfies $S_\alpha'=p_\alpha$.
  Then there are $0<\rho'<\rho$ and $h_0,C>0$ such that, for
  $0<h<h_0$, $\psi_+$ continues analytically to $I_{\rho'}$ and
  \begin{align}
    \psi_+(\alpha,h;z)
      =\frac{e^{iS_\alpha(z)/h}}{\sqrt{\sin p_\alpha(z)}}
       \exp\bigl(h r_+(\alpha,h;z)\bigr),
    \qquad \abs{r_+}\le C,
    \label{eq:uniform-horizontal-continuation}
  \end{align}
  uniformly for $\alpha\in\mathcal Q$ and $z\in I_{\rho'}$.
\end{theorem}
\begin{proof}
We make the continuation principle
\cite[Theorem~3.1]{FedotovKlopp2021Meromorphic} quantitative, retaining
the $O(h)$ error and uniform parameter dependence.

\textbf{A uniform chain of canonical domains.}
First note that \eqref{eq:horizontal-growth-margin} keeps $I_\rho$ regular:
\begin{align}
  \abs{\sin p}^2
  =\sin^2(\Re p)+\sinh^2(\Im p)\ge\sinh^2\eta.
  \label{eq:horizontal-regularity-margin}
\end{align}
Fix $\alpha_0\in\mathcal Q$.  The local strict canonical-domain
construction, compactness of $I$, and the initial neighborhood $V_0$ give a finite
left-to-right chain
\begin{align}
  D'_0\Subset D_0,\ldots,D'_J\Subset D_J\Subset I_\rho.
\end{align}
Here every $D_j$ is canonical for $p_{\alpha_0}$ with strict inequalities, $D'_0\subset V_0$,
and the smaller domains cover $I$.  Restrict these domains to
a common narrow horizontal band. On each horizontal slice used in
$D'_j$, choose the overlap $D'_{j-1}\cap D'_j$ to contain an interval
of fixed positive length to the left of the points being reached.
Strictness persists for $\alpha$ in
a neighborhood of $\alpha_0$.  Finitely many such parameter neighborhoods
cover $\mathcal Q$.  Taking the smallest canonical margins, distances from $\overline{D'_j}$
to $\partial D_j$, and overlap widths,
and the largest chain length, makes all estimates below uniform.  For a fixed
$\alpha$, choose one chain whose parameter neighborhood contains it.
\Cref{thm:uniform-wkb} gives, on each pair of nested domains
$D'_j\Subset D_j$, solutions
$u_{+,j},u_{-,j}$ with the same branches and the same action $S_\alpha$ as
the initial solution.  For all small $h$, the recurrence itself first
continues the exact solution through the chain; the positive distances from $\overline{D'_j}$ to $\partial D_j$ ensure that
every use of $z\mapsto z+h$ stays in the corresponding outer domain.

\textbf{One overlap.}
Suppose inductively that $\psi_+$ has the desired formula on $D'_{j-1}$.
On $D_j$ write
\begin{align}
  \psi_+=a_j u_{+,j}+b_j u_{-,j}.
  \label{eq:rectangle-basis-decomposition}
\end{align}
The Wronskian formula \eqref{eq:wronskian-expansion} shows that $a_j$
and $b_j$ are
$h$-periodic. On $D'_{j-1}\cap D'_j$, the uniform WKB formulas at $z$ and
$z+h$, together with \eqref{eq:uniform-wkb-wronskian}, give
\begin{align}
  a_j=1+O(h),
  \qquad
  b_j=O(h)e^{2iS_\alpha/h}.
  \label{eq:rectangle-coefficients}
\end{align}
The $O(h)$ errors here come from \Cref{thm:uniform-wkb}; the cited
fixed-parameter continuation theorem states only the corresponding
$o(1)$ conclusion.
\textbf{Propagation in the growing direction.}
For a point $z$ to the right of the overlap, choose a point
$\widetilde z$ in the overlap with $z-\widetilde z\in h\Z$; this is
possible for all small $h$ because the overlap has fixed positive width.
The periodicity in \eqref{eq:rectangle-coefficients} makes the unwanted
term in \eqref{eq:rectangle-basis-decomposition}, relative to the plus
term, at most a constant times
\begin{align}
 h\left\lvert \exp\left(-\frac{2i}{h}
           \int_{\widetilde z}^{z}p_\alpha(\zeta)\dd\zeta\right)\right\rvert 
 &=h\exp\left(\frac2h
       \int_{\Re\widetilde z}^{\Re z}
          \Im p_\alpha(x+i\Im z)\dd x\right)
 \le h.
 \label{eq:rectangle-growing-direction}
\end{align}
Induction over the uniformly finite chain retains this relative $O(h)$
error.  Taking the analytic logarithm of its multiplier near one proves
\eqref{eq:uniform-horizontal-continuation}.  The exact continuation
throughout is the recurrence
$\psi(z+h)=-\psi(z-h)-v_\alpha(z)\psi(z)$.
\end{proof}
The other three cases are obtained by changing the sign of the phase or the
direction.  Explicitly, the plus solution propagates to the right when
$\Im p<0$ and to the left when $\Im p>0$; the minus solution propagates to
the right when $\Im p>0$ and to the left when $\Im p<0$.

The complex WKB theory for second-order difference equations is
closely related to the corresponding theory for adiabatic perturbations of a
periodic operator developed in~\cite{MR1903839,MR1858917,MR2097997}.
\subsection{Normalization and growth at imaginary infinity}
\label{app:growth}

This subsection proves the normalization and growth estimates used in
\Cref{sec:wkb-method-second}.

In this and the next subsection, fix $\Im\lambda>0$ and a compact set
$\mathcal K\Subset\{\Im E>-\Im\lambda\}$.
We use the common canonical strip $K_\delta$ from
\Cref{prop:unbounded-canonical-strip} and fix $0<\delta'<\delta$.
We begin by establishing uniform asymptotic estimates for exact
solutions of \eqref{eq:main-cos} in $K_{\delta'}$.

\begin{proof}[Proof of \Cref{prop:entire-WKB-normalization}]

\Cref{thm:uniform-wkb} gives estimates on bounded domains, but does
not immediately give an error bound uniform throughout the infinite
strip $K_{\delta'}$. For this we use the unbounded WKB estimates of
Fedotov and Shchetka
\cite[Theorems~2.2--2.3]{FedotovShchetka2022Berry}.

Write the analytic finite-difference equation as
\begin{align}
 \begin{pmatrix}
  \psi(z+h)\\
  \psi(z)
 \end{pmatrix}
 &=
 M_E(z)
 \begin{pmatrix}
  \psi(z)\\
  \psi(z-h)
 \end{pmatrix},
 &
 M_E(z)&\coloneq
 \begin{pmatrix}
  2(E-\lambda\cos z)&-1\\
  1&0
 \end{pmatrix}.
\end{align}
The matrix $M_E$ satisfies the hypotheses of
\cite[Theorems~2.2--2.3]{FedotovShchetka2022Berry}:
$\det M_E=1$, $M_{12}M_{21}=-1$, and the coefficients of
$e^{\pm iz}$ in $\tr M_E$ are both $-\lambda\ne0$.
Also, $M_{22}=0$. The strip $K_\delta$ is an unbounded canonical
domain, and $K_{\delta'}$ stays a positive distance from its boundary,
as required by Theorem~2.3.

Theorem~2.2 of \cite{FedotovShchetka2022Berry} constructs exact
solutions with locally uniform WKB asymptotics on the unbounded
canonical domain. Theorem~2.3 controls these same solutions at both
imaginary ends. Using the scalar normalization in
\cite[(2.14)]{FedotovShchetka2022Berry}, the plus solution satisfies
\begin{align}
 \phi_+^0(z;E,\lambda,h)
 &=
 \frac{1}{\sqrt{\sin p(z;E)}}
 \exp\left(\frac{i}{h}S_E(z)+g_E(z,h)\right)
 \bigl(1+r_E(z,h)\bigr),\\
 S_E(z)&\coloneq\int_\pi^z p(w;E)\dd w,
\end{align}
for $z\in K_{\delta'}$ with $\lvert \Im z\rvert \ge Y$. Here $Y$ is fixed,
$r_E=O(h)$ uniformly at both ends, and the theorem gives
\begin{align}
 \lvert g_E(z,h)\rvert &\le Ch(1+\lvert z\rvert ).
\end{align}
If $M_{22}/M_{11}\to0$ at both imaginary ends, the stronger conclusion is
\begin{align}
 \lvert g_E(z,h)\rvert &\le Ch.
\end{align}
For our matrix, $M_{22}=0$, so this stronger estimate applies.
In particular, $e^{g_E}=1+O(h)$ uniformly, and hence
\begin{align}
 \phi_+^0(z;E,\lambda,h)
 &=
 \frac{e^{iS_E(z)/h}}{\sqrt{\sin p(z;E)}}\bigl(1+O(h)\bigr)
 \label{eq:rewrite-entire-strip-WKB}
\end{align}
uniformly at both ends of $K_{\delta'}$. The locally uniform estimate of
\cite[Theorem~2.2]{FedotovShchetka2022Berry} for this same solution
covers the remaining bounded portion of the strip.

We also require these bounds to be uniform in $E$.
\Cref{prop:unbounded-canonical-strip} supplies the same regular
canonical strip for all $E$ near $\mathcal K$.
Moreover, on every narrower strip, the complex momentum satisfies
\begin{align}
 p(z;E)
 &=-\operatorname{sgn}(\Im z)(z-\pi)-i\log\lambda
   +O(e^{-\lvert \Im z\rvert }),
 \label{eq:momentum-at-infinity}
\end{align}
as $\lvert \Im z\rvert \to\infty$, uniformly for $E$ in the fixed compact set.
The leading terms are
independent of $E$, and the remaining bounded portion of the strip
stays uniformly away from the turning points. Consequently, after
shrinking the parameter neighborhood, there is $c_0>0$ such that
$c_0\le\Re p(z;E)\le\pi-c_0$ and
$\lvert\sin p(z;E)\rvert\ge c_0e^{\lvert\Im z\rvert}$
on $K_{\delta'}$, uniformly in $E$.
These bounds allow the contours and constants in the construction of
\cite[Section~5]{FedotovShchetka2022Berry} to be chosen uniformly.
With fixed contours and analytic branches, its integral equations
depend analytically on $E$. Their uniformly convergent successive
approximations therefore give a solution analytic in $E$.

For $h$ sufficiently small compared with the strip width, the recurrence
\begin{align}
 \phi_+^0(z+h)
 &=2(E-\lambda\cos z)\phi_+^0(z)-\phi_+^0(z-h)
\end{align}
and its backward form extend the solution horizontally to all of $\C$. These extensions agree on their overlaps by the
difference equation and preserve analytic dependence on $E$.
Thus $\phi_+^0$ is entire in $z$.
Finally, $\cos(2\pi-z)=\cos z$, so the chosen momentum branch
and its action satisfy
\begin{align}
 p(2\pi-z;E)&=p(z;E),&
 S_E(2\pi-z)&=-S_E(z).
\end{align}
The chosen square root of $\sin p(z;E)$ is also invariant under this
reflection: the two branches agree at $z=\pi$.
Thus the substitution $z\mapsto2\pi-z$ gives the minus WKB estimate
for $\phi_-^0(z)=\phi_+^0(2\pi-z)$ on $K_{\delta'}$.

\end{proof}

\subsection{Wronskian coefficients and uniqueness}
\label{app:wronskian-growth}

We next estimate the Wronskians $\mathcal{A} = W[\phi_-^0,\phi_+^0]$ and $\mathcal{B} = W[\phi_-^0,\phi_+^{2\pi}]$.
Broadly speaking, this is done by noting that such Wronskians are entire and $h$-periodic in the real direction, and thus admit a Fourier expansion.
Growth estimates at both imaginary ends of $K_{\delta'}$ show that only the constant Fourier coefficients of $\mathcal{A}$ and $\mathcal{B}$ can be nonzero.

\begin{proof}[Proof of \Cref{prop:constant-AB}]
We first prove that $\mathcal{A}$ is constant.
Take $h<\delta'$ and $\pi\le\Re z\le\pi+h$, so that both $z$
and $z-h$ lie in $K_{\delta'}$.

By \Cref{prop:entire-WKB-normalization}, for $z$ and $z-h$ in $K_{\delta'}$ 
\begin{align}
    \phi_-^0 (z) \phi_+^0(z-h) = \frac{\exp(- \tfrac{i}{h}\int_{z-h}^{z} p(\zeta)\dd\zeta)}{\sqrt{\sin p(z)}\sqrt{\sin p(z-h)}} (1+O(h)).
\end{align}
Put $t=\lvert \Im z\rvert $. Along a horizontal segment from $z-h$ to $z$, as $t\to\infty$, we can write
\begin{align}
    \Im p (\zeta) = -t - \log \lvert \lambda \rvert  + O(e^{-t})
\end{align}
so that
\begin{align}
 \left\lvert \exp\left(-\frac{i}{h}\int_{z-h}^{z}p(\zeta)\dd\zeta\right)\right\rvert 
 &=\exp\left(\frac1h\int_{z-h}^{z}\Im p(\zeta)\dd\zeta\right)
 \le Ce^{-t}.
\end{align}
Also, $\lvert \sin p(z)\rvert \asymp\lvert \sin p(z-h)\rvert \asymp e^t$, so the
amplitude contributes a further factor $O(e^{-t})$. Hence
\begin{align}
 \lvert \phi_-^0(z)\phi_+^0(z-h)\rvert &\le Ce^{-2t}.
 \label{eq:rewrite-A-first-product}
\end{align}
The other term in the Wronskian has the opposite phase:
\begin{align}
 \phi_-^0(z-h)\phi_+^0(z)
 &=\frac{\exp\left(\frac{i}{h}\int_{z-h}^{z}p(\zeta)\dd\zeta\right)}
 {\sqrt{\sin p(z-h)}\sqrt{\sin p(z)}}(1+O(h)).
\end{align}
Its exponential is bounded by $Ce^t$, while its amplitude is again
$O(e^{-t})$. Consequently,
\begin{align}
 \lvert \phi_-^0(z-h)\phi_+^0(z)\rvert &\le C.
 \label{eq:rewrite-A-second-product}
\end{align}
These estimates hold at both imaginary ends. In particular,
\begin{align}
 \mathcal{A}(z-h)
 &=\phi_-^0(z)\phi_+^0(z-h)-\phi_-^0(z-h)\phi_+^0(z)
\end{align}
is bounded there. The remaining bounded-height portion of the period
cell is compact, so $\mathcal{A}$ is bounded on the whole cell.

Because $\mathcal{A}$ is entire and $h$-periodic, it is bounded everywhere and hence constant.
A local Wronskian computation at $z=\pi$ determines its value.
Since $S_E(\pi)=0$ and
$S_E(\pi+h)=hp(\pi)+O(h^2)$, the WKB estimates give
\begin{align}
 \mathcal{A}
 &=\phi_-^0(\pi+h)\phi_+^0(\pi)
   -\phi_-^0(\pi)\phi_+^0(\pi+h)\\
 &=\frac{e^{-ip(\pi)}-e^{ip(\pi)}}{\sin p(\pi)}+O(h)
 =-2i(1+O(h)).
\end{align}
Multiplying both solutions by the analytic branch of
$\sqrt{(-2i)/\mathcal{A}}=1+O(h)$ close to one makes $\mathcal{A}=-2i$ exactly.
This preserves the WKB estimates and the relation
$\phi_-^0(z)=\phi_+^0(2\pi-z)$.

To estimate $\mathcal{B} = W[\phi_-^0,\phi_+^{2\pi}]$, we cannot use the same estimate as we do not have the same asymptotic knowledge of $\phi_+^{2\pi}$ in $K_{\delta ' }$.
Work on the period cell $\pi\le\Re z\le\pi+h$, with $h<\delta'$,
so that $z$ and $z-h$ remain in $K_{\delta'}$.
To estimate the translated terms we iterate the recurrence from the finite difference equation backward.
Since $Nh=2\pi$, we can write
\begin{align}
    \begin{pmatrix}
        \phi_+^0(z - 2\pi)\\ \phi_+^0(z- 2\pi -h   )
    \end{pmatrix}
    &=T_N(z)\begin{pmatrix}\phi_+^0(z)\\\phi_+^0(z-h)\end{pmatrix},
    \label{eq:rewrite-B-transfer}\\
    T_N(z)&\coloneq Q_N(z)\cdots Q_1(z),&
    Q_j(z)&\coloneq
    \begin{pmatrix}0&1\\-1&2(E-\lambda\cos(z-jh))\end{pmatrix}.
\end{align}
The entries of $T_N$ are polynomials in the values
$2(E-\lambda\cos(z-jh))$. For $N\ge2$, induction using $T_{j+1}=Q_{j+1}T_j$ gives the entrywise degree bounds
\begin{align}
    \deg T_N&\le\begin{pmatrix}N-2&N-1\\N-1&N\end{pmatrix}.
\end{align}

Writing $T_N=(t_{ij})$, and using that $2(E-\lambda\cos(z-jh))=O(e^t)$ uniformly in $j$,
the degree bounds give
\begin{align}
    \lvert t_{11}\rvert &\le C_h e^{(N-2)t},&
    \lvert t_{12}\rvert +\lvert t_{21}\rvert &\le C_h e^{(N-1)t},&
    \lvert t_{22}\rvert &\le C_h e^{Nt}.
    \label{eq:rewrite-B-entry-bounds}
\end{align}
Here $C_h$ is independent of $t$.

Expanding the Wronskian and using \eqref{eq:rewrite-B-transfer}, we get
\begin{align}
 \mathcal{B}(z-h)
 ={}&\phi_-^0(z)\phi_+^0(z-h-2\pi)
     -\phi_-^0(z-h)\phi_+^0(z-2\pi)\\
 ={}&t_{21}\phi_-^0(z)\phi_+^0(z)
     +t_{22}\phi_-^0(z)\phi_+^0(z-h)\\
 &-t_{11}\phi_-^0(z-h)\phi_+^0(z)
  -t_{12}\phi_-^0(z-h)\phi_+^0(z-h).
\end{align}
Using the WKB product estimates \eqref{eq:rewrite-A-first-product} and \eqref{eq:rewrite-A-second-product}, we obtain
\begin{align}
 \lvert \phi_-^0(z-h)\phi_+^0(z)\rvert &\le C,&
 \lvert \phi_-^0(z)\phi_+^0(z-h)\rvert &\le Ce^{-2t}.
\end{align}
While writing $\phi_-^0(z) \phi_+^0(z) = (\sin p (z))^{-1} (1 + O(h)) = O(e^{-t})$ gives us
\begin{align}
     \lvert \phi_-^0(z)\phi_+^0(z)\rvert &\le Ce^{-t},& \lvert \phi_-^0(z-h)\phi_+^0(z-h)\rvert &\le Ce^{-t}.
\end{align}
Combining these estimates with
\eqref{eq:rewrite-B-entry-bounds}, each of the four terms is bounded by
$C_h e^{(N-2)t}$. By periodicity, and after increasing $C_h$ to cover
bounded heights, we obtain
\begin{align}
 \lvert \mathcal{B}(x+iy)\rvert &\le C_h e^{(N-2)\lvert y\rvert },\qquad x,y\in\R.
 \label{eq:rewrite-B-growth}
\end{align}

Now, because $\mathcal{B}$ is entire and $h$-periodic in the real direction,
it has the Fourier expansion
\begin{align}
 \mathcal{B}(z)&=\sum_{k\in\Z}\widehat{\mathcal{B}}_k e^{ikNz},&
 \widehat{\mathcal{B}}_k&=\frac1h\int_\pi^{\pi+h}
 \mathcal{B}(x+iy)e^{-ikN(x+iy)}\dd x,
\end{align}
where $y\in\R$ is arbitrary. Cauchy's theorem makes the integral
independent of $y$, since the straight vertical sides of the rectangle cancel
by periodicity.
For $k\ge1$, the growth bound \eqref{eq:rewrite-B-growth}, with $y=-t$, gives
\begin{align}
 \lvert \widehat{\mathcal{B}}_k\rvert &\le C_h e^{(N-2-kN)t}\xrightarrow{t\to\infty}0.
\end{align}
For $k\le-1$, taking $y=t$ instead gives
\begin{align}
 \lvert \widehat{\mathcal{B}}_k\rvert &\le C_h e^{(N-2+kN)t}\xrightarrow{t\to\infty}0.
\end{align}
Thus all coefficients with $k\ne0$ vanish, and $\mathcal{B}$ is constant.
The limits keep $h$ fixed, so the dependence of $C_h$ on $h$ is harmless.
\end{proof}

We next prove that $\mathcal{C}=W[\phi_+^0,\phi_+^{2\pi}]$ has only the
Fourier modes $0$ and $1$, so $\mathcal{C}(z,E)=c(E,h)+a(E,h)e^{2\pi iz/h}$,
with $a$ and $c$ analytic in $E$ near $\mathcal K$.

\begin{proof}[Proof of \Cref{prop:exact-C-modes}]
We use the same recurrence and Fourier argument as in the proof
of \Cref{prop:constant-AB}. The difference is that the Wronskian
has two plus solutions.

Let $\theta\coloneq\arg\lambda\in(0,\pi)$.
Throughout, $h$ is fixed and $t\to\infty$.
We have
\begin{align}
    \mathcal{C}(z-h)
    &=\phi_+^0(z)\phi_+^0(z-h-2\pi)
    -\phi_+^0(z-h)\phi_+^0(z-2\pi).
    \label{eq:C_append_wronskian}
\end{align}
The goal is to prove
\begin{align}
    \lvert \mathcal{C}(x+it)\rvert
    &\le C_h e^{(N-2\theta/h)t},&
    \lvert \mathcal{C}(x-it)\rvert
    &\le C_h e^{(N+2\theta/h)t},
    \qquad x\in\R.
    \label{eq:C_goal}
\end{align}

By \eqref{eq:momentum-at-infinity} and
\eqref{eq:rewrite-entire-strip-WKB}, for $-3h/2\le x\le h/2$,
\begin{align}
    \lvert\phi_+^0(\pi+x\pm it)\rvert
    &\le C_h
    \exp\left(\left(\frac{x\mp\theta}{h}-\frac12\right)t\right).
    \label{eq:281}
\end{align}
This interval lies inside the WKB strip for sufficiently small $h$.
Now write $z=\pi+x\pm it$ with $\lvert x\rvert\le h/2$.
Applying \eqref{eq:281} at $x$ and $x-h$ gives
\begin{align}
    \lvert\phi_+^0(z)\rvert
    &\le C_h e^{\mp\theta t/h},&
    \lvert\phi_+^0(z-h)\rvert
    &\le C_h e^{-t}e^{\mp\theta t/h}.
\end{align}
Consequently,
\begin{align}
    \lvert\phi_+^0(z)\rvert+\lvert\phi_+^0(z-h)\rvert
    &\le C_h e^{\mp\theta t/h}.
\end{align}

The recurrence bounds
\eqref{eq:rewrite-B-transfer} and
\eqref{eq:rewrite-B-entry-bounds} give
\begin{align}
    &\lvert\phi_+^0(z-2\pi)\rvert
    +\lvert\phi_+^0(z-h-2\pi)\rvert\le C_h e^{Nt}
    \bigl(\lvert\phi_+^0(z)\rvert
    +\lvert\phi_+^0(z-h)\rvert\bigr).
\end{align}
Substituting these bounds into \eqref{eq:C_append_wronskian},
we obtain
\begin{align}
    \lvert \mathcal{C}(z-h)\rvert
    &\le C_h e^{Nt}
    \bigl(\lvert\phi_+^0(z)\rvert
          +\lvert\phi_+^0(z-h)\rvert\bigr)^2\le C_h e^{(N\mp2\theta/h)t}.
\end{align}
Since $\lvert x\rvert\le h/2$ covers a full period cell,
$h$-periodicity gives \eqref{eq:C_goal}.

As in the proof for $\mathcal{B}$, shift the contour defining each
Fourier coefficient to obtain
\begin{align}
    \lvert\widehat{\mathcal{C}}_k\rvert
    &\le C_h e^{((k+1)N-2\theta/h)t}
    \xrightarrow{t\to\infty}0,
    &&k\le-1,\\
    \lvert\widehat{\mathcal{C}}_k\rvert
    &\le C_h e^{((1-k)N+2\theta/h)t}
    \xrightarrow{t\to\infty}0,
    &&k\ge2.
\end{align}
Both limits follow from $0<2\theta/h<N$.
Hence
\begin{align}
    \mathcal{C}(z,E)&=c(E,h)+a(E,h)e^{2\pi iz/h},
    &c&\coloneq\widehat{\mathcal{C}}_0,
    &a&\coloneq\widehat{\mathcal{C}}_1.
\end{align}
The coefficients are analytic in $E$, since they are integrals
over a fixed period cell of the analytic function $\mathcal{C}(z,E)$
against the corresponding Fourier factors.
\end{proof}

We lastly prove uniqueness of solutions modulo a constant factor.

\begin{proof}[Proof of \Cref{cor:entire-normalization-unique}]
Let $\phi_+^0$ and $\tilde\phi_+^0$ be two solutions satisfying
the plus WKB estimate of \Cref{prop:entire-WKB-normalization}.
Since $\mathcal{A}=W[\phi_-^0,\phi_+^0]$ is a nonzero constant, the
Wronskian expansion \eqref{eq:wronskian-expansion} gives
\begin{align}
    \tilde\phi_+^0(z)
    &=\alpha(z)\phi_+^0(z)+\beta(z)\phi_-^0(z),\\
    \alpha(z)
    &\coloneq\frac{W[\phi_-^0,\tilde\phi_+^0](z)}{\mathcal{A}},
    &\beta(z)
    &\coloneq\frac{W[\tilde\phi_+^0,\phi_+^0](z)}{\mathcal{A}}.
\end{align}
Both coefficients are entire and $h$-periodic, since the
Wronskians have these properties.

We first show that $\beta=0$. Expanding the Wronskian gives
\begin{align}
    \mathcal{A}\beta(z-h)
    &=\tilde\phi_+^0(z)\phi_+^0(z-h)
    -\tilde\phi_+^0(z-h)\phi_+^0(z).
\end{align}
Applying \eqref{eq:281} to both plus solutions and using
periodicity, we obtain (by a similar computation as above)
\begin{align}
    \lvert\beta(x+it)\rvert
    &\le C_h e^{-(2\theta/h+1)t},&
    \lvert\beta(x-it)\rvert
    &\le C_h e^{(2\theta/h-1)t},
    \qquad x\in\R,
\end{align}
where $\theta=\arg\lambda\in(0,\pi)$.
The same Fourier argument excludes all modes $k\le0$
using the upper bound and all modes $k\ge1$ using the
lower bound, since $0<2\theta/h<N$.
Thus $\beta=0$.

We now have $\tilde\phi_+^0=\alpha\phi_+^0$.
Comparing the two uniform WKB estimates
\eqref{eq:rewrite-entire-strip-WKB} gives
\begin{align}
    \alpha(z)
    &=\frac{\tilde\phi_+^0(z)}{\phi_+^0(z)}
    =1+O(h)
\end{align}
throughout a period cell contained in $K_{\delta'}$. The quotient is well defined for small $h$, since the uniform WKB estimate makes $\phi_+^0$ nonzero there.
By periodicity, $\alpha$ is bounded on $\C$ and hence constant (as it is entire).

Finally, assume that $W[\phi_-^0,\phi_+^0]=-2i$ and
$W[\tilde\phi_-^0,\tilde\phi_+^0]=-2i$, where each minus solution
is defined by reflection of its plus partner.
Since $\alpha$ is constant,
\begin{align}
    \tilde\phi_-^0(z)
    &=\tilde\phi_+^0(2\pi-z)
    =\alpha\phi_+^0(2\pi-z)
    =\alpha\phi_-^0(z).
\end{align}
Thus
\begin{align}
    -2i = W[\tilde \phi_-^0 ,\tilde \phi_+^0 ] = \alpha^2 W[\phi_-^0,\phi_+^0] = \alpha^2 \mathcal{A} = -2i\alpha^2.
\end{align}
Hence $\alpha^2=1$. Since $\alpha=1+O(h)$, we have
$\alpha=1$ for sufficiently small $h$, proving uniqueness.
\end{proof}

\subsection{Comparison regions}
\label{app:geometry}

We prove \Cref{prop:uniform-WKB-comparison-geometry}.  The interior
construction has one global issue: a canonical curve near the
real-momentum arc must actually reach the normalization strip at both
ends.  Negative level sets give barriers that settle this issue.
The open-side and exterior comparisons instead use local endpoint
bases and horizontal continuation.

\subsubsection{Interior energies}
\label{sec:uniform-geometry-interior}

\begin{proof}[Proof of part \textup{(i)} of
\Cref{prop:uniform-WKB-comparison-geometry}]
\textbf{The arc and its right-hand region.}
Write $\lambda=\lambda_r+i\lambda_i$, $\lambda_i>0$, and
$E=s+\lambda t$, where $s,t\in(-1,1)$.  In
\eqref{eq:action-arc}, put $r\coloneq s-u$ and $T(u)\coloneq\lambda\sin z(u)$.
Then
\begin{align}
 T^2&=(1-t^2)\lambda^2-2tr\lambda-r^2,\\
 \Im(T^2)&=2\lambda_i\bigl((1-t^2)\lambda_r-tr\bigr).
\end{align}
This square never lies on the nonnegative real axis.  If $t\ne0$
and its imaginary part vanishes, its real part is
$-(1-t^2)(\lambda_i^2+\lambda_r^2/t^2)<0$; for $t=0$ the same
conclusion is immediate.  Choose the square root with positive
imaginary part.  It agrees with $\lambda\sqrt{1-t^2}$ at $u=s$, so
\begin{align}
 \frac{\dd z}{\dd u}&=\frac1T,&
 \frac{\dd\Im z}{\dd u}&=-\frac{\Im T}{\lvert T\rvert ^2}<0.
 \label{eq:general-arc-vertical}
\end{align}
Thus the real-momentum arc is a vertical graph.  As its height
increases from $\Im\zeta_+<0$ to $\Im\zeta_->0$, its momentum
increases from zero to $\pi$.

Let $D$ be the region between this arc and the line $\Re z=\pi$,
extended above and below the endpoints as follows.  From each turning
point run a horizontal cut to the imaginary axis; beyond these cuts,
use the imaginary axis as the left boundary of $D$.
The region is regular and horizontally connected.  The branch
continued from the $\pi$-strip has $\Im p<0$ throughout $D$.
On the two outer left boundary rays,
\begin{align}
 \Re p(iy)&<0 &&(y<\Im\zeta_+),&
 \Re p(iy)&>\pi &&(y>\Im\zeta_-).
 \label{eq:general-boundary-sheet}
\end{align}
To see the branch signs, extend the parameter $u$ over the real line
in \eqref{eq:general-arc-vertical}.  On each horizontal line there is
exactly one point in $0<\Re z<\pi$ where $f_E(z)$ is real. Below the lower endpoint
it has $f_E>1$, and above the upper endpoint it has $f_E<-1$.
Continuation from the right crosses $\Re p=0$ or $\Re p=\pi$,
respectively, while retaining $\Im p<0$.

\textbf{Negative levels cannot escape to the left.}
Choose a regular point $z_\vartheta$ of the arc with
$p(z_\vartheta)=\vartheta\in(0,\pi)$, and set
\begin{align}
 H_\vartheta(z)
 &\coloneq\Im\int_{z_\vartheta}^z(p(w)-\vartheta)\dd w.
 \label{eq:separated-corridor-level-function}
\end{align}
The entire left boundary of $D$ has $H_\vartheta\geq0$.
On the arc, $\dd H_\vartheta/\dd y=p-\vartheta$, whose sign changes
from negative to positive at $z_\vartheta$.  Thus its minimum there
is zero.  Moving left on either horizontal cut increases
$H_\vartheta$, since $\partial_xH_\vartheta=\Im p<0$.
On the outer boundary rays, \eqref{eq:general-boundary-sheet} shows
that $H_\vartheta$ increases when moving away from the endpoints.

A small negative level $H_\vartheta=-v$, $v>0$, therefore passes
just to the right of $z_\vartheta$.  Since
$\partial_xH_\vartheta<0$ in $D$, its component is a vertical graph
$x=x(y)$, and
\begin{align}
 x'(y)&=\frac{\Re p-\vartheta}{-\Im p},&
 \Re p+\Im p\,x'(y)&=\vartheta.
 \label{eq:separated-level-canonical-derivative}
\end{align}
This is the reason for using negative levels: they are canonical with
a fixed strict margin, and their boundary values prevent them from
hitting either a turning point or the left boundary.  At a finite
height the graph can stop only on $\Re z=\pi$, inside the original
strip.  Otherwise it continues to imaginary infinity.

\textbf{Joining both ends to the original strip.}
We also check the ends that do not reach the strip at finite height.
Put $\theta_\lambda=\arg\lambda\in(0,\pi)$.
Uniformly for bounded $x$ on the continued branch,
\begin{align}
 p(x+iy)&=\pi-x+\theta_\lambda-i(y+\log\lvert \lambda\rvert )
             +O(e^{-y}) &&(y\to+\infty),\\
 p(x+iy)&=x-\pi+\theta_\lambda+i(y-\log\lvert \lambda\rvert )
             +O(e^{y}) &&(y\to-\infty).
\end{align}
Equation \eqref{eq:separated-level-canonical-derivative} consequently
gives
\begin{align}
 x(y)&=\pi+\theta_\lambda-\vartheta+O(1/y)
             &&(y\to+\infty),\\
 x(y)&=\pi-\theta_\lambda+\vartheta+O(1/\lvert y\rvert )
             &&(y\to-\infty),&
 x'(y)&=O(\lvert y\rvert ^{-2}).
\end{align}
Only an end whose limiting line is at or to the left of $\pi$ can
occur before the graph enters the strip.

Such an end can be bent into the strip while preserving the canonical
inequalities. Choose $x_1<x_2<\pi$ inside the normalization strip and
a smooth function $\chi_0\colon\R\to[0,1]$ equal to one for $x\leq x_1$
and zero for $x\geq x_2$. At sufficiently large $\lvert y\rvert $, the momentum
expansions give, throughout the region between the outer curve and the strip,
\begin{align}
 \eta_0\leq\Re p\leq\pi-\eta_0,
 \qquad \lvert \Im p\rvert \leq C_0(1+\lvert y\rvert ).
\end{align}
For $0<\varepsilon<\eta_0/(2C_0)$, choose curves satisfying
\begin{align}
 x'(y)&=\frac{\operatorname{sgn}(y)\varepsilon\chi_0(x)}{1+\lvert y\rvert }.
 \label{eq:tail-curves-into-strip}
\end{align}
Then $\eta_0/2\leq\Re p+\Im p\,x'\leq\pi-\eta_0/2$.
Writing $t=\lvert y\rvert $, as long as $x\leq x_1$ we have
\begin{align}
 \frac{\dd x}{\dd t}&=\frac{\varepsilon}{1+t},
 &x(t)-x(t_0)&=\varepsilon\log\frac{1+t}{1+t_0}.
\end{align}
Thus every curve enters the strip at finite height at both ends.

To obtain a domain, use the full negative sublevel region to the right
of one fixed level, not just a thin neighborhood of that level curve.
Since $\partial_xH_\vartheta<0$, each nonempty horizontal slice is an
interval meeting the original strip. On the compact part, choose curves
satisfying \eqref{eq:separated-level-canonical-derivative}. Inside the
strip, join them to straight vertical lines by smoothly replacing the
right-hand side by zero. On finite transition bands at large
$\lvert y\rvert$, interpolate with the equation
\eqref{eq:tail-curves-into-strip}. The canonical inequalities are affine
in $x'$, so this interpolation preserves a strict margin wherever
both equations give canonical curves.
The bounded horizontal width gives a common height by which every
curve has entered the strip. Beyond that height, a smooth cutoff
makes $x'=0$, leaving straight vertical tails.
Choose the deformed outer boundary to satisfy the same differential
equation. Uniqueness keeps the curves ordered, and they fill the region
between this boundary and a straight vertical line inside the strip.
This gives a horizontally connected unbounded canonical domain.
Finite joins inside the strip can also be rounded using
\cite[Proposition~4.1 and Lemma~4.1]{MR2097997}.

\textbf{Crossing the regular arc and reaching zero.}
At $z_\vartheta$ one has $\Im p'(z_\vartheta)<0$, since $p$ increases
along the vertical arc.  Taylor expansion gives
\begin{align}
 H_\vartheta(z)
 &=\frac12\Im\bigl(p'(z_\vartheta)(z-z_\vartheta)^2\bigr)
      +O(\lvert z-z_\vartheta\rvert ^3).
\end{align}
The negative-level graphs approach this saddle from the right with
bounded derivatives $x'(y)$. A small local bend therefore carries a thin part of
the domain across the arc.  This preserves the strict inequalities:
in a sufficiently small crossing neighborhood, $\Re p$ stays in a
compact subinterval of $(0,\pi)$ and $\Im p$ is small.
The infinite ends and the narrower normalization strip are unchanged.
The adjacent canonical-domain principle consequently gives the WKB
estimate for the same solution $\phi_+^0$ across the arc.

Take $z_\vartheta=z_*=\arccos t$.  On the real segment immediately
to its left,
\begin{align}
 \Im f_E(x)&=(\Im\lambda)(t-\cos x)<0,
 \qquad \Im p(x)>0.
\end{align}
Start inside the crossing neighborhood, strictly to the left of
$z_*$, and use \Cref{thm:finite-segment-continuation} to reach zero.
The remaining segment is regular and has a strict growing-direction
margin.
\Cref{fig:negative-level-construction} distinguishes the negative levels
from the deformations used to join the strip and cross the regular arc.

\negativelevelsfigure

\textbf{Uniformity.}
Make all choices first at $E_*$.  The compact parts of the domains,
bends, and horizontal segments have positive distance from
the turning points and strict canonical or growing-direction margins.
Keep these curves fixed and shrink a complex energy neighborhood.
Regularity and the margins persist there.  On the infinite tails the
displayed momentum expansions are uniform, and changing $E$ changes
their remainders by $O(e^{-\lvert \Im z\rvert })$.  The unbounded-domain hypotheses
and a narrower normalization strip therefore persist as well.
This proves the locally uniform statement.
\end{proof}

\subsubsection{Horizontal comparisons at open sides and exterior energies}
\label{sec:uniform-geometry-open-side}
\label{sec:uniform-geometry-exterior}

\begin{proof}[Proof of part \textup{(ii)} of
\Cref{prop:uniform-WKB-comparison-geometry}]
\textbf{The representative open side.}
Let $E_*=s+\lambda$, $-1<s<1$.  The real period is regular, including
zero, where $p(0)=\arccos s\in(0,\pi)$.
For $0<x<2\pi$,
\begin{align}
 \Im(E_*-\lambda\cos x)
 &=(\Im\lambda)(1-\cos x)>0,
\end{align}
so the continued interior momentum has $\Im p(x)<0$.

Choose a small bounded canonical domain about zero.  For example,
a narrow lens swept out by curves close to straight vertical lines,
with common upper and lower endpoints, is canonical, since $\Re p$ is close to
$\arccos s\in(0,\pi)$ and $\Im p$ is small there.
Choose $\delta>0$ so that a smaller domain contains
$[-2\delta,2\delta]$.  The bounded WKB theorem gives an exact pair
$\psi_\pm^0$ there, and translation gives a pair near $2\pi$.

The compact middle interval $[\delta,2\pi-\delta]$ has
$\Im p\leq-\eta<0$ on a fixed thin neighborhood.  The finite-segment
continuation theorem carries the plus estimate from the left endpoint
rightward and the translated minus estimate from the right endpoint
leftward.  Both resulting overlaps can be chosen a fixed positive distance from
the boundaries of the local endpoint domains.  Moreover,
\begin{align}
 \Im\int_0^\delta p(x)\dd x&<0.
\end{align}
Fix $\delta$ and the domains before shrinking the complex energy
neighborhood.  Regularity, the strict middle-interval inequality,
and this strict integral inequality all persist.  The sign of
$\Im p$ close to zero need not persist, and is not used for the
horizontal continuation.

\textbf{Exterior energies.}
For an exterior energy, $E-\lambda[-1,1]$ is disjoint from $[-1,1]$.
The even periodic branch $p$ has $\Im p<0$ on a fixed tube about
the real period.  A local canonical domain about zero again gives
an exact WKB basis. If $\Re p(0)\notin(0,\pi)$,
choose a slanted central curve with
$\Re p+\Im p\,x'=\pi/2$ and a sufficiently thin lens around it;
$\Im p\ne0$ permits this choice.  A smaller lens contains a short
real interval.  The preceding horizontal continuation and endpoint
integral estimates apply without change.

In both cases the local solutions extend by the recurrence along a
horizontal strip containing the period, so their Wronskians are
$h$-periodic wherever compared.  All domains are fixed away from the
turning points.  Their strict margins give analytic, uniform WKB
estimates for nearby energies.
\end{proof}

\section{Fourier transformation and duality}
\label{app:fourier-normalization}

\subsection{The Fourier transform}

The two normalized solutions need not have convergent Fourier
integrals on the same straight vertical line. We first define a
single linear transform by continuing the integrals over the upper
and lower half-lines separately.

\begin{lemma}
\label{lem:fourier-transform}
Fix $\lambda\ne0$ and $h=2\pi/N$. Let $\psi$ be an entire solution
of \eqref{eq:main-cos} such that, for every bounded interval
$I\subset\R$, there are constants $C_I,c_I$ with
\begin{align}
 \lvert\psi(x+iy)\rvert&\le C_I e^{c_I\lvert y\rvert},
 &x&\in I,\quad y\in\R.
 \label{eq:fourier-transform-growth-class}
\end{align}
There is a linear transform $\mathcal F_h\psi$, entire in $\rho$,
which agrees with the upward-oriented integral
\begin{align}
 (\mathcal F_h\psi)(\rho)
 &=\int_{x-i\infty}^{x+i\infty}e^{-i\rho z/h}\psi(z)\dd z
 \label{eq:fourier-bilateral-integral}
\end{align}
whenever both half-line integrals converge on an open set of $\rho$.
It solves
\begin{align}
 \lambda\bigl((\mathcal F_h\psi)(\rho+h)
       +(\mathcal F_h\psi)(\rho-h)\bigr)
 +2\cos\rho\,(\mathcal F_h\psi)(\rho)
 &=2E(\mathcal F_h\psi)(\rho).
 \label{eq:fourier-dual-equation}
\end{align}
The following identities hold as identities of entire functions:
\begin{align}
 \mathcal F_h[\psi(\,\cdot-2\pi)](\rho)
 &=e^{-2\pi i\rho/h}(\mathcal F_h\psi)(\rho),\\
 \mathcal F_h[e^{2\pi i(\cdot)/h}\psi](\rho)
 &=(\mathcal F_h\psi)(\rho-2\pi),\\
 \mathcal F_h[\psi(2\pi-\cdot)](\rho)
 &=e^{-2\pi i\rho/h}(\mathcal F_h\psi)(-\rho).
 \label{eq:fourier-intertwining-identities}
\end{align}
For analytic families of solutions with locally uniform bounds
\eqref{eq:fourier-transform-growth-class}, the transform is analytic
in the parameters.
\end{lemma}

\begin{proof}
\textbf{The two half-line integrals.}
For $x\in\R$, put
\begin{align}
 U_x(\rho)&\coloneq\int_x^{x+i\infty}
                  e^{-i\rho z/h}\psi(z)\dd z,&
 D_x(\rho)&\coloneq\int_{x-i\infty}^x
                  e^{-i\rho z/h}\psi(z)\dd z.
\end{align}
The first integral converges for $\Re\rho$ sufficiently negative,
and the second for $\Re\rho$ sufficiently positive. Define the
entire function
\begin{align}
 R_x(\rho)&\coloneq
 e^{i\rho}\int_x^{x+h}e^{-i\rho z/h}\psi(z)\dd z
 -e^{-i\rho}\int_{x-h}^x e^{-i\rho z/h}\psi(z)\dd z.
\end{align}
Integrate \eqref{eq:main-cos} on either half-line and shift the
contours by $h$. On a sufficiently far left or right half-plane,
\eqref{eq:fourier-transform-growth-class} makes the horizontal
segments at imaginary infinity tend to zero. The finite segments
at $x$ give
\begin{align}
 \lambda\bigl(U_x(\rho+h)+U_x(\rho-h)\bigr)
 +(2\cos\rho-2E)U_x(\rho)&=R_x(\rho),\\
 \lambda\bigl(D_x(\rho+h)+D_x(\rho-h)\bigr)
 +(2\cos\rho-2E)D_x(\rho)&=-R_x(\rho).
 \label{eq:fourier-half-line-recurrences}
\end{align}
Since $\lambda\ne0$, the first recurrence extends $U_x$ to the
right, one strip of width $h$ at a time. The second extends $D_x$
to the left. Each extension agrees with the previous function on
an open overlap. Thus both extend uniquely to entire functions,
with no division by a function of $\rho$.

\textbf{Independence of the straight vertical line.}
For real $x_1,x_2$, contour deformation gives
\begin{align}
 U_{x_2}(\rho)-U_{x_1}(\rho)
 &=-\int_{x_1}^{x_2}e^{-i\rho z/h}\psi(z)\dd z,\\
 D_{x_2}(\rho)-D_{x_1}(\rho)
 &=\int_{x_1}^{x_2}e^{-i\rho z/h}\psi(z)\dd z.
\end{align}
These identities first hold on the respective half-planes of
convergence and then everywhere by analytic continuation.
Consequently
\begin{align}
 \mathcal F_h\psi&\coloneq U_x+D_x
\end{align}
is independent of $x$ and is linear in $\psi$. Adding
\eqref{eq:fourier-half-line-recurrences} gives
\eqref{eq:fourier-dual-equation}. It also gives
\eqref{eq:fourier-bilateral-integral} wherever that integral
converges.

\textbf{Changes of variable.}
Translation by $2\pi$ changes the straight vertical line and
multiplies both half-line integrals by $e^{-2\pi i\rho/h}$.
Multiplication by $e^{2\pi iz/h}$ replaces $\rho$ by $\rho-2\pi$.
The substitution $z\mapsto2\pi-z$ exchanges the upper and lower
half-lines and replaces $\rho$ by $-\rho$, with the same factor
$e^{-2\pi i\rho/h}$. These substitutions are valid on the initial
half-planes and hence for their entire continuations. They prove
\eqref{eq:fourier-intertwining-identities}. The transformed
functions on the left satisfy the original equation because
$Nh=2\pi$. Finally, locally uniform exponential bounds give
locally uniform convergence of the initial integrals. The
recurrences preserve analytic dependence on the parameters.
\end{proof}

\subsection{Identification of the transformed solutions}

We now identify the transform of the normalized plus solution.
The bounds at both imaginary infinities determine it up to a scalar.

\begin{lemma}
\label{lem:fourier-normalization}
Fix $\Im\lambda>0$ and a compact set
$\mathcal K\Subset\operatorname{int}\Sigma_\lambda$.
For the parameters in \eqref{eq:coefficient-dual-parameters}, set
\begin{align}
 \widehat\phi_+^0(\rho;E,\lambda)
 &\coloneq\overline{\phi_+^0(\overline\rho;
                                  E^\dagger,\lambda^\dagger)},&
 \widehat\phi_-^0(\rho)&\coloneq\widehat\phi_+^0(2\pi-\rho).
 \label{eq:fourier-dual-normalized-solutions}
\end{align}
Here the suppressed parameter $h$ is unchanged, and the solution
under the conjugation is normalized by $\mathcal A=-2i$ for the
coupling $\lambda^\dagger$. For all sufficiently large $N$, there
is a nonvanishing function $\kappa(E,h)$, analytic near $\mathcal K$,
such that
\begin{align}
 (\mathcal F_h\phi_+^0)(\rho)
 &=\kappa(E,h)\widehat\phi_+^0(\rho),
 &\rho&\in\C.
 \label{eq:fourier-normalized-identification}
\end{align}
\end{lemma}

\begin{proof}
\textbf{A convergent integral.}
The strip estimates of \Cref{prop:entire-WKB-normalization} and
finite iteration of \eqref{eq:main-cos} give
\eqref{eq:fourier-transform-growth-class} for $\phi_+^0$ on every
bounded-width vertical strip. Thus \Cref{lem:fourier-transform}
applies. Write $\theta\coloneq\arg\lambda\in(0,\pi)$ and choose
$0<d<\min(\delta',\theta,\pi-\theta)/4$, where $K_{\delta'}$
is a fixed normalization strip valid near $\mathcal K$.
By \eqref{eq:momentum-at-infinity} and
\eqref{eq:rewrite-entire-strip-WKB}, for $t\ge0$,
\begin{align}
 \lvert\phi_+^0(\pi-d\pm it)\rvert
 &\le C_h\exp\left(
       \left(\frac{-d\mp\theta}{h}-\frac12\right)t\right).
\end{align}
Hence $F\coloneq\mathcal F_h\phi_+^0$ is given by the convergent
integral on $\Re z=\pi-d$ when $\lvert\Re\rho-\theta\rvert<d$.
On $\lvert s-\theta\rvert\le d/2$ this gives
\begin{align}
 \lvert F(s+it)\rvert&\le C_h e^{(\pi-d)t/h},&
 \lvert F(s-it)\rvert&\le C_h e^{-(\pi-d)t/h}.
 \label{eq:fourier-initial-growth}
\end{align}
Throughout this proof $h$ is fixed when $t\to\infty$. The
constants $C_h$ may depend on $h$, but are uniform for $E$ near
the fixed compact set.

\textbf{Identification by the growth at both ends.}
Choose an integer $m\ge0$ with $\pi-mh=\theta+O(h)$.
To reach an $h$-cell about $\pi$, iterate
\eqref{eq:fourier-dual-equation} $m+O(1)$ times from
\eqref{eq:fourier-initial-growth}. Each transfer matrix has norm
at most $C e^t$ on a bounded real interval. Thus, uniformly for
$\lvert x\rvert\le2h$,
\begin{align}
 \lvert F(\pi+x+it)\rvert
 &\le C_h e^{((2\pi-\theta-d)/h+C_0)t},\\
 \lvert F(\pi+x-it)\rvert
 &\le C_h e^{((-\theta+d)/h+C_0)t},
 \label{eq:fourier-propagated-growth}
\end{align}
where $C_0$ is independent of $h$. Since
$\arg\lambda^\dagger=\theta$, the dual strip estimates give
two-sided bounds for $\lvert\widehat\phi_+^0(\pi+x\pm it)\rvert$
by constant multiples of
\begin{align}
 \exp\left(\left(\frac{x\pm\theta}{h}-\frac12\right)t\right).
\end{align}
Both functions solve the dual equation. Their Wronskian is entire
and $h$-periodic, and the preceding bounds imply, for $\lvert x\rvert\le h/2$,
\begin{align}
 \lvert W[\widehat\phi_+^0,F](\pi+x+it)\rvert
 &\le C_h e^{((2\pi-d)/h+C_1)t},\\
 \lvert W[\widehat\phi_+^0,F](\pi+x-it)\rvert
 &\le C_h e^{((-2\theta+d)/h+C_1)t}.
\end{align}
For small $h$, the upper rate is strictly less than $2\pi/h$,
and the lower rate is negative. Integration over an $h$-period
at height $t$ therefore excludes every negative Fourier mode.
Integration at height $-t$ excludes every nonnegative mode.
Thus $W[\widehat\phi_+^0,F]=0$.

Conjugation changes the normalized Wronskian to
$W[\widehat\phi_-^0,\widehat\phi_+^0]=2i$. The Wronskian
expansion consequently gives $F=\alpha\widehat\phi_+^0$, where
$\alpha$ is entire and $h$-periodic. Dividing
\eqref{eq:fourier-propagated-growth} by the dual lower bounds gives
\begin{align}
 \lvert\alpha(\pi+x+it)\rvert
 &\le C_h e^{((2\pi-2\theta-d)/h+C_2)t},&
 \lvert\alpha(\pi+x-it)\rvert
 &\le C_h e^{(d/h+C_2)t}.
\end{align}
Both rates are strictly less than $2\pi/h$. The same Fourier
argument excludes all nonconstant modes, so $\alpha=\kappa(E,h)$.

\textbf{Nonvanishing and analytic dependence.}
Put $L\coloneq\pi-d$. For real $\eta$, the integral at
$\rho=\theta+i\eta$ is
\begin{align}
 F(\theta+i\eta)
 &=i e^{-i\theta L/h}e^{\eta L/h}
   \int_\R\phi_+^0(L+iy)e^{\theta y/h}e^{i\eta y/h}\dd y.
\end{align}
The function $y\mapsto\phi_+^0(L+iy)e^{\theta y/h}$ is nonzero
and belongs to $L^1(\R)$. Injectivity of the ordinary Fourier
transform therefore gives $F\not\equiv0$, so $\kappa\ne0$.
The construction is analytic in $E$. Finally,
$\widehat\phi_+^0(\pi)$ is nonzero by its uniform WKB estimate,
so $\kappa=F(\pi)/\widehat\phi_+^0(\pi)$ is analytic near
$\mathcal K$. All strip widths and strict growth inequalities
above can be fixed on that neighborhood.
\end{proof}

\section{Exterior equal-growth curves}
\label{app:exterior-cuts}

In this appendix, we prove \Cref{prop:general-exterior-cuts}. Our first goal is to show
that the four curves composing $\mathcal W_\lambda$ are disjoint
and extend from their respective corners to infinity. We work first
at $1+\lambda$ with $\Im\lambda>0$ and will use symmetry to get statements for the other three corners.

\subsection{The vanishing action at \texorpdfstring{$1+\lambda$}{1+lambda}}

Fix $\Im\lambda>0$ and use the action $V_+$ normalized in
\eqref{eq:upper-right-vanishing-action-normalization}.
Throughout this appendix, $\sqrt{\lambda}$ denotes the principal
square root. Recall from \eqref{eq:exterior-log-period} that
$c_\lambda=i\log\lambda$, with the logarithm fixed in
\eqref{eq:action-identity}. By \Cref{lem:vanishing-action-two-periods},
\begin{align}
 V_+&=P+Z+c_\lambda.
\end{align}

Write $E=1+\lambda+u+\lambda v$, with $u,v$ real.
We will relate $V_+$ to the complete elliptic integral of the first kind:
\begin{align}
 \mathrm K(k)&\coloneq
 \int_0^{\pi/2}\frac{\dd\vartheta}
 {\sqrt{1-k\sin^2\vartheta}},
 \qquad k\in\C\setminus[1,\infty).
\end{align}
Here the square root is the principal branch.

\begin{lemma}
\label{lem:vanishing-period-derivative}
With all branches continued from the corner,
\begin{align}
 V_+'(E)
 &=-\frac2{\pi\sqrt{\lambda}}\mathrm K(k(E)),&
 k(E)
 &=\frac{(1+\lambda)^2-E^2}{4\lambda}.
 \label{eq:vanishing-action-derivative}
\end{align}
\end{lemma}
\begin{proof}
We first consider $E$ near $E_+=1+\lambda$, with $E\ne E_+$.
By \Cref{lem:vanishing-action-two-periods},
\begin{align}
 V_+(E)&=\frac2\pi\int_0^{\zeta_+(E)}p(z)\dd z.
\end{align}
Using $\partial_Ep(z)=-(\sin p(z))^{-1}$, we obtain
\begin{align}
 V_+'(E)
 &=-\frac2\pi\int_0^{\zeta_+}
       \frac{\dd z}{\sqrt{1-(E-\lambda\cos z)^2}}. \label{eq:4405}
\end{align}
The endpoint term vanishes because $p(\zeta_+)=0$.
The square root is $\sin p(z)$ on the chosen branch, and its
reciprocal has an integrable singularity at
$\zeta_+$.
Setting $u=\cos z$ gives
\begin{align}
 V_+'(E)
 &=-\frac2\pi\int_{(E-1)/\lambda}^{1}
 \frac{\dd u}{\sqrt{(1-u^2)(1-(E-\lambda u)^2)}}.
\end{align}

Making the substitution
\begin{align}
 t^2&=\frac{2\lambda(1-u)}
 {(1+\lambda-E)(E+1-\lambda u)},
\end{align}
we get
\begin{align}
 V_+'(E)
 &=-\frac2{\pi\sqrt{\lambda}}\int_0^1
 \frac{\dd t}{\sqrt{(1-t^2)(1-k(E)t^2)}}=-\frac2{\pi\sqrt{\lambda}}\mathrm K(k(E)),
\end{align}
where $k(E)= (4\lambda)^{-1}((1+\lambda)^2-E^2)$.

Near the corner, the transformed contour can be deformed to $[0,1]$.
The square-root branch is fixed by $V_+'(1+\lambda)=-1/\sqrt{\lambda}$,
and the identity then extends by analytic continuation.

\end{proof}

\begin{prop}
\label{prop:general-outward-wall-branches}
Suppose $\Im\lambda>0$.
\begin{enumerate}
\item In the outward cone, there are real-analytic functions
$u,v\colon[0,\infty)\to[0,\infty)$ satisfying
\begin{align}
 u(0)&=v(0)=0,&
 u(r)+v(r)&=r,&
 u'(r),v'(r)&>0,
\end{align}
such that
\begin{align}
 \Gamma_\lambda
 &=\{1+\lambda+u(r)+\lambda v(r):r\geq0\}.
\end{align}
This is a simple curve from $1+\lambda$ to infinity.
It divides the outward cone $\R_{\geq1}+\lambda\R_{\geq1}$
into two connected regions. The region containing the boundary ray
$1+\lambda + \R_{>0}$ has $\Im V_+>0$, while the other
region has $\Im V_+<0$.

\item Outside the cone, continue $V_+$ along the straight segment
from $1+\lambda$ to $E=1+\lambda+u+\lambda v$. Then
\begin{align}
 \Im V_+(E)&<0 &&(u<0<v),\\
 \Im V_+(E)&>0 &&(v<0<u).
 \label{eq:general-wall-side-strip-signs}
\end{align}
\end{enumerate}
\end{prop}

\begin{proof}
\textbf{Derivative signs.}
Let $E=1+\lambda+u+\lambda v$, $u,v\geq0$ and define $H(u,v)\coloneq\Im V_+(1+\lambda+u+\lambda v)$.
Initially near the corner, \Cref{lem:vanishing-period-derivative} gives
\begin{equation} \label{eq:4476}
  \begin{split}
   H_u &= \Im V_+' (E) =\Im\left(-\frac2{\pi\sqrt{\lambda}}\mathrm K(k(E))\right), \\
H_v &= \Im (\lambda V_+'(E)) =\Im\left(-\frac{2\sqrt{\lambda}}{\pi}\mathrm K(k(E))\right).
\end{split}
\end{equation}

Observe that
\begin{align}
 1-k(E)
 &=\frac{(u+\lambda(v+2))(u+2+\lambda v)}{4\lambda},&
 \lvert\arg(1-k(E))\rvert&<\arg\lambda. \label{eq:4487}
\end{align}
Thus $k(E)\notin[1,\infty)$.
The function $-2\mathrm K(k(E))/(\pi\sqrt{\lambda})$ is therefore
holomorphic throughout the cone and agrees with $V_+'$ near the
corner by \Cref{lem:vanishing-period-derivative}. Since the cone is
convex, integrating this function from $1+\lambda$, with initial
value zero, gives a single-valued analytic continuation of $V_+$
throughout the cone. Thus \eqref{eq:4476} holds throughout the cone.

Observe that
\begin{align}
 \mathrm K(k(E))
 &=\int_0^{\pi/2}
 \frac{\dd\vartheta}
 {\sqrt{\cos^2\vartheta+(1-k(E))\sin^2\vartheta}},&
 \lvert\arg\mathrm K(k(E))\rvert&<\frac12\arg\lambda.
\end{align}
The argument of the expression under the square root lies between
$0$ and $\arg(1-k(E))$. By \eqref{eq:4487}, the integrand therefore has argument between
$0$ and $-\frac12\arg(1-k(E))$, as does its integral. Therefore,
by \eqref{eq:4476}, $H_u >0$ and $H_v<0$.
These inequalities hold also on the cone edges. At the corner,
$V_+'=-1/\sqrt{\lambda}$.

\textbf{The zero curve.}
For every $r>0$, the derivative signs and $H(0,0)=0$ give
\begin{align}
 H(0,r)&<0,&H(r,0)&>0,&
 \frac{\dd}{\dd u}H(u,r-u)&=H_u-H_v>0.
\end{align}
Thus each segment $u+v=r$, $u,v\geq0$, contains exactly one zero.
The implicit function theorem gives real-analytic functions
$u(r),v(r)$, including at $r=0$, and differentiation gives
\begin{align}
 u'(r)&=\frac{-H_v}{H_u-H_v}>0,&
 v'(r)&=\frac{H_u}{H_u-H_v}>0.
\end{align}
All derivatives on the right are evaluated at $(u(r),v(r))$.
Since $u(r)+v(r)=r$ and $1,\lambda$ are linearly independent over
$\R$, we have $\lvert E(r)\rvert\to\infty$.

On each segment, $H$ is negative before its zero and positive after
it. These portions form two connected regions, containing the edges
$u=0$ and $v=0$, respectively. This proves the separation and sign
statements.

\textbf{Signs in the adjacent regions.}
Put $d\coloneq u+\lambda v$ and $E(t)\coloneq1+\lambda+td$, $0\leq t\leq1$.
Using prime to denote differentiation with respect to $E$, we have
\begin{align}
  \frac{\dd}{\dd t}V_+(E(t))=dV_+'(E(t)).
\end{align}
Hence $V_+(E)=\int_0^1 dV_+'(E(t))\dd t$,
since $E(1)=E$ and $V_+(E(0))=V_+(1+\lambda)=0$.

In the integral for $V_+'$ in \eqref{eq:4405}, substituting
$\cos z=1+(td/\lambda)\sin^2\vartheta$ gives us
\begin{align}
 dV_+'(E(t))
 &=-\frac{2d}{\pi\sqrt{\lambda}}\int_0^{\pi/2}
 \frac{\dd\vartheta}{\sqrt{A_t(\vartheta)B_t(\vartheta)}},
 \label{eq:side-strip-radial-derivative}\\
 A_t(\vartheta)&=1+\frac t2\frac d{\lambda}\sin^2\vartheta,&
 B_t(\vartheta)&=1+\frac t2d\cos^2\vartheta.
\end{align}
Suppose $u<0<v$. Both $d$ and $d/\lambda$ lie in the
upper half-plane, and
\begin{align}
 \arg(d/\sqrt{\lambda})
 &=\frac12\bigl(\arg d+\arg(d/\lambda)\bigr),\\
 0\leq\arg A_t&<\arg(d/\lambda),&
 0\leq\arg B_t&<\arg d.
\end{align}
Taking each square root by continuation from one at $t=0$ gives
\begin{align}
 0
 &<\arg\left(
 \frac{d}{\sqrt{\lambda}\sqrt{A_t}\sqrt{B_t}}
 \right)  =\frac12\bigl(
 \arg(d/\lambda)-\arg A_t+\arg d-\arg B_t
 \bigr)
 <\pi.
\end{align}
Thus \eqref{eq:side-strip-radial-derivative} gives
$\Im(dV_+'(E(t)))<0$ so that
\begin{align}
 \Im V_+(E)
 &=\int_0^1\Im\bigl(dV_+'(E(t))\bigr)\dd t<0.
\end{align}

For $u>0>v$, both $d$ and $d/\lambda$ lie in the
lower half-plane. The same argument places the integrand's
argument in $(-\pi,0)$, giving $\Im V_+(E)>0$.
In both cases, $A_t$ and $B_t$ never vanish, so the
square roots remain well defined along the continuation.
\end{proof}

We will show that the four cuts divide the exterior of $\Sigma_\lambda$
into four components, as illustrated in
\Cref{fig:four-exterior-regions}. For $\Im\lambda>0$, we denote these
components by $\mathsf N,\mathsf E,\mathsf S,\mathsf W$, according to
their adjacent sides $\lambda+(-1,1)$, $1+\lambda(-1,1)$,
$-\lambda+(-1,1)$, and $-1+\lambda(-1,1)$, respectively.
The following corollary establishes this description and the action
inequalities in each region.

\fourexteriorregionsfigure

\begin{corr}
\label{cor:four-corner-curves-signs}
For every nonreal $\lambda$, the four curves
$\Gamma_\lambda,-\Gamma_\lambda,\Gamma_{-\lambda},-\Gamma_{-\lambda}$
are simple, disjoint, and extend from their respective corners to
infinity within their outward cones. They divide the exterior of
$\Sigma_\lambda$ into four components, one adjacent to each open side.

Suppose $\Im\lambda>0$. With the branches of $P,Z$ continued together
from the interior without crossing the cuts, we have
\begin{equation}
\begin{array}{c|cc}
 \text{Region}&\multicolumn{2}{c}{\text{Action inequalities}}\\ \hline
 \mathsf N&\Im(P+Z+c_\lambda)<0&\Im(P-Z+c_\lambda)<0\\
 \mathsf E&\Im(P+Z+c_\lambda)>0&\Im(-P+Z+c_\lambda)>0\\
 \mathsf S&\Im(-P+Z+c_\lambda)<0&\Im(-P-Z+c_\lambda)<0\\
 \mathsf W&\Im(P-Z+c_\lambda)>0&\Im(-P-Z+c_\lambda)>0.
\end{array}
\label{eq:exterior-region-action-signs}
\end{equation}
Complex conjugation gives the corresponding statements for
$\Im\lambda<0$.
\end{corr}

\begin{proof}
By conjugation,
\begin{align}
 V_+(E;\lambda)&=\overline{V_+(\overline E;\overline\lambda)}.
\end{align}
Thus \Cref{prop:general-outward-wall-branches} also describes the
curve for $\Im \lambda < 0$, with the
imaginary-part signs reversed. Applying this with couplings $\lambda$
and $-\lambda$, and then replacing $E$ by $-E$, gives the four
curves in the statement. They lie in disjoint outward cones.

For $\Im\lambda>0$, write $E=s+\lambda t$, with $s,t\in\R$.
The four \emph{side strips} are the half-infinite strips
\begin{equation}
\begin{array}{c|cc}
 \text{Strip}&s&t\\ \hline
 \text{Upper}&-1<s<1&t>1\\
 \text{Right}&s>1&-1<t<1\\
 \text{Lower}&-1<s<1&t<-1\\
 \text{Left}&s<-1&-1<t<1.
\end{array}
\end{equation}
Each cut divides its corner cone into two regions. Each exterior
component consists of one side strip and the adjoining portions of
the two neighboring corner cones, including their common edges.
In particular, the upper strip is only part of $\mathsf N$:
$\mathsf N$ also extends into the upper-left and upper-right cones,
up to their cuts.

By
\eqref{eq:upper-right-action-difference} and the corner symmetries
\eqref{eq:corner-action-definition}, we have, modulo $2\pi\Z$,
\begin{equation}
\label{eq:corner-action-difference}
\begin{array}{r@{}l@{\qquad}r@{}l}
 P+Z+c_\lambda&{}=V_+(E;\lambda),&
 P-Z+c_\lambda&{}=-V_+(-E;-\lambda),\\
 -P+Z+c_\lambda&{}=-V_+(E;-\lambda),&
 -P-Z+c_\lambda&{}=V_+(-E;\lambda).
\end{array}
\end{equation}

In the upper side strip, we have
\begin{align}
 E&=1+\lambda+(s-1)+\lambda(t-1),\\
 -E&=1-\lambda+(-s-1)+(-\lambda)(t-1).
\end{align}
Since $-1<s<1$ and $t>1$, both decompositions satisfy $u<0<v$.
Hence \eqref{eq:general-wall-side-strip-signs}, with conjugation
for the coupling $-\lambda$, gives
\begin{align}
 \Im V_+(E;\lambda)&<0,&
 \Im V_+(-E;-\lambda)&>0.
\end{align}
The first two identities in \eqref{eq:corner-action-difference}
therefore give the two inequalities for $\mathsf N$.

To cover the whole of $\mathsf N$, consider its portions in the
neighboring corner cones. In the upper-right cone, the sign of
$V_+(E;\lambda)$ follows from the position relative to
$\Gamma_\lambda$: the portion attached to the upper side strip has
negative imaginary part. The second sign still follows from
$-s-1<0<t-1$. In the upper-left cone, the roles reverse: the sign
of $-V_+(-E;-\lambda)$ follows from its own cut, while
$s-1<0<t-1$ gives the first sign. On the cone edges, the signs
follow from the strict derivative inequalities in
\Cref{prop:general-outward-wall-branches}.

For the other regions, the same substitutions are summarized below.
The entries are the displacements $u,v$ from the upper-right corner
of the problem appearing in each action:
\begin{equation}
\begin{array}{c|cc}
 \text{Action}&u&v\\ \hline
 V_+(E;\lambda)&s-1&t-1\\
 -V_+(-E;-\lambda)&-s-1&t-1\\
 -V_+(E;-\lambda)&s-1&-t-1\\
 V_+(-E;\lambda)&-s-1&-t-1.
\end{array}
\end{equation}
For each action in this table, $u<0<v$ gives negative imaginary
part and $v<0<u$ gives positive imaginary part. In the two rows
with coupling $-\lambda$, this follows by conjugation and the
displayed minus sign.

In the lower side strip, $-1<s<1$ and $t<-1$, so the last two
rows give the two negative signs for $\mathsf S$.
In the right side strip, $s>1$ and $-1<t<1$, the first and third
rows give the two positive signs for $\mathsf E$.
In the left side strip, $s<-1$ and $-1<t<1$, the second and fourth
rows give the two positive signs for $\mathsf W$.
In each neighboring corner cone, the action belonging to that
corner has the required sign on the side of its cut attached to
the strip. For the action belonging to the other endpoint of the
side, the table still gives displacements of opposite signs.
This extends the remaining six inequalities throughout their respective
regions. The latter inequalities remain strict on the cuts as well.

All these branches agree with continuation from the side strips.
Indeed, the regions with $u<0<v$ or $v<0<u$ are simply connected,
and the continued integral \eqref{eq:side-strip-radial-derivative}
is regular there. The continuations therefore agree with the
branches fixed near the corresponding corners. This proves
\eqref{eq:exterior-region-action-signs}.
\end{proof}

\begin{proof}[Proof of \Cref{prop:general-exterior-cuts}]
\textbf{Continuation of the actions.}
By \Cref{cor:four-corner-curves-signs}, the four cuts are disjoint
curves from the corners to infinity. After adding infinity, they
form a tree, whose complement is simply connected. The action
contours can be continued along paths avoiding the corners, since
these are the only energies where turning points coalesce.
Thus the interior actions $P,Z$ have single-valued analytic
continuations to $\C\setminus\mathcal W_\lambda$.

\textbf{Identification of $Q$.}
Assume $\Im\lambda>0$. Complex conjugation treats $\Im \lambda < 0$.
Comparing the momentum branches in the integrals defining $Q$ and
$P$ gives $Q=P$ near the side adjacent to $\mathsf N$ and $Q=-P$
near the side adjacent to $\mathsf S$.

Near the sides adjacent to $\mathsf E$ and $\mathsf W$, apply the
same comparison to the dual equation. By
\Cref{lem:duality-swaps-actions}, the dual exterior action is $-Z$
near $\mathsf E$ and $Z$ near $\mathsf W$. The signs are fixed by
the choice $\Im p<0$ for the dual exterior momentum.
Apply \Cref{thm:exterior-quantization} to the original and dual
equations in \eqref{eq:mN-duality}. Since
$\lambda^N=e^{-iNc_\lambda}$, comparison gives
\begin{align}
 e^{iN(Q+Z+c_\lambda)}&=1+O(N^{-1})
 &&\text{near }\mathsf E,\\
 e^{iN(Q-Z+c_\lambda)}&=1+O(N^{-1})
 &&\text{near }\mathsf W.
\end{align}
For fixed $E$, divide each relation at $N+1$ by the relation at $N$
and let $N\to\infty$. This gives $Q=-Z-c_\lambda$ near $\mathsf E$
and $Q=Z-c_\lambda$ near $\mathsf W$, modulo $2\pi\Z$.
Analytic continuation extends these identities
throughout the respective regions:
\begin{equation}
\begin{array}{c|cccc}
 \text{Region}&\mathsf N&\mathsf E&\mathsf S&\mathsf W\\ \hline
 Q&P&-Z-c_\lambda&-P&Z-c_\lambda.
\end{array}
\label{eq:exterior-chamber-gap-table}
\end{equation}
All action identities here are understood modulo $2\pi\Z$.

\textbf{Comparison of the four exponentials.}
The four terms are
\begin{align}
 e^{iNP},\qquad e^{-iNP},\qquad
 e^{iN(Z-c_\lambda)},\qquad e^{iN(-Z-c_\lambda)}.
\end{align}
Since $\lvert e^{iNA}\rvert=e^{-N\Im A}$, it remains to show that
each of the other three actions has larger imaginary part than $Q$.

In $\mathsf N$, we have $Q=P$. Subtracting $Q$ from the other
actions gives
\begin{align}
 -P-Q&=-2Q,\\
 (Z-c_\lambda)-Q&=-(P-Z+c_\lambda),\\
 (-Z-c_\lambda)-Q&=-(P+Z+c_\lambda).
\end{align}
The first difference has positive imaginary part because the
exterior normalization gives $\Im Q<0$. The other two have positive
imaginary parts by the $\mathsf N$ row of
\eqref{eq:exterior-region-action-signs}.
In $\mathsf S$, where $Q=-P$, the differences are
$-2Q,-(-P+Z+c_\lambda),-(-P-Z+c_\lambda)$.
The $\mathsf S$ row gives the same conclusion.

In $\mathsf E$, we have $Q=-Z-c_\lambda$, so
\begin{align}
 (Z-c_\lambda)-Q&=2Z,\\
 P-Q&=P+Z+c_\lambda,\\
 -P-Q&=-P+Z+c_\lambda.
\end{align}
The dual exterior action is $-Z$ modulo a real period and has
negative imaginary part. Hence $\Im(2Z)>0$. The other two
differences have positive imaginary parts by the $\mathsf E$ row of
\eqref{eq:exterior-region-action-signs}.
In $\mathsf W$, the differences are
$-2Z,P-Z+c_\lambda,-P-Z+c_\lambda$.
Now the dual exterior action is $Z$, and the $\mathsf W$ row
gives the two remaining signs.

Thus $Q$ is the unique action with smallest imaginary part off the
cuts. On each cut the two actions associated with that corner have
equal imaginary parts, while the adjacent-corner signs keep the
other two strictly above them. This completes the proof.
\end{proof}

For purely imaginary $\lambda$, the curves have additional symmetry.

\begin{corr}
\label{cor:imaginary-exterior-walls}
For $\lambda=i\rho$, $\rho>0$, the upper-right curve extends from
$1+i\rho$ to infinity with strictly increasing real and imaginary
parts. Its tangent at the corner is parallel to $1+i$.
The other curves are its reflections in the coordinate axes.
If $\rho=1$, then
\begin{align}
 \Gamma_i&=\{t(1+i):t\geq1\}.
\end{align}
Complex conjugation gives the curves for $\lambda=-i\rho$.
\end{corr}
\begin{proof}
The monotonicity follows from \Cref{prop:general-outward-wall-branches}.
At the corner, $V_+'=-1/\sqrt{\lambda}$ gives $H_u/H_v=-1/\rho$.
Writing $x=1+u$, $y=\rho(1+v)$ therefore gives
$\dd y/\dd x=-\rho H_u/H_v=1$.
If $\rho=1$, then $k(t(1+i))=(1-t^2)/2\leq0$.
By \eqref{eq:vanishing-action-derivative},
$(1+i)V_+'(t(1+i))$ is real, so $V_+(t(1+i))$ is real for $t\geq1$.
Uniqueness of the zero on each segment $u+v=r$ gives the asserted ray.
\end{proof}

\subsection{Why the cuts matter}
\label{sec:exterior-obstruction}

Away from the cuts, one of the four exponentials is strictly larger
in modulus than the others. On a cut, two have the same modulus.
The following example shows why the cuts must be excluded when
extending the quantization condition of \Cref{thm:main}\textup{(i)}.

\begin{prop}
\label{prop:exterior-obstruction}
Let $\lambda=i$, $E=t(1+i)$, $t>1$, and $E_*=1/5$.  Continue the
interior actions along either path
\begin{align}
 E_*\longrightarrow3/2\longrightarrow E,
 \qquad E_*\longrightarrow E_*+3i/2\longrightarrow E.
 \label{eq:exterior-continuation-paths}
\end{align}
The first path approaches the diagonal cut through $\mathsf E$, and
the second through $\mathsf N$; see
\Cref{fig:exterior-continuation-paths}.
For either continuation, any corrections
$P_N-P,Z_N-Z=O(N^{-2})$ satisfy
\begin{align}
 \cos(NP_N(E))+i^N\cos(NZ_N(E))
 &\ne\tfrac12m_N(E,i)
\end{align}
for all sufficiently large $N$. In particular, if $E$ is an
eigenvalue of $H^N_{p_0,x_0,i}$, then
\begin{align}
 \cos(NP_N(E))+i^N\cos(NZ_N(E))
 &\ne R_N(p_0,x_0,i).
\end{align}
Thus the quantization condition of \Cref{thm:main}\textup{(i)}
does not extend to this cut using either continued pair of actions
with $O(N^{-2})$ corrections.
\end{prop}

\exteriorcontinuationpathsfigure

\begin{proof}
The vanishing action $V=2\widetilde\chi_c=P+Z-\pi/2$, with
$V(1+i)=0$, is analytic through the corner and in the first quadrant
by the derivative formula \eqref{eq:vanishing-action-derivative}.
By \Cref{cor:imaginary-exterior-walls}, $t(1+i)\in\Gamma_i$ for
$t>1$, so $V(t(1+i))$ is real. Both paths
give the same $V$: their difference encloses only a point through
which $V$ is analytic.

By \eqref{eq:exterior-chamber-gap-table}, with $c_i=-\pi/2$,
the first continuation gives $Z=\pi/2-Q$, while the second
gives $P=Q$. These identities persist on approaching the cut.
Together with $P+Z-\pi/2=V$, they give
\begin{align}
 (P,Z)&=(Q+V,\pi/2-Q),&
 (P,Z)&=(Q,\pi/2-Q+V),
\end{align}
respectively, modulo $2\pi\Z$ in each component.
These real periods do not affect the cosine terms.
Since $i^N=e^{i\pi N/2}$, $V$ is real, and $\Im Q<0$, division of
$2\cos(NP)+2i^N\cos(NZ)$ by $e^{iNQ}$ gives respectively
\begin{align}
 1+e^{iNV}+O(e^{2N\Im Q}),\qquad
 1+e^{-iNV}+O(e^{2N\Im Q}).
\end{align}
In contrast, \Cref{thm:exterior-quantization} gives
\begin{align}
 \frac{m_N(E,i)}{e^{iNQ(E)}}&=1+O(N^{-1}).
\end{align}
Since $V$ is real, each of the preceding expressions differs
from this quotient by a term of modulus $1+O(N^{-1})$.
Corrections $P_N-P,Z_N-Z=O(N^{-2})$ change those expressions
only by $O(N^{-1})$, so they cannot remove the discrepancy.
The spectral conclusion follows from
\eqref{eq:finite-spectral-condition}, since
$\tfrac12m_N(E,i)=R_N(p_0,x_0,i)$ whenever
$E\in\Spec(H^N_{p_0,x_0,i})$.
\end{proof}

\section{Proof of the uniform corner comparison}
\label{app:uniform-corner-comparison}

Throughout, $\Im\lambda>0$, $\theta=\arg\lambda\in(0,\pi)$,
$E_+=1+\lambda$, and $\alpha^2=\lambda$ with $0<\arg\alpha<\pi/2$.
We retain the normalization of Section~\ref{sec:exact-reduction}:
$\mathcal{A}=-2i$ and $b=\mathcal{B}/\mathcal{A}$.  All energy neighborhoods are
fixed independently of $h$; we use smaller nested neighborhoods to keep
the contours regular and to apply Cauchy estimates a fixed positive
distance from the boundary.

\subsection{The stable domain and the upper comparison}
\label{app:corner-stable-domain}

\begin{proof}[Proof of \Cref{prop:stable-corner-comparison}]
\textbf{A canonical curve at the corner.}
Write $\lambda=R e^{i\theta}$.  At $E=E_+$ the upper real-momentum arc is
\begin{align}
 z(u)&=\arccos\left(1+\frac{1-u}{\lambda}\right),
 &-1<u<1.
\end{align}
The branch in $0<\Re z<\pi$, $\Im z>0$ joins 0 to the point $\zeta_-$
where $p=\pi$.  For $r=1-u>0$,
\begin{align}
 (\lambda\sin z(u))^2&=-2r\lambda-r^2.
\end{align}
Its root with positive imaginary part has negative real part, so
$z'(u)=1/(\lambda\sin z(u))$ has both components negative.  Thus the
arc is a regular vertical graph, with both coordinates increasing
with $p=\arccos u$.

Take the region between the arc and the normalization strip near
$\Re z=\pi$.  Complete its left boundary by a horizontal cut from
$\zeta_-$ to the upward imaginary ray and by the downward imaginary
ray from zero.  The strip branch has $\Im p<0$ inside and, as in
\eqref{eq:general-boundary-sheet},
\begin{align}
 \Re p(iy)&<0 &&(y<0),&
 \Re p(iy)&>\pi &&(y>\Im\zeta_-).
\end{align}
The lower sign also follows from $p(z)=-i\alpha z+O(z^3)$; the upper
one follows by continuing above $\zeta_-$.

As in Appendix~\ref{sec:uniform-geometry-interior}, at a
regular point $p(z_\vartheta)=\vartheta\in(0,\pi)$, define
$H_\vartheta(z)$ by~\eqref{eq:separated-corridor-level-function}.
Thus, $\dd H_\vartheta/\dd y=p-\vartheta$ on the real-momentum arc,
and the minimum is zero attained at $z_\vartheta$.  The horizontal-cut
and ray inequalities make $H_\vartheta$ increase along the rest of the
boundary towards infinity.  Since $H_\vartheta(0)>0$, a fixed negative
level avoids the double point.  Its regular component is a vertical
graph satisfying~\eqref{eq:separated-level-canonical-derivative}.

Choose $0<\varepsilon<\min(\theta,\pi-\theta)$.  Offset the compact
real-momentum arc $p^{-1}([\varepsilon,\pi-\varepsilon])$ slightly to
the right, attaching a negative $H_\varepsilon$ level below and a
negative $H_{\pi-\varepsilon}$ level above. The curves and their
tangents approach the arc as the offset tends to zero, so the joins
can be rounded while preserving a fixed strict canonical margin:
$\Re p+\Im p\,x'$ is close to $p$ on the middle part and equals
$\varepsilon,\pi-\varepsilon$ on the lower and upper attached levels,
respectively.

The infinite ends require a check.  On the continued branch,
\begin{align}
 p(x+iy)&=\pi-x+\theta-i(y+\log R)+O(e^{-y}),&&y\longrightarrow+\infty,\\
 p(x+iy)&=x-\pi+\theta+i(y-\log R)+O(e^y),&&y\longrightarrow-\infty.
\end{align}
An attached level that does not first reach the original strip has
the asymptotes
\begin{align}
 x(y)&=\theta+\varepsilon+O(1/y),&&y\longrightarrow+\infty,\\
 x(y)&=\pi-\theta+\varepsilon+O(1/\lvert y\rvert ),&&y\longrightarrow-\infty.
\end{align}
The upper asymptote lies in $\theta<x<\pi+\theta$, and the lower one
in $\pi-\theta<x<2\pi-\theta$. The displayed expansions therefore
give $0<\Re p<\pi$ with a strict margin between each tail and the strip.
The curves in \eqref{eq:tail-curves-into-strip} bring both ends into the
normalization strip while preserving this margin.

\textbf{A domain carrying the original solution.}
Fill the region between this curve and the normalization strip.
Write the rounded left boundary as $x=x_b(y)$ and set
\begin{align}
 v(y)&\coloneq\Re p(x_b(y)+iy)
       +\Im p(x_b(y)+iy)x_b'(y).
\end{align}
For some $\varepsilon_0>0$, one has
$\varepsilon_0\leq v(y)\leq\pi-\varepsilon_0$. On the compact part, use
\begin{equation}
  \label{eq:2}  x'=\frac{\Re p-v(y)}{-\Im p}
\end{equation}
away from the strip. Join these curves to straight vertical lines
inside the strip by interpolating the right-hand sides only where both
differential equations give canonical curves. On the infinite ends use
\eqref{eq:tail-curves-into-strip}, with the transition bands and final
cutoff described in Appendix~\ref{app:geometry}. Choose the outer
boundary to satisfy the differential equation~\eqref{eq:2}. Uniqueness
keeps the curves ordered, and the tail construction gives straight
vertical ends inside the strip. These curves fill a horizontally
connected canonical domain.

Choose two such domains so that the closure of the inner domain stays
at a fixed positive distance from the boundary of the outer domain
(see Fig.~\ref{fig:stable-corner-domain}).  The inner domain contains
the entire normalization strip $K_{\delta'}$.  The unbounded WKB theorem of
Appendix~\ref{app:growth} gives the estimate on the inner domain;
after scalar normalization to $\mathcal{A}=-2i$,
\Cref{cor:entire-normalization-unique} identifies its solution with
$\phi_+^0$.

To reach an arbitrarily small fixed positive real point, write
$\alpha_1=\Re\alpha$.  Then
$z_\varepsilon=i\varepsilon/\alpha+O(\varepsilon^3)$ and, for
real $x=O(\varepsilon)$,
\begin{align}
 H_\varepsilon(x)
 &=-\frac{\alpha_1x^2}{2}
   +\frac{\alpha_1\varepsilon^2}{2\lvert \alpha\rvert ^2}
   +O(\varepsilon^4).
\end{align}
The level $-k\varepsilon^2$, with fixed small $k>0$, crosses the
positive real axis at
\begin{align}
 x_\varepsilon
 &=\varepsilon\sqrt{\lvert \alpha\rvert ^{-2}+2k/\alpha_1}+O(\varepsilon^3).
\end{align}
Given a small fixed $r>0$, choose $\varepsilon$ so that
$x_\varepsilon<r/2$.  The domain contains
a fixed neighborhood of $[r,\pi]$ whose closure lies inside the domain,
with phase the direct real
integral from $\pi$; no estimate at zero is used.

\textbf{The upper overlap.}
To cross a regular piece of the arc, write it as $a(y)+iy$, with
momentum $q(y)\in(0,\pi)$.
For a perturbed curve $x=a(y)+d(y)$, Taylor expansion and
$\dd\Im p(a(y)+iy)/\dd y=0$ give
\begin{align}
 \Re p+\Im p\,x'
 &=q(y)+\Im p'(a(y)+iy)d(y)d'(y)
       +O\bigl(d(y)^2(1+\lvert d'(y)\rvert )\bigr).
\end{align}
A small smooth bump therefore crosses the arc and returns with a
strict margin.  Blend a bounded field on this crossing into the
previous field where both are canonical.  This supplies a fixed
starting region to the left of the arc, where $\Im p>0$.

Choose a fixed height $0<Y<\Im\zeta_-$, and put $w_0\coloneq-\delta+iY$ for
a fixed small $\delta>0$.  The preceding crossing and
\Cref{thm:finite-segment-continuation} give the WKB expansion for
$\phi_+^0(w)$ near $w_0$.  To the left of the arc the continuation is
in the leftward growing direction.  From the original strip continue
$\phi_+^0$ rightward to a neighborhood of $2\pi+w_0$.  At the corner
$E_+=1+\lambda$, the segment $[\pi+iY,2\pi-\delta+iY]$ contains no
real-momentum point, and its branch has $\Im p<0$.  Both compact
continuation rectangles are regular and have strict growing-direction
margins.

On a fixed neighborhood in the overlap, let $p^{(0)},p^{(2\pi)}$ and
$S^{(0)},S^{(2\pi)}$ be the momentum branches and phase integrals
continued along these paths from $\pi$ and $3\pi$, respectively.  The
branches satisfy $p^{(0)}=-p^{(2\pi)}$.  For nearby interior energies,
the comparison loop goes clockwise around $2\pi-\zeta_+$.  The
continued amplitude $(\sin p)^{-1/2}$ is locally a nonzero analytic
factor times $(z-(2\pi-\zeta_+))^{-1/4}$, so this loop contributes a
factor $+i$.  Consequently the product of the two continued amplitudes
is $i/\sin p^{(2\pi)}$.  The same comparison of paths gives
\begin{align}
 S^{(0)}(z)+S^{(2\pi)}(z)
 &=2\int_\pi^{2\pi-\zeta_+}p(w)\dd w-2\pi P
 =-2\pi\widetilde\chi_c=-\pi V_+.
\end{align}
The amplitude and phase identities persist throughout the energy
neighborhood by analytic continuation of the fixed regular paths.
Expanding the Wronskian and using the one-step WKB expansion gives
\begin{align}
 \mathcal{C}(z)
 &=\phi_+^0(z+h)\phi_+^{2\pi}(z)
   -\phi_+^0(z)\phi_+^{2\pi}(z+h)\\
 &=\frac{i e^{-i\pi V_+/h}}{\sin p^{(2\pi)}(z)}
   \left[e^{-ip^{(2\pi)}(z)}-e^{ip^{(2\pi)}(z)}+O(h)\right]\\
 &=2e^{-i\pi V_+/h}(1+O(h)).
\end{align}
The relative error is uniform because the fixed overlap stays away
from the turning points, so $\sin p^{(2\pi)}$ is bounded away from zero.

Average $\mathcal{C}(z)=c+a e^{2\pi iz/h}$ over a horizontal $h$-segment
in the overlap.  Its nonconstant term has zero average, and hence
\begin{align}
 c&=2e^{-i\pi V_+/h}(1+O(h)),&
 V_h&=V_++\frac{ih}\pi\log\left(\frac c2 e^{i\pi V_+/h}\right).
 \label{eq:corner-upper-coefficient-proof}
\end{align}
The logarithm is the standard branch. Thus, $V_h=V_++O(h^2)$ and
$c=2e^{-i\pi V_h/h}$ exactly.

All domains and paths now stay fixed as $E$ varies.  Their compact
parts have positive turning-point distance and strict margins; on
their infinite ends,
$\partial_Ep=-1/\sin p=O(e^{-\lvert \Im z\rvert })$, so the same
margins persist.  The WKB construction is therefore uniform and
analytic in $E$, as in Appendix~\ref{app:growth}.

\textbf{A rough bound for $b$.}
At $E_+$ put
\begin{align}
 \kappa(x)&\coloneq\operatorname{arcosh}\bigl(1+\lambda(1-\cos x)\bigr),&
 \rho&\coloneq\Re\int_0^{2\pi}\kappa(x)\dd x>0.
\end{align}
The branch satisfying $\kappa(x)=\alpha x+O(x^3)$ has $\Re\kappa(x)>0$ for
$x\in(0,2\pi)$, $p=-i\kappa$, and $\kappa(2\pi-x)=\kappa(x)$.  Thus
$\Re\int_r^\pi\kappa\dd x=\rho/2-O(r^2)$.  For $r_h\in h\Z$ within $h$ of
the fixed $r$, WKB gives
\begin{align}
  \label{eq:3}
  \lvert \phi_+^0(r_h)\rvert +\lvert \phi_+^0(r_h+h)\rvert
 &\leq C_r\exp\left(-\frac{\rho/2-O(r^2)}h\right).
\end{align}
The phase from $\pi$ to $r$ is regular, so shrinking the energy
neighborhood changes the term $\rho/2-O(r^2)$ in~\eqref{eq:3} by an
arbitrarily small amount.

The one-step recurrence matrices are bounded by a common $K$ for small
$r,h$ and energies in a fixed compact set.  Propagation from $r_h+h$
to $-h$ costs at most $K^{r/h+O(1)}$, while
\begin{align}
 \mathcal{B}&=\phi_+^0(0)\bigl(\phi_+^0(-h)-\phi_+^0(h)\bigr).
\end{align}
Hence $\lvert b\rvert \leq C_\gamma e^{-\gamma/h}$ for any $\gamma<\rho$, by choosing $r$
small and then shrinking the energy neighborhood.  The point $r$ is
independent of $h$ sufficiently small.

\textbf{The exponent for large couplings.}  To make $2R_Nb$
exponentially small also for $\lvert \lambda\rvert >1$, put
$\Lambda\coloneq\max(1,\lvert \lambda\rvert )$. Since
$\lvert 2R_N\rvert \leq4\Lambda^N$ for real phases $(x_0,p_0)$, the
decay exponent must satisfy $\gamma>2\pi\log\Lambda$. Write
$\rho_\lambda=\rho$. We show
\begin{align}
 \rho_\lambda
 &=\rho_{1/\lambda}+2\pi\log\lvert \lambda\rvert 
 >2\pi\log\Lambda.
 \label{eq:corner-classical-exponent-duality}
\end{align}
Both integrals giving $\rho_\lambda$ and $\rho_{1/\lambda}$ are positive
for every nonreal coupling.  For the equality, use
\begin{align}
  \label{eq:4}
 \frac1{2\pi}\int_0^{2\pi}\log\lvert t-\cos y\rvert \dd y
 &=\Re\operatorname{arcosh}t-\log2,\qquad t\notin[-1,1].
\end{align}
This follows by writing $t=(q+q^{-1})/2$, $\lvert q\rvert >1$, and
averaging the logarithms of the three factors of $t-\cos y$ expressed
in terms of $q$.

Apply~\eqref{eq:4} in both
variables to $(1-\cos p)+\lambda(1-\cos x)$ to get
\begin{align}
  \label{eq:5}
  &\int_0^{2\pi}\left[\Re\operatorname{arcosh}(1+\lambda(1-\cos x))-\log2\right]\dd x\\
  &=\frac1{2\pi}\int_0^{2\pi}\int_0^{2\pi}\log|(1-\cos p)+\lambda(1-\cos x)|\dd p\dd x
  \\&=2\pi\log|\lambda|+\frac1{2\pi}\int_0^{2\pi}\int_0^{2\pi}\log|(1-\cos x)+\lambda^{-1}(1-\cos p)|\dd x\dd p
  \\&=2\pi\log|\lambda|+    \int_0^{2\pi}\left[\Re\operatorname{arcosh}(1+\lambda^{-1}(1-\cos p))-\log2\right]\dd p.
\end{align}
This proves~\eqref{eq:corner-classical-exponent-duality}.
\end{proof}

\subsection{The zeros of the reflected Wronskian}
\label{app:corner-B-zeros}

\begin{proof}[Proof of \Cref{lem:corner-B-zeros}]
\textbf{A zero forces parity.}
If $b=0$, then $a=0$, $c^2=\mathcal{A}^2$, and the period relation becomes
\begin{align}
 \phi_+^0(z-2\pi)&=-\frac c{\mathcal{A}} \phi_+^0(2\pi-z).
\end{align}
Hence $\phi_+^0$ is even or odd on $\C$, with a constant sign multiplier.

\textbf{Decay on the imaginary axis.}
Choose $\varepsilon>0$ so small that $L=\pi-\varepsilon$ lies in
the original normalization strip.  Its unbounded WKB estimate and
\eqref{eq:momentum-at-infinity} give, for $y>0$,
\begin{align}
 \lvert \phi_+^0(L+iy)\rvert &\leq C_h e^{-((\varepsilon+\theta)/h+1/2)y},\\
 \lvert \phi_+^0(L-iy)\rvert &\leq C_h e^{((\theta-\varepsilon)/h-1/2)y}.
\end{align}
Constants may depend on $h$; only decay for each small $h$ is needed.
Parity gives the left boundary estimate.  Put $d\coloneq\varepsilon/h+1/2$ and
\begin{align}
 g(z)&\coloneq\phi_+^0(z)\exp\left(-idz-\frac{i\theta z^2}{2Lh}\right).
\end{align}
It is bounded on the sides and bottom of the upper half-strip.
Finitely many recurrence steps from the normalization strip give
\begin{align}
 \log\lvert g(x+iy)\rvert &\leq C_h(1+y),&\lvert x\rvert &\leq L,\quad y\geq0.
\end{align}
For $0<k<\pi/(2L)$, the positive harmonic barrier
$\delta\cos(kx)\sinh(ky)$ dominates this linear bound on the top of
a sufficiently tall rectangle, is nonnegative on its sides, and
vanishes on its bottom.  The maximum principle, first sending the top
to infinity and then $\delta$ to zero, shows that $g$ is bounded.
Consequently
\begin{align}
 \lvert \phi_+^0(x+iy)\rvert &\leq C_h e^{-dy-\theta xy/(Lh)},&y&\geq0.
\end{align}
Parity gives the lower-half-strip estimate.  In particular, on
$\lvert x\rvert \leq h$,
\begin{align}
 \lvert \phi_+^0(x+iy)\rvert &\leq C_h e^{-\beta_h\lvert y\rvert },&
 \beta_h&=\frac\varepsilon h+\frac12-\frac\theta L>1
\end{align}
for small $h$. 

Set $u(t)\coloneq\phi_+^0(it)$.  These bounds put $u$, $\cosh(t)u$, and
$u(\cdot\pm ih)$ in $L^2(\R)$ and justify Fourier contour shifts.
With $D_t=-i\partial_t$,
\begin{align}
 \widehat{u(\cdot\pm ih)}(\xi)&=e^{\mp h\xi}\widehat u(\xi),&
 \bigl(\cosh(hD_t)+\lambda\cosh t\bigr)u&=Eu.
\end{align}
Taking the inner product with $u$ gives
\begin{align}
 E&=\frac{\int_\R\cosh(h\xi)\lvert \widehat u(\xi)\rvert ^2\dd\xi}
          {\lVert u\rVert ^2}
   +\lambda\frac{\int_\R\cosh t\,\lvert u(t)\rvert ^2\dd t}{\lVert u\rVert ^2}.
 \label{eq:corner-outward-zero-identity}
\end{align}
Both real quotients are strictly greater than one: a nonzero $L^2$
function cannot be supported at a single point.  Every zero of $b$
therefore lies in
\begin{align}
 E_++(0,\infty)+\lambda(0,\infty).
 \label{eq:corner-outward-zero-cone}
\end{align}

\textbf{Which outward points are zeros.}  Recall that $w_h=V_h/h$.
The identity~\eqref{eq:corner-division-free-trace} yields
\begin{align}
 b(m_N-2R_N)
 &=1+e^{-2\pi iw_h}
   +\bigl(1+\lambda^{2N}e^{2\pi iw_h}\bigr)b^2-2R_Nb.
 \label{eq:corner-reference-scalar}
\end{align}
Use $p_0=x_0=0$, for which $\lvert 2R_N\rvert \leq4\Lambda^N$, with
$\Lambda=\max(1,\lvert \lambda\rvert )$.  Choose $\gamma>2\pi\log\Lambda$
and shrink the neighborhood so that
$2\pi\lvert \Im V_h\rvert
<\gamma-2\pi\log\Lambda$. By~\eqref{eq:corner-corrected-action}, the last
two terms in~\eqref{eq:corner-reference-scalar} are then uniformly
$O(e^{-\eta/h})$.  By Cauchy's estimate, $V_h'=V_+'+O(h^2)$ on a
smaller disk.  It is uniformly close to $V_+'(E_+)=-1/\alpha$ there;
integrating along segments proves injectivity.  Thus $V_h$ is a
coordinate on a fixed disk.

Around each half-integer level $V_h=h(k+1/2)$ in a smaller
neighborhood, take a disk of radius $Mhe^{-\eta/h}$ in this
coordinate.  On its boundary
$\lvert 1+e^{-2\pi iV_h/h}\rvert \asymp M e^{-\eta/h}$, uniformly in $k$.
For fixed large $M$, Rouch{\'e}'s theorem gives one simple zero of
\eqref{eq:corner-reference-scalar} in each disk.

The ray $\alpha(0,\infty)$ lies strictly in
$(0,\infty)+\lambda(0,\infty)$ because $\arg\alpha=\theta/2$.
Since $E=E_+-\alpha V_h+O(V_h^2+h^2)$, the centers with $V_h<0$ lie
outward, a distance comparable to $\lvert h(k+1/2)\rvert $ from the cone boundary.
Here $\lvert h(k+1/2)\rvert \geq h/2$ absorbs the $O(h^2)$ error, and reducing
the fixed neighborhood absorbs the quadratic error.  Their disks
are outside $\Sigma_\lambda$, whereas the reference spectrum is
inside: its operator is a sum of two self-adjoint contractions with
coefficients $1,\lambda$.  Thus each disk contains a simple zero of
$b$, not of $m_N-2R_N$.  The exact identity then forces
$1+e^{-2\pi iw_h}=0$, so that zero is exactly the center.

Conversely, every zero of $b$ is one of these centers. Those with $V_h>0$
lie inward and are excluded by \eqref{eq:corner-outward-zero-cone}.
This proves the claimed zero set.  The nested neighborhoods ensure
that all disks stay within the region of uniform estimates.

At centers with $V_h>0$ the count instead gives a simple zero of
$m_N-2R_N$.  Every such zero satisfies
$1+e^{-2\pi iw_h}=O(e^{-\eta/h})$, hence belongs to one of these disks.
For any real phases the same spectral inclusion and bound on $2R_N$
hold.  This proves \eqref{eq:corner-uniform-Bohr-Sommerfeld}, uniformly
also for indices of order $h^{-1}$.
\end{proof}

\begin{proof}[Proof of \Cref{cor:real-phase-corner-quantization}]
By \Cref{lem:trace-polynomial}, the eigenvalues of
$H^N_{p_0,x_0,\lambda}$ are precisely the solutions of
\begin{align}
 m_N(E,\lambda)-2R_N(p_0,x_0,\lambda)&=0.
\end{align}
Put $h=2\pi/N$. Since $p_0,x_0\in\R$,
\begin{align}
 \lvert R_N\rvert
 &=\lvert\cos(Np_0)+\lambda^N\cos(Nx_0)\rvert
 \le 1+\lvert\lambda\rvert^N.
\end{align}
Thus the bound on $b$ in \Cref{prop:stable-corner-comparison} and
\eqref{eq:corner-reference-scalar} give, on a sufficiently small
fixed neighborhood of $E_+$,
\begin{align}
 b(E)(m_N-2R_N)
 &=1+e^{-2\pi iV_h(E)/h}+O(e^{-\eta/h}),
\end{align}
where the remainder is analytic in $E$ and uniform in the real phases.

It remains to locate the zeros of the right-hand side that are not
zeros of $b$. Since $V_h=V_++O(h^2)$ and
$V_+'(E_+)=-1/\alpha\ne0$, we may choose $\varepsilon>0$ so that,
for all sufficiently large $N$, each equation
\begin{align}
 V_h(E)&=h(n+\tfrac12),&0\le n&\le\varepsilon N,
\end{align}
has a unique solution in a fixed smaller neighborhood of $E_+$.
By \Cref{lem:corner-B-zeros}, $b$ does not vanish near these solutions.
The function $1+e^{-2\pi iV_h/h}$ has a simple zero at each of them.
The preceding Rouch{\'e} estimate therefore gives exactly one simple
eigenvalue $E_{n,N}$ near each solution, with
\begin{align}
 V_h(E_{n,N})&=h(n+\tfrac12)+O(he^{-\eta/h}),
\end{align}
uniformly in $n,p_0,x_0$.
Using \eqref{eq:corner-corrected-action} and
\eqref{eq:upper-right-action-difference}, we obtain
\begin{align}
 P(E_{n,N})+Z(E_{n,N})+i\log\lambda
 =V_h(E_{n,N})+O(h^2)=h(n+\tfrac12)+O(h^2).
\end{align}
For fixed $n$, we have $E_{n,N}-E_+=O(h)$, so
\eqref{eq:corner-vanishing-action-linearization} gives
\begin{align}
 E_{n,N}&=1+\lambda-\alpha h(n+\tfrac12)+O(h^2).
\end{align}
Substituting $h=2\pi/N$ and $\alpha=\sqrt{\lambda}$ proves the
corollary.
\end{proof}

\subsection{The regularized classical action}
\label{app:corner-regularized-action}

\begin{proof}[Proof of \Cref{lem:corner-regularized-action}]
To normalize $V_+$, put $\delta\coloneq E-E_+$, $t\coloneq2\sin(z/2)$, and
$s\coloneq2\sin(p/2)$.  The local equation is
\begin{align}
 s^2+\lambda t^2&=-2\delta.
\end{align}
On the outgoing part of the short cycle substitute
$t=\sqrt{-2\delta/\lambda}\,\xi$, $-1\leq\xi\leq1$, and
$s=\sqrt{-2\delta}\sqrt{1-\xi^2}$.  Expand
$p=2\arcsin(s/2)$ and $\dd z=\dd t/\sqrt{1-t^2/4}$.
The return piece gives the same integral.  The integrand is $\delta$
times an analytic function of $\delta$, uniformly integrable in $\xi$;
thus
\begin{align}
 V_+(E)&=-\frac{2\delta}{\pi\alpha}
          \int_{-1}^1\sqrt{1-\xi^2}\dd\xi+O(\delta^2)
       =-\frac\delta\alpha+O(\delta^2).
\end{align}
This fixes the analytic extension and orientation.  Choose
$\log V_+$ real on the inward curve $V_+>0$.
Below the negative $V_+$-axis the ordinary contour deforms to the
exterior real-period contour.  Above it, the deformation adds the
short cycle once.  Hence
\begin{align}
 P&=Q,&\log V_+&=\log(-V_+)-i\pi&&\text{below the cut},\\
 P&=Q+V_+,&\log V_+&=\log(-V_+)+i\pi&&\text{above the cut}.
 \label{eq:corner-two-contour-continuations}
\end{align}
The difference of the two deformations goes down one side of the short
cut between the points $\pm\zeta_+$ and back along the other: it is
the defining contour of $V_+$.

For the sign, the outward branch is
$p=-i\alpha\sqrt{z^2+2\delta/\lambda}$ to leading order.  Its two
endpoint integrals contribute the logarithmic term
$(i\delta/(2\pi\alpha))\log\delta$ to $P$, whose increment is
$-\delta/\alpha$, the leading term of $V_+$.  This fixes the
orientation; the exact contour difference is the full $V_+$.

Thus $J=-2\pi iP+V_+(\log V_+-1)$ is single-valued on the punctured
neighborhood.  It is bounded there: the ordinary contour is regular
outside a small fixed $z$-disk, and inside it has bounded-length
representatives avoiding the cut, with
$p=O(\lvert z\rvert +\lvert \delta\rvert ^{1/2})$.  The local integrals stay bounded, while
$V_+\log V_+\to0$.  The removable singularity theorem extends $J$ through
$E_+$.

Either row of \eqref{eq:corner-two-contour-continuations} gives
\eqref{eq:corner-exterior-regularization}, with $\log(-V_+)$ real on
the outward curve $V_+<0$.
\end{proof}

\subsection{Matching on a fixed contour around \texorpdfstring{$E_+$}{1+lambda} in the energy
  plane}
\label{app:corner-boundary-normalization}

For comparison on the energy contour, write
\begin{align}
 F_h&=b\,\Gamma(1/2+w_h),&
 F_{\mathrm{app}}&\coloneq\sqrt{2\pi}\,e^{J/h}h^{-w_h},&
 w_h&=\frac{V_h}h.
\end{align}
\begin{lemma}
\label{lem:corner-boundary-normalization}
There is a fixed simple closed contour around $E_+$, inside the
neighborhood of the preceding estimates, on which
\begin{align}
 \frac{F_h}{F_{\mathrm{app}}}&=1+O(h)
 \label{eq:corner-boundary-ratio-appendix}
\end{align}
uniformly as $h\longrightarrow0$.
\end{lemma}
\begin{proof}
Take the inverse image of a small fixed circle $\lvert V_+\rvert =r$ under
$V_+$. The images of the two adjacent sides are tangent at
zero to the rays $\arg V_+=\pm\theta/2$. For small $r$, the exterior
arc therefore stays a fixed angle from the positive $V_+$-axis.
The locally uniform estimates in \Cref{thm:interior-side-quantization}
cover neighborhoods of the two side crossings.

We use Stirling's formula \cite[Section~5.11]{NISTDLMF} only in closed
sectors away from its negative-real cut:
\begin{align}
 \Gamma(1/2+z)
 &=\sqrt{2\pi}\,e^{z(\log z-1)}(1+O(\lvert z\rvert ^{-1})).
\end{align}
Since $\lvert w_h\rvert \asymp h^{-1}$, the interior estimate
$b=e^{-2\pi iP/h}(1+O(h))$ and Stirling with $z=w_h$ give
\begin{align}
 F_h&=\sqrt{2\pi}\exp\left(-\frac{2\pi iP}{h}
                          +w_h(\log w_h-1)\right)(1+O(h)).
\end{align}
Since $V_h=hw_h=V_++O(h^2)$, the logarithm of the exponential ratio to
$F_{\mathrm{app}}$ is
\begin{align}
 \frac{V_h(\log V_h-1)-V_+(\log V_+-1)}h=O(h),
\end{align}
because $V_+$ stays away from zero.  This proves the interior estimate.

On the exterior arcs and the neighborhoods of the side crossings,
\Cref{prop:exterior-Wronskian-estimates} and
\Cref{thm:interior-side-quantization}
and their stated symmetry reductions give
\begin{align}
 m_N&=e^{2\pi iQ/h}(1+O(h)),&\lvert e^{2\pi iQ/h}\rvert &\geq e^{\eta_0/h}.
\end{align}
The exterior action $Q$ continues across each side. Near
$E=s+\lambda$, $-1<s<1$, it equals $P$; near $E=1+\lambda t$,
$-1<t<1$, it equals $-Z-c_\lambda$ modulo $2\pi\Z$, as in
\eqref{eq:exterior-chamber-gap-table}. The contour relations
\eqref{eq:corner-two-contour-continuations} fix the corresponding branches.

We must show that $d_hb/m_N$ is exponentially small, where
$d_h=1+\lambda^{2N}e^{2\pi iw_h}$.
Recall the strict gap
\begin{align}
 \rho&=\Re\int_0^{2\pi}\kappa\dd x>\ell,
 &\ell&=2\pi\log\max(1,\lvert \lambda\rvert ).
\end{align}
The normalized period tends to $-i(2\pi)^{-1}\int_0^{2\pi}\kappa$, so its exponential growth rate tends
to $\rho$.  Choose $\gamma<\rho$ close to $\rho$ and then shrink the
neighborhood until
$\lvert d_h\rvert \leq C\exp((2\ell+\epsilon)/h)$ and
$\lvert m_N\rvert \geq C^{-1}\exp((\rho-\epsilon)/h)$.
Together with $\lvert b\rvert \leq Ce^{-\gamma/h}$ from
\Cref{prop:stable-corner-comparison}, these bounds give
\begin{align}
 \left\lvert \frac{d_hb}{m_N}\right\rvert 
 &\leq C\exp\left(\frac{2\ell-\gamma-\rho+2\epsilon}{h}\right).
\end{align}
Since $\rho>\ell$, choosing $\gamma<\rho$ sufficiently close to $\rho$
and then $\epsilon>0$ sufficiently small makes the exponent strictly
negative.  This proves the required decay, also for large couplings.
The trace identity and Gamma reflection
\cite[Section~5.5]{NISTDLMF} give
\begin{align}
 F_h(m_N-d_hb)&=\frac{2\pi e^{-\pi iw_h}}{\Gamma(1/2-w_h)}.
\end{align}
Divide by the nonzero $m_N-d_hb$ before using Stirling with $z=-w_h$.
This keeps the calculation valid at zeros of $b$ and gives
\begin{align}
 F_h&=\sqrt{2\pi}\exp\left(-\frac{2\pi iQ}{h}-\pi iw_h
                         +w_h(\log(-w_h)-1)\right)(1+O(h)).
\end{align}
By the exterior expression for $J$, the logarithm of its exponential
ratio to $F_{\mathrm{app}}$ is
\begin{align}
 \frac{V_h(\log(-V_h)-1)-V_+(\log(-V_+)-1)
             -i\pi(V_h-V_+)}h=O(h).
\end{align}
The logarithm is regular through the negative $V_+$-axis.  This proves
the estimate on the exterior and open-side arcs; a finite cover makes
it uniform on the contour.  The main proof now applies the maximum
principle inside.
\end{proof}

\section{Action identities and interior Scottish-flag eigenvalues}
\label{app:action-identity}

\subsection{Action integral identities}

The following identity relates the position period to the half-period
integrals in \Cref{sec:classical-actions}.
We recall the integrals from \eqref{eq:action-definitions}, with the same branches and contours:
\begin{align}
   \widetilde\chi\coloneq\frac1\pi\int_{\zeta_+}^\pi p (z) \dd z, \qquad \chi_c\coloneq\frac1\pi\int_0^{\zeta_-} (\pi - p(z)) \dd z.
\end{align}

\begin{prop}
\label{prop:half-period-action-identity}
For $\Im \lambda >0 $ and $E \in \Sigma _\lambda$, we have
\begin{align}
 \widetilde\chi+\chi_c&=-i\log\lambda.
\end{align}
This identity remains valid under simultaneous analytic continuations of $\widetilde\chi$ and $\chi_c$.
\end{prop}
\begin{proof}
Let $I(E,\lambda)\coloneq\widetilde\chi+\chi_c$. We first show that
$I$ is independent of $E$, then evaluate it at $E=0$.

\textbf{Independence of $E$.}
Fix $\lambda$ and first take $E$ in the interior of $\Sigma_\lambda$,
so the turning points are simple.
Differentiating
$\cos p(z;E)+\lambda\cos z=E$ at fixed $z$ gives
\begin{align}
 -\sin p(z;E)\,\partial_Ep(z;E)&=1,
 &\partial_Ep(z;E)&=-\frac{1}{\sin p(z;E)}.
\end{align}
Thus
\begin{align}\label{eq:half-period-tilde-derivative}
 \partial_E\widetilde\chi
 &=-\frac1\pi p(\zeta_+;E)\zeta_+'
   +\frac1\pi\int_{\zeta_+}^{\pi}\partial_Ep(z;E)\dd z
   =-\frac1\pi\int_{\zeta_+}^{\pi}\frac{\dd z}{\sin p(z;E)},
\end{align}
since $p(\zeta_+;E)=0$. Similarly,
\begin{align}\label{eq:half-period-chi-c-derivative}
 \partial_E\chi_c
 &=\frac1\pi\bigl(\pi-p(\zeta_-;E)\bigr)\zeta_-'
   -\frac1\pi\int_0^{\zeta_-}\partial_Ep(z;E)\dd z
   =\frac1\pi\int_0^{\zeta_-}\frac{\dd z}{\sin p(z;E)},
\end{align}
since $p(\zeta_-;E)=\pi$.

Combining \eqref{eq:half-period-tilde-derivative} and
\eqref{eq:half-period-chi-c-derivative}, we obtain
\begin{align}
 \partial_E I
 &=-\frac1\pi\int_{\zeta_+}^{\pi}\frac{\dd z}{\sin p(z;E)}
   +\frac1\pi\int_0^{\zeta_-}\frac{\dd z}{\sin p(z;E)}.
\end{align}
Making the substitution $u=\cos z$, we get
\begin{align}\label{eq:sfsfsdfsew}
 \pi\lambda\partial_E I
 &=\int_{a_-}^{-1}\frac{\dd u}{w(u)}
   +\int_{a_+}^1\frac{\dd u}{w(u)},
\end{align}
where
\begin{align}
 a_\pm&\coloneq\frac{E\pm1}{\lambda},
 &w(u)&\coloneq\lambda^{-1}
 \sqrt{\bigl(1-(E-\lambda u)^2\bigr)(1-u^2)},
\end{align}
so that $\cos\zeta_+=a_-$ and $\cos\zeta_-=a_+$.
We choose the sign of the square root so that, along the integration paths,
\begin{align}
 \lambda w(\cos z)=\sin p(z;E)\sin z.
\end{align}

Choose the branch of $w$ with cuts along the straight line segments from $a_-$ to $-1$ and from $a_+$ to $1$.
$w$ is analytic on the complement of these two segments.
The contours in \eqref{eq:sfsfsdfsew} can be deformed onto the sides
of the cuts facing the origin, oriented from $a_-$ to $-1$ and from
$a_+$ to $1$. These orientations agree with clockwise traversal
around each cut. See \Cref{fig:contourdefo}.

Let $\gamma_1,\gamma_2$ be clockwise contours around the respective
cuts. The values of $w$ on the two sides of each cut have opposite
signs. The direction of integration also reverses, so the two sides
contribute equally, thus
\begin{align}
 \int_{a_-}^{-1}\frac{\dd u}{w(u)}
 &=\frac12\int_{\gamma_1}\frac{\dd u}{w(u)},&
 \int_{a_+}^{1}\frac{\dd u}{w(u)}
 &=\frac12\int_{\gamma_2}\frac{\dd u}{w(u)}.
\end{align}

Let $C_R$ be a counterclockwise circle of radius $R>0$ centered at the origin enclosing both cuts. Cauchy's
theorem on the region between $C_R$ and $\gamma_1,\gamma_2$ gives
\begin{align}
 \int_{\gamma_1}\frac{\dd u}{w(u)}
 +\int_{\gamma_2}\frac{\dd u}{w(u)}
 &=-\int_{C_R}\frac{\dd u}{w(u)}.
\end{align}
Since $1/w(u)=O(u^{-2})$ as $u\to\infty$, the right-hand side goes to zero as $R\to\infty$. Using \eqref{eq:sfsfsdfsew}, we conclude that
\begin{align}
 \pi\lambda\partial_E I
 &=\int_{a_-}^{-1}\frac{\dd u}{w(u)}
   +\int_{a_+}^{1}\frac{\dd u}{w(u)}=\frac12\left(
   \int_{\gamma_1}\frac{\dd u}{w(u)}
   +\int_{\gamma_2}\frac{\dd u}{w(u)}\right)=0.
\end{align}
Thus $I(E,\lambda)=I(0,\lambda)$ throughout the interior of
$\Sigma_\lambda$.

\begin{figure}[htbp]
 \centering
 \includegraphics[width=0.98\linewidth]{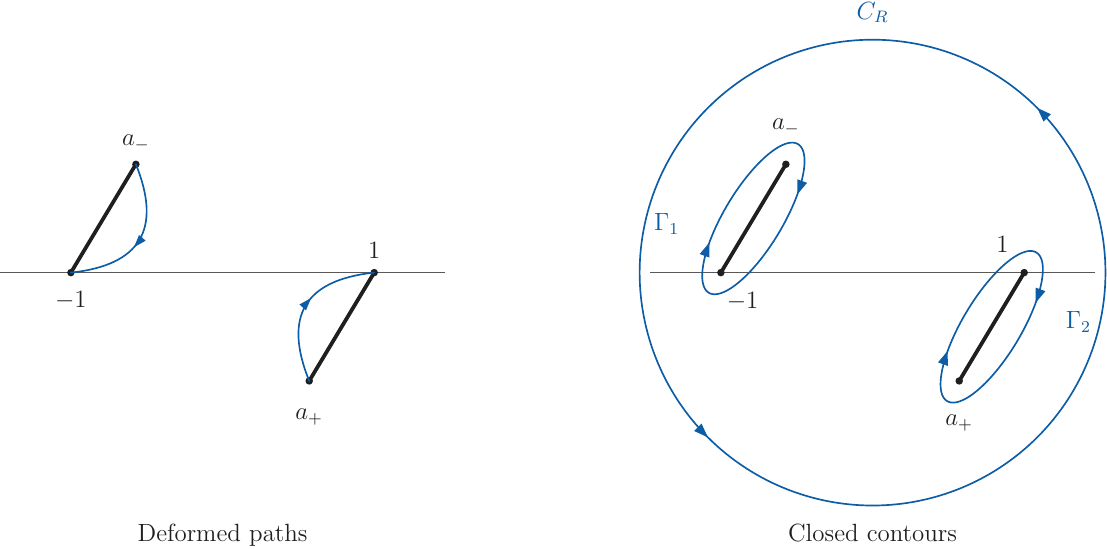}
 \caption{Deformed paths from $a_-$ to $-1$ and from $a_+$ to $1$
 (left) and the contour cancellation (right). The black segments are
 the two cuts. The inner contours are clockwise and the outer circle
 is counterclockwise.}
 \label{fig:contourdefo}
\end{figure}

\textbf{Evaluation at zero energy.}
At $E=0$, we have $p(z)=\arccos(-\lambda\cos z)$ and
$p(\pi/2)=\pi/2$. Split the contour at $\pi/2$:
\begin{align}\label{eq:half-period-zero-split}
 \widetilde\chi(0)
 &=\frac1\pi\int_{\zeta_+}^{\pi/2}p(z)\dd z
   +\frac1\pi\int_{\pi/2}^{\pi}p(z)\dd z.
\end{align}
For the first integral, integration by parts and
$z(p)=\arccos(-\lambda^{-1}\cos p)$ give
\begin{align}
 \int_{\zeta_+}^{\pi/2}p(z)\dd z
 &=\bigl[zp(z)\bigr]_{\zeta_+}^{\pi/2}
   -\int_0^{\pi/2}z(p)\dd p\\
 &=\frac{\pi^2}{4}
   -\int_0^{\pi/2}\arccos(-\lambda^{-1}\cos p)\dd p\\
 &=-\frac{\pi^2}{4}
   +\int_0^{\pi/2}\arccos(\lambda^{-1}\cos p)\dd p.
 \label{eq:half-period-zero-first-integral}
\end{align}
The last equality uses $\arccos(-v)=\pi-\arccos v$ on our branches.
In the second integral, the substitution $x=\pi-z$ gives
\begin{align}\label{eq:half-period-zero-second-integral}
 \int_{\pi/2}^{\pi}p(z)\dd z
 &=\int_0^{\pi/2}\arccos(\lambda\cos x)\dd x.
\end{align}
Substituting \eqref{eq:half-period-zero-first-integral} and
\eqref{eq:half-period-zero-second-integral} into
\eqref{eq:half-period-zero-split}, then setting $u=\cos x$, gives
\begin{align}
 \widetilde\chi(0)
 &=-\frac{\pi}{4}
   +\frac1\pi\int_0^1\frac{\arccos(\lambda u)+\arccos(\lambda^{-1}u)}
                    {\sqrt{1-u^2}}\dd u.
 \label{eq:half-period-tilde-at-zero}
\end{align}
At $E=0$, $\zeta_-=\pi-\zeta_+$ and $p(\pi-z)=\pi-p(z)$.
The substitution $x=\pi-z$ therefore also gives
\begin{align}
 \chi_c(0)
 &=\frac1\pi\int_0^{\zeta_-}(\pi-p(z))\dd z
   =\frac1\pi\int_{\zeta_+}^{\pi}p(x)\dd x=\widetilde\chi(0).
\end{align}
Hence
\begin{align}
 I(0,\lambda)
 &=-\frac{\pi}{2}
   +\frac2\pi\int_0^1\frac{\arccos(\lambda u)+\arccos(\lambda^{-1}u)}
                     {\sqrt{1-u^2}}\dd u.
 \label{eq:half-period-at-zero}
\end{align}

First take $\lambda=i\rho$, $\rho>0$. Since
$\arccos(ix)=\pi/2-i\operatorname{arsinh}x$ for real $x$,
\begin{align}
 I(0,i\rho)&=\frac{\pi}{2}-\frac{2i}{\pi}F(\rho),\\
 F(\rho)&\coloneq\int_0^1
 \frac{\operatorname{arsinh}(\rho u)-\operatorname{arsinh}(u/\rho)}
      {\sqrt{1-u^2}}\dd u.
\end{align}
Differentiating and then substituting $v=\sqrt{1-u^2}$, we obtain
\begin{align}
 F'(\rho)
 &=\int_0^1\left(
 \frac{u}{\sqrt{1+\rho^2u^2}}
 +\frac{u}{\rho^2\sqrt{1+u^2/\rho^2}}\right)
 \frac{\dd u}{\sqrt{1-u^2}}\\
 &=\frac1\rho\int_0^1\left(
 \frac1{\sqrt{1+\rho^{-2}-v^2}}
 +\frac1{\sqrt{1+\rho^2-v^2}}\right)\dd v\\
 &=\frac{\arctan\rho+\arctan(\rho^{-1})}{\rho}
 =\frac{\pi}{2\rho}.
\end{align}
Differentiation under the integral is justified because the differentiated
integrand is bounded by $C(1-u^2)^{-1/2}\in L^1(0,1)$, uniformly for
$\rho$ in compact subsets of $(0,\infty)$.
Since $F(1)=0$, it follows that
\begin{align}
 F(\rho)&=\frac\pi2\log\rho,&
 I(0,i\rho)&=\frac{\pi}{2}-i\log\rho=-i\log(i\rho).
\end{align}
Both sides of $I(0,\lambda)=-i\log\lambda$ are analytic for
$\Im\lambda>0$, so equality on the positive imaginary axis proves
it throughout the upper half-plane. Together with independence of $E$,
this gives
\begin{align}
 \widetilde\chi+\chi_c=I(E,\lambda)=I(0,\lambda)=-i\log\lambda.
\end{align}
The identity extends to the boundary of $\Sigma_\lambda$ by taking
the continuous limits of the action integrals. Uniqueness of analytic
continuation then preserves it whenever the two actions are continued
together.
\end{proof}

\subsection{The Scottish-flag computations}
\label{app:scottish-bulk}

Here we justify the properties of the action integrals for the
Scottish-flag matrix $H^N_{0,0,i}$ used in the quantization rules for the interior energies away from zero
\eqref{eq:scottish-classical-diagonal-quantization} and the explicit
formula \eqref{eq:scottish-interior-action-explicit}.
Throughout, $\lambda=i$. We take $N=4m$ in the final eigenvalue
lemma. The branches of $P,Z$ are those in
\eqref{eq:intro-classical-actions}.

\begin{lemma}
\label{lem:scottish-bulk-actions}
On the upper-right diagonal $E=t(1+i)$, $0<t<1$, we have
\begin{align}
 Z(E)&=\overline{P(E)},& \Im P(E)&<0.
\end{align}
On the upper-left diagonal $E=t(-1+i)$, $0<t<1$, we have
\begin{align}
 Z(E)&=\pi-\overline{P(E)},& \Im P(E)&<0.
\end{align}
Thus $P+Z$ is real on the upper-right diagonal and $Z-P$ is real on
the upper-left diagonal, while
\begin{align}
 \Im(P-Z)&=2\Im P<0 &&\text{on the upper-right diagonal},\\
 \Im(P+Z)&=2\Im P<0 &&\text{on the upper-left diagonal}.
\end{align}
\end{lemma}
\begin{proof} 
\textbf{Upper-right diagonal.} We have, for $t \in (0,1)$
\begin{align}
    Z(t(1+i)) &= \frac1{2\pi}\int_0^{2\pi} z(p)\dd p = \frac1{2\pi}\int_{0}^{2\pi} \arccos(-it+t+i \cos(p)) \dd p\\
    &= \frac1{2\pi}\int_{0}^{2\pi} \overline{\arccos(it+ t - i \cos (p))}\dd p = \overline{P(t(1+i))}.
\end{align}
To prove the sign of $\Im P$, put
$f_t(y)\coloneq\Im\arccos(t+iy)$. We then have
\begin{align}
 f_t'(y)
 &=-\Re\frac{1}{\sqrt{1-(t+iy)^2}}<0,
\end{align}
since $\Re(1-(t+iy)^2)=1-t^2+y^2>0$ and our square root
has positive real part. Therefore
\begin{align}
 \Im P(t(1+i))
 &=\frac1{2\pi}\int_0^{2\pi}f_t(t-\cos x)\dd x<\frac1{2\pi}\int_0^{2\pi}f_t(-\cos x)\dd x=0.
\end{align}
To see the last equality, set $x=u+\pi/2$ and note that $f_t$ is odd so that
\begin{align}
 \int_0^{2\pi}f_t(-\cos x)\dd x
 &=2\int_0^\pi f_t(-\cos x)\dd x=2\int_{-\pi/2}^{\pi/2}f_t(\sin u)\dd u=0.
\end{align}

\textbf{Upper-left diagonal.} For $E=t(-1+i)$, $0<t<1$, we have
\begin{align}
 Z(E)
 &=\frac1{2\pi}\int_0^{2\pi}\arccos(t+it+i\cos x)\dd x\\
 &=\frac1{2\pi}\int_0^{2\pi}\arccos(t+it-i\cos x)\dd x\\
 &=\frac1{2\pi}\int_0^{2\pi}
   \left[\pi-\overline{\arccos(-t+it-i\cos x)}\right]\dd x\\
 &=\pi-\overline{P(E)}.
\end{align}
We used translation by $\pi$ in the second step, then conjugation
and $\arccos(-w)=\pi-\arccos w$ in the third.
For the sign, the same identities give
\begin{align}
 P(t(-1+i))
 &=\frac1{2\pi}\int_0^{2\pi}
   \left[\pi-\overline{\arccos(t+it-i\cos x)}\right]\dd x\\
 &=\pi-\overline{P(t(1+i))}.
\end{align}
Consequently, the upper-right calculation implies
\begin{align}
 \Im P(t(-1+i))=\Im P(t(1+i))<0.
\end{align}
\end{proof}

We next express the action combinations $P+Z$ and $Z-P$ on the two diagonals through a single real integral.

\begin{lemma}
\label{lem:scottish-explicit-action}
For $0<t<\sqrt2$, the two diagonal actions satisfy
\begin{equation} \label{eq:scottish-radial-action}
    \begin{split}
     S(t)&\coloneq P(te^{i\pi/4})+Z(te^{i\pi/4})-\frac\pi2\\
 &=Z(te^{3i\pi/4})-P(te^{3i\pi/4})+\frac\pi2\\
 &=\frac2\pi\int_0^{\arccos(t^2/2)}
 \arccos\left(\frac{t}{\sqrt{2\cos u}}\right)\dd u.
    \end{split}
\end{equation}
The function $S$ extends continuously to $[0,\sqrt2]$, with
$S(0)=\pi/2$ and $S(\sqrt2)=0$, and satisfies $S'(t)<0$
for $0<t<\sqrt2$.
\end{lemma}
\begin{proof}
The conjugation identities proved in \Cref{lem:scottish-bulk-actions}
give the equality of the two action combinations.
Let $g(t)$ denote the final expression in
\eqref{eq:scottish-radial-action}. We show that $S=g$ by comparing
their derivatives and their values at $t=\sqrt2$.
By \Cref{lem:vanishing-action-two-periods},
$S(t)=V_+(te^{i\pi/4})$ near $t=\sqrt2$, and this identity continues
along the open diagonal. Thus $S(\sqrt2)=0$ by
\eqref{eq:upper-right-vanishing-action-normalization}.
For $\lambda=i$ and $E=te^{i\pi/4}$, the integral parameter in
\eqref{eq:vanishing-action-derivative} is $(2-t^2)/4\in(0,1/2)$,
so the continued square root remains positive.
Using $S'(t)=e^{i\pi/4}V_+'(te^{i\pi/4})$ gives
\begin{align}
 S'(t)&=-\frac2\pi\int_0^{\pi/2}
 \frac{\dd\vartheta}{\sqrt{1-\frac{2-t^2}{4}\sin^2\vartheta}}<0.
\end{align}
Differentiating $g$ and setting
$\sin(u/2)=\sqrt{(2-t^2)/4}\sin\vartheta$, we obtain
\begin{align}
 g'(t)&=-\frac2\pi\int_0^{\arccos(t^2/2)}\frac{\dd u}{\sqrt{2\cos u-t^2}}\\
 &=-\frac2\pi\int_0^{\pi/2}
 \frac{\dd\vartheta}{\sqrt{1-\frac{2-t^2}{4}\sin^2\vartheta}}=S'(t).
\end{align}
The integrand defining $g$ vanishes at the upper limit.
The function $g$ extends continuously to both endpoints, with
$g(\sqrt2)=0$ and $g(0)=(2/\pi)(\pi/2)(\pi/2)=\pi/2$.
Since $S'=g'$ and $S(\sqrt2)=g(\sqrt2)=0$, we have $S=g$,
which proves the formula and the endpoint values.
\end{proof}

We lastly show how we can use the actions $P$ and $Z$ to estimate eigenvalues of the Scottish-flag matrix $H^N_{0,0,i}$ along the diagonal segments, away from the origin and corners.

\begin{lemma}
\label{lem:scottish-interior-eigenvalues}
Let $K$ be a compact subsegment of
$\{te^{i\pi/4}:0<t<\sqrt2\}$.
For $N=4m$, let $t_{n,N}\in(0,\sqrt2)$ denote the solutions of
\begin{align}
 S(t_{n,N})=\frac{2\pi}{N}\left(n+\frac12\right),
 \qquad n=0,\ldots,m-1,
\end{align}
for $S$ defined in \eqref{eq:scottish-radial-action}.
These solutions exist uniquely because $S$ decreases from $\pi/2$ to zero.
There exist $C>0$ and a neighborhood $U$ of $K$, independent of $N$,
such that, for all sufficiently large $N=4m$,
\begin{enumerate}
 \item For each $n$ with $e^{i\pi/4}t_{n,N}\in K$, the disc
 \[
   \left\{z:\left\lvert z-e^{i\pi/4}t_{n,N}\right\rvert <\frac{C}{N^2}\right\}
 \]
 contains exactly one simple eigenvalue of $H^N_{0,0,i}$.
 \item Every eigenvalue of $H^N_{0,0,i}$ in $U$ lies within $C/N^2$
 of a candidate $e^{i\pi/4}t_{n,N}$ on the same diagonal segment.
\end{enumerate}
The analogous statement holds on the other three diagonal segments.
\end{lemma}
\begin{proof}
\textbf{The quantization condition.}
Recall the quantization condition \eqref{eq:scottish-balanced-factorization}, which we rewrite here for the convenience of the reader,
\begin{align}
 \cos\left(\frac{N(P_N+Z_N)}2\right)
 \cos\left(\frac{N(P_N-Z_N)}2\right)&=1.
 \label{eq:scottish-local-quantization}
\end{align}
Choose a compact subsegment $K'$ of the same open diagonal segment
containing $K$ in its relative interior. By \Cref{lem:scottish-bulk-actions},
$\Im(P-Z)<0$ on $K'$. Compactness and continuity give a fixed complex
neighborhood of $K'$ on which $\Im(P-Z)\le-\delta$ for some $\delta>0$.
On a finite cover of $K'$ by neighborhoods from
\Cref{thm:main}\textup{(i)}, the corrected actions satisfy
$P=P_N+O(N^{-2})$ and $Z=Z_N+O(N^{-2})$, uniformly after shrinking
the neighborhoods. Consequently, uniformly on a smaller neighborhood
of $K'$,
\begin{align}\label{eq:FNExpon}
 \left\lvert \cos\left(\frac{N(P_N-Z_N)}2\right)\right\rvert
 &\ge c_0e^{\delta_0N}
\end{align}
for some $c_0,\delta_0>0$ and all large $N$.
The quantization condition is therefore equivalent to $F_N(E)=0$, where
\begin{align}
 F_N(E)&\coloneq\cos\left(\frac{N(P_N+Z_N)}2\right)
 - \left( \cos\left(\frac{N(P_N-Z_N)}2\right)\right)^{-1}.
\end{align}
These functions are defined on the individual neighborhoods of the
finite cover, and all estimates below are uniform.

\textbf{Every eigenvalue is near a candidate.}
At an eigenvalue, \eqref{eq:scottish-local-quantization} and \eqref{eq:FNExpon} imply
\begin{align}
 \cos\left(\frac{N(P_N+Z_N)}2\right)&=O(e^{-\delta_0N}).
\end{align}
Using $\lvert \cos(a+ib)\rvert ^2=\cos^2a+\sinh^2b$ (for $a,b \in \mathbb R$), we get that there exists $k\in\mathbb Z$ such that
\begin{align}
 P_N+Z_N
 &=\frac{2\pi}{N}\left(k+\frac12\right)
   +O(N^{-1}e^{-\delta_0N}). \label{eq:will1}
\end{align}
Using the corrected-action estimates \eqref{eq:intro-corrected-actions},
\eqref{eq:will1} becomes
\begin{align}
    P+Z
    &=\frac{2\pi}{N}\left(k+\frac12\right)
    +O(N^{-2}). \label{eq:will2}
\end{align}
We have therefore proved that if $E$ is an eigenvalue, then \eqref{eq:will2} holds.
We want to now find a nearby point where \eqref{eq:will2} holds exactly (i.e., without the $O(N^{-2})$ term).

By \Cref{lem:scottish-explicit-action},
\begin{align}
 (P+Z)(te^{i\pi/4})&=S(t)+\frac\pi2.
\end{align}
Since $S'(t)<0$, $P+Z$ is one-to-one on $K'$, and
$(P+Z)(K)$ lies in the interior of the real interval $(P+Z)(K')$.
We may therefore choose a fixed neighborhood $U$ of $K$ small enough
that, for all sufficiently large $N$, \eqref{eq:will2} places
$\frac{2\pi}{N}(k+\frac12)$ in $(P+Z)(K')$ whenever $E\in U$
is an eigenvalue. Thus there is a unique $E_\ast=t_\ast e^{i\pi/4}\in K'$
such that
\begin{align}
 (P+Z)(E_\ast)&=\frac{2\pi}{N}\left(k+\frac12\right).
\end{align}
Moreover,
\begin{align}
 (P+Z)'(te^{i\pi/4})&=e^{-i\pi/4}S'(t)\ne0.
\end{align}
Since $P+Z$ is injective on $K'$ and its derivative does not vanish
there, it is injective on a sufficiently small complex neighborhood
of $K'$. Choose a smaller convex neighborhood of the real interval
$(P+Z)(K')$ on which the analytic inverse has bounded derivative.
After shrinking $U$, \eqref{eq:will2} therefore gives
\begin{align}
 \lvert E-E_\ast\rvert
 &\le C\bigl\lvert(P+Z)(E)-(P+Z)(E_\ast)\bigr\rvert
 =O(N^{-2}).
\end{align}
Lastly,
\begin{align}
    S(t_\ast) = \frac{2\pi}{N} (k+\frac 12) - \frac\pi2= \frac{2\pi}{N }(k - m +\frac{1}{2})
\end{align} 
using that $N = 4m$.
Setting $n=k-m$ gives $t_\ast=t_{n,N}$. Define
$E_{n,N}\coloneq e^{i\pi/4}t_{n,N}$, so $E_\ast=E_{n,N}$.
This proves the second assertion of the lemma.

\textbf{One eigenvalue near each candidate.}
Put
\begin{align}
 F_N^{\mathrm{cl}}(E)
 &\coloneq\cos\left(\frac{N(P+Z)}2\right).
\end{align}
Each $E_{n,N}\in K$ is a simple zero of $F_N^{\mathrm{cl}}$.
The action derivative is bounded above and away from zero on $K'$,
so the candidate spacing is bounded below by a positive constant
times $N^{-1}$.
Hence, for fixed $C>0$ and all sufficiently large $N$,
$F_N^{\mathrm{cl}}$ has exactly one zero in
$\lvert E-E_{n,N}\rvert <C/N^2$. On the boundary of this disc, Taylor expansion gives
\begin{align}
 \lvert F_N^{\mathrm{cl}}(E)\rvert &\ge\frac{c_1C}{N},\\
 \lvert F_N(E)-F_N^{\mathrm{cl}}(E)\rvert 
 &\le\frac{C_1}{N}+C_2e^{-\delta_0N}.
\end{align}
The second estimate follows from the corrected action identity \eqref{eq:intro-corrected-actions}
and the lower bound on the other cosine. The constants are uniform
for $E_{n,N}\in K$. Choose $C$ large enough for this comparison and
for the second assertion, then take $N$ sufficiently large.
Rouch\'e's theorem gives exactly one zero of $F_N$ in each disc,
counted with multiplicity. The exact trace identity
\eqref{eq:regular-side-two-cosine} and the determinant identity
\eqref{eq:finite-characteristic-determinant} give
\begin{align}
 2^{N-2}\det(EI_N-H^N_{0,0,i})
 &=\cos\left(\frac{N(P_N-Z_N)}2\right)F_N(E).
\end{align}
The cosine never vanishes here, so this zero is exactly one eigenvalue
counted with algebraic multiplicity.

On the upper-left diagonal segment, the same proof uses $Z-P$ in place of $P+Z$
and the other conclusions of \Cref{lem:scottish-bulk-actions,lem:scottish-explicit-action}.
The opposite diagonal segments follow from the symmetry $E\mapsto-E$ in
\eqref{eq:mN-sign-symmetries}.
\end{proof}

\clearpage

\bibliographystyle{alpha}  
\bibliography{references} 

\end{document}